\documentclass[a4paper,12pt,twoside]{article}
\usepackage[pdftex,colorlinks,bookmarksnumbered,bookmarksopen,citecolor=red,urlcolor=red]{hyperref}
\usepackage{graphicx}\usepackage{a4wide}
\usepackage{amssymb,amsmath,amsthm}
\usepackage{siunitx}
\usepackage{subfigure,color,colordvi,graphicx,xcolor}
\usepackage[only,llbracket,rrbracket,interleave]{stmaryrd}
\usepackage[sort,compress]{cite}
\usepackage{enumitem}
\usepackage{sidecap, caption}
\usepackage{multirow}

\makeatletter
\newcommand*\wt[2][0.2ex]{%
        \begingroup
        \mathchoice{\wt@helper{#1}{#2}{\displaystyle}{\textfont}}
                   {\wt@helper{#1}{#2}{\textstyle}{\textfont}}
                   {\wt@helper{#1}{#2}{\scriptstyle}{\scriptfont}}
                   {\wt@helper{#1}{#2}{\scriptscriptstyle}{\scriptscriptfont}}%
        \endgroup
        #2%
}
\newcommand*\wt@helper[4]{%
        \def\currentfont{\the#41}%
        \def\currentskewchar{\char\the\skewchar\currentfont}%
        \setbox\tw@\hbox{\currentfont#2\currentskewchar}%
        \dimen@ii\wd\tw@
        \setbox\tw@\hbox{\currentfont#2{}\currentskewchar}%
        \advance\dimen@ii-\wd\tw@
        \rlap{\raisebox{-#1}{$\m@th#3\kern\dimen@ii\widetilde{\phantom{#2}}$}}%
}
\makeatother

\usepackage{booktabs,array,dcolumn}
\newcommand{\PreserveBackslash}[1]{\let\temp=\\#1\let\\=\temp}
\newcolumntype{C}[1]{>{\PreserveBackslash\centering}p{#1}}
\newcolumntype{R}[1]{>{\PreserveBackslash\raggedleft}p{#1}}
\newcolumntype{L}[1]{>{\PreserveBackslash\raggedright}p{#1}}

\newcommand{\bm}[1]{\text{\boldmath $#1$\unboldmath}}
\newcommand{\abs}[1]{\lvert#1\rvert}

\newcommand{\vect}[1]{\mathbf{#1}}
\newcommand{\mat}[1]{\mathbf{#1}}

\newcommand{\Div}{{\bm{\nabla}{\cdot}}}
\newcommand{\grad}{\bm{\nabla}}
\newcommand{\defo}{\bm{\nabla}^{\texttt{s}}}
\newcommand{\RR}{\mathbb{R}}

\newcommand{\Ga}[1][\;]{\Gamma_{\!\!#1}}

\newcommand{\eltwo}{\ensuremath{\mathcal{L}_2}}
\newcommand{\elinf}{\ensuremath{\mathcal{L}_\infty}}

\newcommand{\nsd}  {\ensuremath{\texttt{n}_{\texttt{sd}}}}
\newcommand{\msd}  {\ensuremath{\texttt{m}_{\texttt{sd}}}}
\newcommand{\nrr}  {\ensuremath{\texttt{n}_{\texttt{rr}}}}
\newcommand{\numel}{\ensuremath{\texttt{n}_{\texttt{el}}}}
\newcommand{\numfa}{\ensuremath{\texttt{n}_{\texttt{fa}}}}

\newcommand{\numdof}[1]{\ensuremath{\texttt{n}_{\texttt{dof}}^{#1}}}

\newcommand{\hu}{\hat{u}}
\newcommand{\hp}{\hat{p}}
\newcommand{\bu}{\bm{u}}

\newcommand{\bhu}{\widehat{\bu}}

\newcommand{\bs}{\bm{s}}

\newcommand{\bsigma}{\bm{\sigma}}

\newcommand{\bsigmaD}{\bm{\sigma}^d}
\newcommand{\bepsD}{\bm{\varepsilon}^d}
\newcommand{\bL}{\bm{L}}
\newcommand{\bg}{\bm{g}}
\newcommand{\bn}{\bm{n}}
\newcommand{\bt}{\bm{t}}

\newcommand{\bNv}{\mat{N}_{\texttt{V}}}
\newcommand{\bNvi}{\mat{N}_{\texttt{V},i}}
\newcommand{\bNvj}{\mat{N}_{\texttt{V},j}}
\newcommand{\bDv}{\mat{D}_{\texttt{V}}}

\newcommand{\Lcomp}[1]{\mathrm{L}_{#1}}
\newcommand{\ncomp}[1]{\mathrm{n}_{#1}}

\newcommand{\tangK}{\bm{t}_{\!k}}

\newcommand{\bB}{\bm{B}}
\newcommand{\bhB}{\bm{\hat{B}}}
\newcommand{\bhBi}{\vect{\hat{B}}_i}
\newcommand{\bR}{\mat{R}}

\newcommand{\tauP}{\tau^p}
\newcommand{\btau}{\bm{\tau}}
\newcommand{\btauD}{\btau^d}
\newcommand{\btauA}{\btau^a}
\newcommand{\btauLF}{\bm{\tau}^a_\texttt{LF}}

\newcommand{\btauHLL}{\bm{\tau}^a_\texttt{HLL}}

\newcommand{\fluxEl}{\mathrm{J}_e}

\newcommand{\bue}{\vect{u}_e}
\newcommand{\buV}{\vect{u}}

\newcommand{\pe}{\mathrm{p}_e}
\newcommand{\bpV}{\vect{p}}

\newcommand{\bLVe}{\vect{L}_e}
\newcommand{\bLV}{\vect{L}}
\newcommand{\bLVv}{\vect{L}_{\texttt{V}}}
\newcommand{\bLVve}{\vect{L}_{\texttt{V},e}}

\newcommand{\bhuj}{\vect{\hat{u}}_j}
\newcommand{\bhui}{\vect{\hat{u}}_i}
\newcommand{\bhuV}{\hat{\vect{u}}}

\newcommand{\hpj}{\hat{\mathrm{p}}_j}
\newcommand{\hpi}{\hat{\mathrm{p}}_i}
\newcommand{\bhpV}{\hat{\vect{p}}}

\newcommand{\bse}{\vect{s}_e}
\newcommand{\bnV}{\vect{n}}

\newcommand{\buDj}{\vect{u}_{D,j}}
\newcommand{\bni}{\vect{n}_i}
\newcommand{\bnj}{\vect{n}_j}
\newcommand{\bgi}{\vect{g}_i}
\newcommand{\bgj}{\vect{g}_j}

\newcommand{\btki}{\vect{t}_{k,i}}
\newcommand{\btkj}{\vect{t}_{k,j}}

\newcommand{\bT}{\mat{T}}
\newcommand{\TLL}{\bT_{LL}}
\newcommand{\TLhu}{\bT_{L\hu}}
\newcommand{\Tuu}{\bT_{uu}}
\newcommand{\Tuhu}{\bT_{u\hu}}
\newcommand{\Tuhp}{\bT_{u\hp}}
\newcommand{\Tpp}{\bT_{pp}}
\newcommand{\Tphu}{\bT_{p\hu}}
\newcommand{\Tphp}{\bT_{p\hp}}
\newcommand{\ThuL}{\bT_{\hu L}}
\newcommand{\Thuu}{\bT_{\hu u}}
\newcommand{\Thuhu}{\bT_{\hu\hu}}
\newcommand{\Thuhp}{\bT_{\hu\hp}}
\newcommand{\Thpp}{\bT_{\hp p}}
\newcommand{\Thphp}{\bT_{\hp\hp}}

\newcommand{\bK}{\mat{K}}
\newcommand{\Khuhu}{\bK_{\hu\hu}}
\newcommand{\Khuhp}{\bK_{\hu\hp}}
\newcommand{\Khphu}{\bK_{\hp\hu}}
\newcommand{\Khphp}{\bK_{\hp\hp}}
\newcommand{\Klhu}{\bK_{\lambda\hu}}
\newcommand{\Klhp}{\bK_{\lambda\hp}}
\newcommand{\fV}[1]{\mathrm{f}_{\!#1}}
\newcommand{\bfV}[1]{\vect{f}_{\!#1}}

\newcommand{\TUU}{\bT_{UU}}
\newcommand{\TULambda}{\bT_{U\Lambda}}
\newcommand{\TLambdaU}{\bT_{\Lambda U}}
\newcommand{\TLambdaLambda}{\bT_{\Lambda\Lambda}}

\newcommand{\bU}{\vect{U}}
\newcommand{\bLambda}{\bm{\Lambda}}

\newcommand{\Id}[1]{\mat{I}_{#1}}
\newcommand{\Jd}[1]{\mat{J}_{#1}}
\newcommand{\Insd}{\mat{I}_{\nsd}}

\newcommand{\jump}[1]{\llbracket #1\rrbracket}

\newcommand{\bigjump}[1]{\bigl\llbracket #1\bigr\rrbracket}

\newcommand{\Aset}{\mathcal{A}_e}
\newcommand{\Bset}{\mathcal{B}_e}
\newcommand{\Cset}{\mathcal{C}_e}
\newcommand{\Dset}{\mathcal{D}_e}
\newcommand{\Eset}{\mathcal{E}_e}
\newcommand{\Iset}{\mathcal{I}_e}

\newcommand{\Nset}{\mathcal{N}_e}
\newcommand{\Sset}{\mathcal{S}_e}

\newcommand{\volE}{|\Omega_e|}
\newcommand{\areaFi}{|\Ga[e,i]|}
\newcommand{\areaFj}{|\Ga[e,j]|}

\newenvironment{keywords}{\begin{quote}\emph{\textbf{Keywords:}}}{\end{quote}}
\newtheorem{remark}{Remark}

\newcommand{\etal}{\textit{et al.}\ }

\RequirePackage{algorithm}[0.1] 
\usepackage{algorithmic} 

\graphicspath{{figsPDF/}{figsPNG/}{figs/}}

\usepackage{hyphenat}

\begin{document}
\title{A hybrid pressure formulation of the face-centred finite volume method for viscous laminar incompressible flows}

\author{
\renewcommand{\thefootnote}{\arabic{footnote}}
			  Matteo Giacomini\footnotemark[1]\textsuperscript{ \ ,}\footnotemark[2] ,
			  Davide Cortellessa$^\dagger$\textsuperscript{,}\footnotemark[1] ,
			  Luan M. Vieira$^\dagger$\textsuperscript{,}\footnotemark[1]\textsuperscript{ \ ,}\footnotemark[2]\textsuperscript{ \ ,}\footnotemark[3] ,  \\
\renewcommand{\thefootnote}{\arabic{footnote}}
			  Ruben Sevilla\footnotemark[3] \  and
			  Antonio Huerta\footnotemark[1]\textsuperscript{ \ ,}\footnotemark[2]\textsuperscript{ \ ,}*
}

\date{\today}
\maketitle

\renewcommand{\thefootnote}{\arabic{footnote}}

\footnotetext[1]{Laboratori de C\`alcul Num\`eric (LaC\`aN), ETS de Ingenier\'ia de Caminos, Canales y Puertos, Universitat Polit\`ecnica de Catalunya, Barcelona, Spain.}
\footnotetext[2]{Centre Internacional de M\`etodes Num\`erics en Enginyeria (CIMNE), Barcelona, Spain.}
\footnotetext[3]{Zienkiewicz Centre for Computational Engineering, Faculty of Science and Engineering, Swansea University, Swansea, SA1 8EN, Wales, UK.
\vspace{5pt}\\
$^\dagger$ These authors are listed in alphabetical order and contributed equally to the work. \\
$^*$ Corresponding author: Antonio Huerta \textit{E-mail:} \texttt{antonio.huerta@upc.edu}
}

\begin{abstract}
This work presents a hybrid pressure face-centred finite volume (FCFV) solver to simulate steady-state incompressible Navier-Stokes flows.
The method leverages the robustness, in the incompressible limit, of the hybridisable discontinuous Galerkin paradigm for compressible and weakly compressible flows to derive the formulation of a novel, low-order face-based discretisation.
The incompressibility constraint is enforced in a weak sense, by introducing an inter-cell mass flux defined in terms of a new, hybrid variable, representing the pressure at the cell faces.
This results in a new hybridisation strategy where cell variables (velocity, pressure and deviatoric strain rate tensor) are expressed as a function of velocity and pressure at the barycentre of the cell faces.
The hybrid pressure formulation provides first-order convergence of all variables, including the stress, without the need for gradient reconstruction, thus being less sensitive to cell type, stretching, distortion, and skewness than traditional low-order finite volume solvers.
Numerical benchmarks of Navier-Stokes flows at low and moderate Reynolds numbers, in two and three dimensions, are presented to evaluate accuracy and robustness of the method. 
In particular, the hybrid pressure formulation outperforms the FCFV method when convective effects are relevant,  achieving accurate predictions on significantly coarser meshes.
\end{abstract}


\section{Introduction}
\label{sc:Intro}

Despite the increasing interest towards high-order methods for the simulation of incompressible flows~\cite{Wang2013,Abgrall-AR-17,Fidkowski-FPN-23}, low-order approaches such as finite volume (FV) methods~\cite{Leveque2013,Barth-BHO-17} still represent the majority of available computational fluid dynamics (CFD) technologies for industrial applications.

Most existing FV codes can be classified in two groups,  namely, vertex-centred (VCFV) and cell-centred (CCFV) schemes.  Whilst these approaches differ in the positioning of the degrees of freedom (respectively, at mesh nodes and at cell barycentres) and in the definition of the control volume (respectively, using a dual mesh and exploiting existing mesh cells),  they share common drawbacks and limitations. In particular, both techniques rely on a reconstruction of the gradient of the unknown velocity to achieve first-order convergence of the stress, and such a reconstruction is prone to experience a significant loss of accuracy in the presence of distorted or stretched cells~\cite{diskin2010comparison,diskin2011comparison}.

To circumvent this issue, a new FV paradigm~\cite{RS-SGH:2018_FCFV1} was proposed. The face-centred finite volume (FCFV) method relies on a mixed formulation of the governing equations in each cell (which acts as a control volume), localising the degrees of freedom at the barycentre of the cell faces. The resulting method achieves first-order convergence of all the physical variables involved in the formulation, including the stress tensor, without the need to perform gradient reconstruction, thus being less sensitive to cell stretching and distortion with respect to traditional FV approaches.

Besides Poisson and Stokes equations~\cite{RS-SGH:2018_FCFV1}, FCFV solvers have been successfully developed for linear elasticity~\cite{RS-SGH:2019_FCFV2},  high-contrast and variable viscosity Stokes flows~\cite{sevilla2023face,sevilla2024face}, viscous laminar and inviscid compressible flows~\cite{vila2022non,vila2023benchmarking}, and viscous laminar and turbulent incompressible flows~\cite{Vieira-VGSH-24,Cortellessa-CGH-25}. Second-order versions of the FCFV method enriched by automatic mesh adaptation strategies were also presented for linear elliptic problems~\cite{RS-VGSH:20,MG-RS-20}.

It is worth recalling that second-order CCFV schemes can be seen as the lowest order version of a discontinuous Galerkin (DG) method, whereas second-order VCFV methods on simplicial meshes can be interpreted as the piecewise linear approximation of a continuous Galerkin (CG) formulation on a dual mesh. In a similar fashion, the FCFV paradigm can be seen as a hybridisable discontinuous Galerkin (HDG) method~\cite{cockburn2004characterization,Jay-CGL:09,Cockburn:16,HDGlab-GSH-20} with piecewise constant approximation in each mesh cell and on each mesh face.
Indeed, existing FCFV formulations for Stokes and incompressible Navier-Stokes equations were derived leveraging the knowledge on HDG for incompressible flows~\cite{Nguyen-NPC:10,Nguyen-NPC:11,Nguyen-NPC:11,GG-GFH:14,giacomini2018superconvergent}. 
Hence, the FCFV scheme inherits from HDG its properties, including the suitability to handle different cell types, the robustness to cell stretching and distortion, as well as the stability for perfectly incompressible flows, circumventing the Ladyzhenskaya-Babu\v{s}ka-Brezzi (LBB) condition~\cite{cockburn2009derivation},  without the need for pressure correction algorithms.
In addition,  FCFV mimics the hybridisation procedure of HDG to express the cell unknowns (i.e., velocity, pressure and gradient of velocity) as functions of the face velocities and the average pressure in the cell via the solution of a set of independent problems cell-by-cell.

Previous works on FCFV for incompressible flows highlighted that although the method achieves first-order convergence for velocity, pressure, and gradient of velocity, the approximation of pressure is prone to experience the largest values of the error~\cite{RS-SGH:2018_FCFV1,sevilla2023face,sevilla2024face,Vieira-VGSH-24}.
On the contrary, HDG methods devised for compressible and weakly compressible flows~\cite{AlS-SKGWH:20,JVP_HDG-VGSH:20,AlS-SF-22} were shown to be robust, accurate, and stable in the incompressible limit~\cite{Tutorial-GSH:2020}, also when using low-order polynomials.
Stemming from these observations,  this work presents a novel formulation of a face-based FV solver to simulate laminar incompressible Navier-Stokes flows. 
The method relies on a new hybridisation procedure such that (i) a new face variable is introduced to represent the trace of pressure on the cell faces,  (ii) mass conservation is enforced weakly (instead of strongly cell-by-cell), leading to an extra integral equation in each cell, and (iii) a suitable definition of the inter-cell mass flux in terms of the newly introduced \emph{hybrid pressure} is provided.

The three points mentioned above mark the difference between the hybrid pressure formulation and the FCFV method for laminar incompressible Navier-Stokes proposed by Vieira \etal~\cite{Vieira-VGSH-24}, where pressure is defined exclusively inside each cell and incompressibility is enforced strongly.
On the contrary, in order to introduce the new hybrid pressure,  mass conservation can only be imposed weakly. Nonetheless, it is worth noticing that upon mesh refinement the pressure defined in each cell converges to the hybrid pressure, thus guaranteeing that the flow is perfectly incompressible and mass conservation is verified up to machine precision.
The resulting hybrid pressure formulation thus outperforms the FCFV method, achieving superior accuracy for pressure and aerodynamic forces using coarser meshes, especially when convection phenomena are involved.
The proposed method is also closely related to the high-order HDG formulations that first introduced face unknowns for pressure~\cite{Rhebergen-RC-12,Wells-LW-12,Rhebergen-RW-18}. 
The novelties of the presented approach with respect to these works are the use of a mixed formulation with a symmetric mixed variable, the discretisation of all variables with the same degree of approximation (i.e., constant functions in the cells and on the faces), and the inclusion of Riemann solvers in the definition of the stabilisation strategy for convection.

The remainder of this article is organised as follows. 
Section~\ref{sc:NS} recalls the governing equations for a viscous laminar incompressible flow.
Section~\ref{sc:Foundation} introduces the new hybridisation procedure and derives the integral form of the hybrid pressure FCFV method for the incompressible Navier-Stokes equations.
The discrete problem and the corresponding nonlinear solver are described in Section~\ref{sc:Discrete}, whereas Section~\ref{sc:Computational} provides some considerations on the computational aspects and the numerical properties of the method.
Numerical examples, including academic test cases, as well as two- and three-dimensional benchmark problems for Navier-Stokes flows at low and moderate Reynolds numbers, are presented in Section~\ref{sc:Numerical}.
Finally, Section~\ref{sc:Conclusion} summarises the results of this work, whereas technical implementation aspects and the hybrid pressure formulation for the Stokes equations are presented in Appendix~\ref{app:Voigt} and~\ref{app:Stokes}, respectively.

\section{Governing equations for viscous laminar incompressible flows} 
\label{sc:NS}

Let $\Omega\subset\RR^{\nsd}$ be an open, bounded computational domain in $\nsd$ spatial dimensions, with boundary $\partial\Omega = \Ga[D] \cup \Ga[N] \cup \Ga[S]$ such that $\Ga[D]$, $\Ga[N]$, and $\Ga[S]$ are disjoint by pairs.  
Denoting by $(\bu,p)$ the velocity-pressure pair, the steady-state incompressible Navier-Stokes equations are given by
\begin{equation} \label{eq:NS}
 \left\{\begin{aligned}
 -\Div\bsigma + \Div(\bu \otimes \bu) &= \bs       &&\text{in $\Omega$,}\\
   \Div\bu &= 0  &&\text{in $\Omega$,} \\
   \bB(\bepsD,\bu,p) &= \bm{0} &&\text{on $\partial\Omega$,}
 \end{aligned}\right.
\end{equation}
where $\bs$ is the body force, $\bB$ denotes an operator enforcing the boundary conditions, and $\bsigma$ is the Cauchy stress tensor. For a Newtonian fluid under Stokes' hypothesis, the stress tensor is written as $\bsigma = -p \Insd + \bsigmaD$, where $\Insd$ is the identity matrix of dimension $\nsd \times \nsd$,  $\bsigmaD = \nu \bepsD$ is the deviatoric stress tensor and
\begin{equation}\label{eq:epsD}
\bepsD := 2 \defo\bu - \frac{2}{3}(\Div\bu)\Insd 
\end{equation}
is the deviatoric strain rate tensor,  $\defo :=(\grad+\grad^T)/2$ being the symmetric part of the gradient operator. 
Finally,  $\nu$ represents the kinematic viscosity of the fluid.  Under the assumption of unitary characteristic length and unitary characteristic velocity, the viscosity is set equal to the inverse of the Reynolds number $Re$, that is, $\nu = 1/Re$.
\begin{remark}\label{rmrk:relaxDiv}
Equation~\eqref{eq:epsD} is usually simplified by imposing that the velocity is pointwise divergence-free, thus neglecting the last term in the definition of $\bepsD$.  Nonetheless, the numerical formulation presented in this work enforces the incompressibility equation only in a weak sense and the term proportional to $\Div\bu$ is thus preserved in~\eqref{eq:epsD} and in the following derivation.
\end{remark}

The governing equations~\eqref{eq:NS} are finally closed by means of a set of boundary conditions, see equation~\eqref{eq:BC}, setting the value $\bu_D$ of the velocity on the Dirichlet boundary $\Ga[D]$, the traction $\bg$ on the Neumann portion $\Ga[N]$, whereas a perfectly slip and a non-penetration conditions are imposed on the symmetry contour $\Ga[S]$, that is,
\begin{align} \label{eq:BC}
\bB(\bepsD,\bu,p) & := 
\begin{cases}
	\bu - \bu_D
	& \text{ on  $\Ga[D]$,} \\
	-p\bn + \nu \bepsD\bn - \bg
	& \text{ on  $\Ga[N]$,} \\
	\left\{\begin{aligned} & \tangK \cdot \nu\bepsD \bn \\ & \bu \cdot \bn \end{aligned} \right\}
	& \text{ on  $\Ga[S]$, for $k=1,\ldots,\nsd-1$,}
\end{cases}
\end{align}
where the outward unit normal vector $\bn$ and the unit tangential vectors $\tangK,  k=1,\ldots,\nsd-1$ form an orthonormal system of vectors.
In this work, outflow boundaries are modeled using homogeneous Neumann conditions. For further details, Giacomini \etal~\cite{Tutorial-GSH:2020} present a numerical study of different conditions suitable for outlet boundaries.

Finally,  recall that the compatibility condition
\begin{equation} \label{eq:compatibilityCondition}
\int_{\Ga[D]} \bu_D \cdot \bn\, d\Gamma + \int_{\partial\Omega \setminus\Ga[D]} \bu \cdot \bn\, d\Gamma = 0 
\end{equation}
stems from the divergence-free condition on the velocity field in $\Omega$.
Moreover, if $\Ga[N] = \emptyset$, an additional constraint, e.g., 
\begin{equation} \label{eq:constraintDirichlet}
\int_{\Omega} p\, d\Omega = 0
\end{equation}
needs to be introduced to eliminate the indeterminacy of the pressure field~\cite{Donea-Huerta-2003}.

\section{Foundations of the FCFV method with hybrid pressure}
\label{sc:Foundation}

The domain $\Omega$ is subdivided into a set of $\numel$ non-overlapping cells $\Omega_e,  e=1,\ldots,\numel$.  The boundary $\partial \Omega_e$ of each cell consists of a set of $\numfa^e$ faces $\Ga[e,j], j=1,\ldots,\numfa^e$, with $\numfa^e$ varying depending on the type of cell $\Omega_e$. 
Moreover, the collection of all the faces excluding the ones on the boundary of the domain defines the internal interface (also known as \emph{internal mesh skeleton})
\begin{equation}\label{eq:Gamma}
 \Ga := \Bigg[ \bigcup_{e=1}^{\numel} \partial\Omega_e \Bigg]\setminus\partial\Omega .
\end{equation}

\begin{remark}
Note that previous works~\cite{RS-SGH:2018_FCFV1,MG-RS-20} showed the ability of the FCFV framework to seamlessly handle meshes of triangles, quadrilaterals, tetrahedra, hexahedra, prisms, and pyramids, as well as hybrid meshes consisting of different cell types.
For the sake of brevity, in this work, meshes of quadrilaterals and triangles will be employed in 2D, whereas only tetrahedral cells will be considered in 3D.
\end{remark}

Similarly to the FCFV method introduced by Sevilla \etal~\cite{RS-SGH:2018_FCFV1} and successively studied for different physical problems~\cite{RS-SGH:2019_FCFV2,vila2022non,vila2023benchmarking,sevilla2023face,sevilla2024face,Vieira-VGSH-24}, the FCFV formulation with hybrid pressure is constructed starting from a series of ingredients:
\begin{itemize}
\item a mixed, hybrid formulation of the equations;
\item a piecewise constant approximation of the unknowns in the cells and on the faces;
\item a definition of the numerical fluxes to handle inter-cell communication;
\item a hybridisation (or \emph{static condensation}) step to eliminate cell unknowns, expressing them in terms of face unknowns.
\end{itemize}

\begin{figure}[!htb]
	\centering
	\subfigure[]{\includegraphics[width=0.4\textwidth]{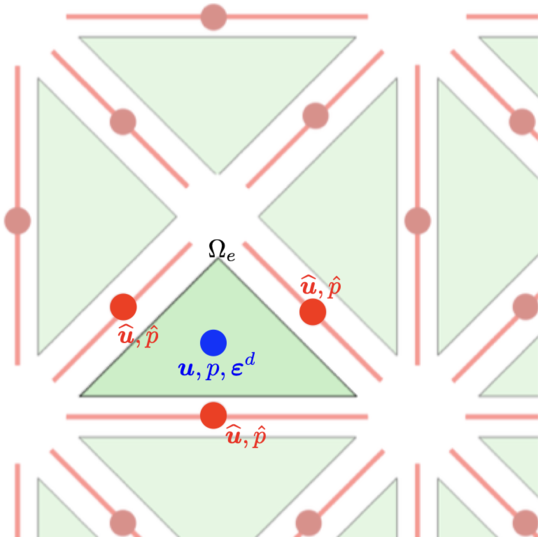}}
	\hspace{5pt}
	\subfigure[]{\includegraphics[width=0.4\textwidth]{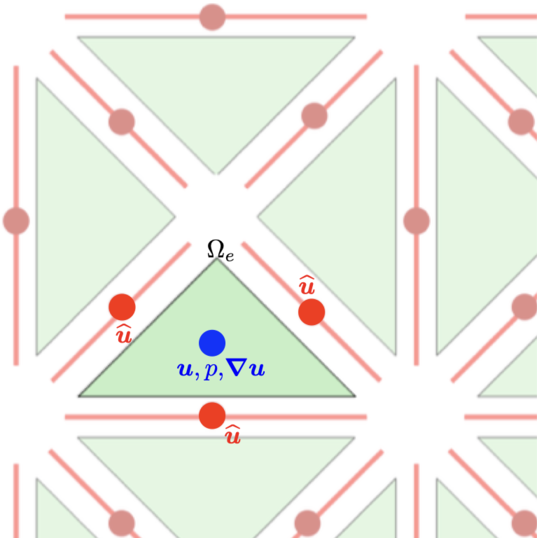}}
	
	\caption{Degrees of freedom of (a) hybrid pressure and (b) FCFV methods for cell $\Omega_e$. $\bu$: Cell velocity. $p$: Cell pressure. $\bepsD$: Cell deviatoric strain rate tensor.  $\grad\bu$: Cell velocity gradient tensor.  $\bhu$: Face velocity. $\hp$: Face pressure.}
	\label{fig:DOF}
\end{figure}

At the same time, several novelties are introduced by the hybrid pressure FCFV method.
First, whilst most of the above cited references rely on a velocity-pressure formulation of the Stokes and Navier-Stokes equations,  the present work derives the mixed problem starting from a Cauchy formulation of the momentum equation and imposes the symmetry of the mixed variable pointwise~\cite{HDG-Elasticity2018,giacomini2018superconvergent,Tutorial-GSH:2020}.
Second,  besides the piecewise constant approximation of velocity, pressure and deviatoric strain rate tensor in each cell, and the piecewise constant approximation of velocity on the cell faces, the proposed approach introduces a new, additional hybrid variable representing the trace of the pressure, also defined on the cell faces. 
A schematic representation of the degrees of freedom of the hybrid pressure formulation and FCFV method is displayed in Figure~\ref{fig:DOF}.
Moreover,  mass conservation equation is imposed in a weak sense, entailing the need to appropriately define an inter-cell mass flux,  which is expressed in terms of the newly introduced hybrid pressure variable.
Finally,  a new hybridisation strategy is proposed to express the cell unknowns as functions of the hybrid variables (i.e., face velocity and face pressure) via the solution of a set of independent problems cell-by-cell.
The following subsections detail the building blocks of the hybrid pressure formulation.

\subsection{Mixed hybrid formulation}
\label{sc:MixedHybrid}

Stemming from the rationale of HDG proposed by Cockburn and co-workers~\cite{cockburn2004characterization,Jay-CG:05,Jay-CGL:09,Nguyen-CNP:10,Nguyen-NPC:11}, the FCFV method with hybrid pressure employs a mixed formulation of the problem via the introduction of  $\bL=-\bepsD$ as mixed variable.
\begin{remark}
By construction, the mixed variable $\bL$ is symmetric. 
For convenience, Voigt notation is employed to enforce pointwise the symmetry at the discrete level,  storing only the non-redundant components of the tensor. 
An overview of the technical details is provided in Appendix~\ref{app:Voigt}, whereas a more extensive discussion on the topic was presented by Giacomini \etal~\cite{Tutorial-GSH:2020}.
\end{remark}

Upon introducing the new, hybrid variables $\bhu$ on $\Ga \cup \Ga[N] \cup \Ga[S]$ and $\hp$ on $\Ga \cup \partial\Omega$, respectively representing the trace of velocity and pressure at the cell faces, problem~\eqref{eq:NS} can be rewritten in each cell $\Omega_e$, for $e=1,\ldots, \numel$, as
\begin{equation}\label{eq:CellEq}
\left\{\begin{aligned}
  \bL + 2 \defo\bu - \dfrac{2}{3}(\Div\bu)\Insd       &= \bm{0}    &&\text{in $\Omega_e$,}\\	
  \Div(\nu \bL) + \Div (p \Insd + \bu\otimes\bu)                      &= \bs    &&\text{in $\Omega_e$, }\\
  \Div\bu                                                                                                        &= 0            &&\text{in $\Omega_e$, }\\
  \bu                                                                                                              &= \bu_D    &&\text{on $\partial\Omega_e \cap \Ga[D]$,}\\
  \bu               &= \bhu         &&\text{on $\partial\Omega_e \setminus \Ga[D]$,}\\
\end{aligned} \right.
\end{equation}
while imposing the boundary conditions on $\Ga[N] \cup \Ga[S]$ and the inter-cell continuity of momentum flux on the interior faces $\Ga$,
\begin{equation}\label{eq:FaceEq1}
\left\{\begin{aligned}
  \bB(-\bL,\bu,p) &= \bm{0} &&\text{on $\partial\Omega \setminus \Ga[D]$,} \\
  \bigjump{\bigl(\nu \bL + p \Insd + \bu\otimes\bu\bigr)\bn} &= \bm{0}  &&\text{on $\Gamma$,}
\end{aligned} \right.
\end{equation}
as well as the inter-cell continuity of velocity on $\Ga$, that is,
\begin{equation}\label{eq:FaceEq2}
\begin{aligned}
  \jump{\bu \otimes \bn}                                                                    &= \bm{0}   &&\text{on $\Gamma$.}\\  
\end{aligned}
\end{equation}
Note that, in~\eqref{eq:FaceEq1} and~\eqref{eq:FaceEq2},  the jump operator is defined as $\jump{\odot} := \odot_{\texttt{l}} + \odot_{\texttt{r}}$ by summing the contributions $\odot_{\texttt{l}}$ and $\odot_{\texttt{r}}$ of the cells $\Omega_{\texttt{l}}$ and $\Omega_{\texttt{r}}$ sharing the interface $\Ga$,  including the outward normal vectors $\bn_{\texttt{l}}$ and $\bn_{\texttt{r}}$, respectively~\cite{AdM-MFH:08}. It is straightforward to verify that this notation allows to rewrite~\eqref{eq:FaceEq2} on each face in $\Gamma$ as
\begin{equation}\label{eq:SolContinuityDetails}
\bm{0}= \jump{\bu \otimes \bn}  = \bu_{\texttt{l}} \otimes \bn_{\texttt{l}} + \bu_{\texttt{r}} \otimes \bn_{\texttt{r}} = (\bu_{\texttt{l}} - \bu_{\texttt{r}}) \otimes \bn_{\texttt{l}} ,
\end{equation}
where the last equality holds because $\bn_{\texttt{r}} = - \bn_{\texttt{l}}$, thus enforcing the expected inter-cell continuity of the velocity field on any arbitrarily oriented mesh face.

In order to have a unique solution, problem~\eqref{eq:CellEq} requires, for each cell, compatible velocities along its faces and an additional condition to handle the indeterminacy of pressure. Both restrictions are imposed weakly by means of~\eqref{eq:FaceEq1} and~\eqref{eq:FaceEq2}, as it will be detailed in the upcoming subsections. 
Consequently, for every cell $\Omega_e$, it is assumed that the hybrid velocity is compatible with the cell velocity via the boundary conditions in~\eqref{eq:CellEq},  and that the hybrid pressure, employed to weakly enforce cell-by-cell incompressibility,  is such that $\hp = p$ on $\partial\Omega_e$ removing the indeterminacy of the cell pressure. 

Moreover, it is worth noticing that equation~\eqref{eq:FaceEq2} can be rewritten in terms of the orthonormal system of vectors $\{\bn,\bt_{\!1}, \ldots , \bt_{\!\nsd-1} \}$ representing the normal and tangential directions to the internal mesh edges ($\nsd=2$) or faces ($\nsd=3$). More specifically,  the inter-cell continuity of the velocity is equivalent to 
\begin{equation}\label{eq:continuityNormTang}
\left\{\begin{aligned}
 \jump{\bu\cdot\bn} &=     0                    &&\text{on $\Gamma$,}\\
(\bu_{\texttt{l}}-\bu_{\texttt{r}})\cdot \tangK &=0    &&\text{on $\Gamma$, for $k = 1,\dotsc ,\nsd {-}1$.}
\end{aligned} \right.
\end{equation}
The first equation in~\eqref{eq:continuityNormTang} enforces the continuity of the normal component of the velocity of two neighbouring cells across the interface $\Gamma$ (i.e., mass conservation), whereas the remaining $\nsd-1$ equations enforce continuity of the tangential components of the velocity, that is, absence of shear effects on internal mesh faces. 

Given the uniqueness of $\bhu$ on the faces of the internal mesh skeleton and the condition imposed on $\bu$ at the cell boundaries in~\eqref{eq:CellEq}, the continuity of the tangential component of the velocity is automatically fulfilled.
On the contrary,  the continuity of the mass flux across $\Ga$, represented by the first equation in~\eqref{eq:continuityNormTang}, does not follow straightforwardly as in traditional formulations of hybrid methods~\cite{Tutorial-GSH:2020,RS-SGH:2018_FCFV1}. Indeed, the relaxation of the incompressiblity constraint (see Remark~\ref{rmrk:relaxDiv}) is responsible for the velocity to not fulfill strongly such equation and for the normal component of $\bu$ to violate the condition prescribed at the cell boundaries in problem~\eqref{eq:CellEq}.
Hence,  the hybrid pressure FCFV method needs to explicitly enforce the continuity of the mass flux on $\Ga$ and the appropriate value of the normal component of the velocity on $\partial\Omega$, that is, equation~\eqref{eq:FaceEq1} is to be complemented with
\begin{equation}\label{eq:FaceEq2revisited}
\left\{\begin{aligned}
  \bu \cdot \bn &= \bu_D \cdot \bn &&\text{on $\partial\Omega \cap \Ga[D]$,} \\
  \bu \cdot \bn &= \bhu \cdot \bn &&\text{on $\partial\Omega \setminus \Ga[D]$,} \\
  \jump{\bu\cdot\bn}                                                                           &= \bm{0}   &&\text{on $\Gamma$.}\\
\end{aligned} \right.
\end{equation}

Following the FCFV rationale, the solution strategy for the hybrid pressure formulation relies on two steps: first,  the so-called local equation~\eqref{eq:CellEq} is solved cell-by-cell to express the cell unknowns as a function of the face unknowns, i.e.\ $\bhu$ and $\hp$ (\emph{hybridisation} step); then,  the face variables are determined solving the so-called global problem to enforce equations~\eqref{eq:FaceEq1} and~\eqref{eq:FaceEq2revisited} (\emph{equilibration} step).

\subsection{Integral form of the local problem}
\label{sc:NSweakLoc}

The hybridisation step is obtained by applying the divergence theorem to~\eqref{eq:CellEq}, yielding the integral form of the local problem
\begin{equation}\label{eq:FCFV_Local}
\left\{\begin{aligned} 
 \int_{\Omega_e} \bL \, d\Omega  
 + \int_{\partial\Omega_e\cap\Ga[D]} \left(\bu_D{\otimes}\bn{+}\bn{\otimes}\bu_D {-} \frac{2}{3} (\bu_D{\cdot}\bn)\Insd \right) d\Gamma \hspace{70pt}& \\
 + \int_{\partial \Omega_e \setminus \Ga[D] } \left(\bhu   {\otimes}\bn{+}\bn{\otimes}\bhu {-} \frac{2}{3}(\bhu   {\cdot}\bn)\Insd \right) d\Gamma 
 &= \bm{0},
\\[1ex]
\int_{\partial \Omega_e} (\widehat{\nu\bL\bn}{-}\nu\bL\bn) \, d\Gamma
+ \int_{\partial \Omega_e} \hp\bn \, d\Gamma
+ \int_{\partial \Omega_e} (\widehat{\bu \otimes \bu})\bn \, d\Gamma 
&= \int_{\Omega_e}  \bs \, d\Omega ,
\\[1ex]
\int_{\partial \Omega_e} \widehat{\bu {\cdot} \bn} \, d\Gamma
&=0 ,
\end{aligned}\right.
\end{equation}
where the hybrid variable $\bhu$ and $\hp$ represent the traces of velocity and pressure on the cell faces. Moreover, $\widehat{\odot}$ denotes the trace, defined on the mesh skeleton, of the numerical flux of quantity $\odot$, and it will be detailed in Section~\ref{sc:Fluxes}.

Note that, in equation~\eqref{eq:FCFV_Local}, the divergence theorem is applied twice to the viscous flux in the momentum equation and once to the gradient of pressure and the convective fluxes in the momentum equation and to the mass flux in the continuity equation. This choice yields the appearance of the difference between the numerical and the physical fluxes when the divergence theorem is applied twice, allowing to significantly simplify the hybridisation step by decoupling the computation of the three variables involved in the local problem (see Section~\ref{sc:NR}).
A detailed discussion on the double application of the divergence theorem to the viscous part of the momentum equation was presented for the FCFV method by Sevilla and Duretz \cite{sevilla2024face}.

As mentioned in the previous section,  although problem~\eqref{eq:CellEq} features only Dirichlet conditions on the boundary of the cell, pressure is not underdetermined in $\Omega_e$ using this formulation, under the assumption that the face pressure $\hp$ matches the cell pressure $p$ on $\partial\Omega_e$. 
This requirement arises from the compatibility condition enforcing incompressibility at the cell level, that is, the third equation in~\eqref{eq:FCFV_Local}.
It is worth noticing that the weak imposition of mass conservation in such an equation needs to be complemented by an appropriate definition of the inter-cell mass flux $\widehat{\bu {\cdot} \bn}$, as detailed in Section~\ref{sc:Fluxes}.

Hence,  the solution of the local problem~\eqref{eq:FCFV_Local} provides the cell variables $(\bL,\bu,p)$ in terms of the velocity $\bhu$ and pressure $\hp$ defined on the mesh faces, which are then determined by solving the global problem described in the following subsection.

\subsection{Integral form of the global problem}
\label{sc:NSweakGlob}

The equilibration step determines the hybrid velocity-pressure pair $(\bhu,\hp)$ on the mesh faces by solving equations~\eqref{eq:FaceEq1} and~\eqref{eq:FaceEq2revisited}. More precisely, the global problem enforces the boundary conditions and the inter-cell continuity of momentum flux, as in traditional hybrid methods,
\begin{subequations}\label{eq:FCFV_Global}
\begin{equation}\label{eq:globMom}
\begin{aligned}
  \sum_{e=1}^{\numel} \Bigg\{ 
		\int_{\partial\Omega_e \setminus \partial\Omega} \widehat{\nu\bL\bn} \, d\Gamma
	+ \int_{\partial\Omega_e \setminus \partial\Omega} \hp\bn \, d\Gamma
	&+ \int_{\partial\Omega_e \setminus \partial\Omega} (\widehat{\bu \otimes \bu})\bn \, d\Gamma  \Bigg. \\
	&+ \Bigg. \int_{\partial\Omega_e \cap \partial\Omega\setminus\Ga[D]} \bhB(\bL,\bu,\bhu,\hp) \, d\Gamma   
\Bigg\} = \bm{0} ,
\end{aligned}
\end{equation}
together with the inter-cell continuity of mass flux and the normal component of the velocity on the external boundary, that is,
\begin{equation}\label{eq:globMass}
  \sum_{e=1}^{\numel} \Bigg\{ 
		\int_{\partial\Omega_e \setminus \partial\Omega} \widehat{\bu {\cdot} \bn} \, d\Gamma
	+ \int_{\partial\Omega_e \cap \Ga[D]} (\widehat{\bu {\cdot} \bn}{-}\bu_D{\cdot}\bn) d\Gamma
	+ \int_{\partial\Omega_e \cap \partial\Omega\setminus\Ga[D]} (\widehat{\bu {\cdot} \bn}{-}\bhu{\cdot}\bn) d\Gamma
\Bigg\} = 0 .
\end{equation}
\end{subequations}

The last term in equation~\eqref{eq:globMom} represents the hybrid version of the boundary operator $\bB$ defined in~\eqref{eq:BC}, and it will be detailed in the following subsection. 
Indeed, to fully characterise the local and global problems presented above, appropriate definitions of all traces of the numerical fluxes are to be introduced.

\subsection{Numerical fluxes and convection stabilisation inspired by Riemann solvers}
\label{sc:Fluxes}

The numerical traces of the viscous and convective fluxes~\cite{vila2022non,vila2023benchmarking,Vieira-VGSH-24} are defined as
\begin{subequations}\label{eq:NumFlux}
\begin{equation} \label{eq:traceDiffusion}
  \widehat{\nu\bL \bn} := 
  \begin{cases}
    \nu\bL\bn + \btauD (\bu - \bu_D) & \text{on $\partial\Omega_e\cap\Ga[D]$,} \\
    \nu\bL\bn + \btauD (\bu - \bhu) & \text{elsewhere,}  
  \end{cases}
\end{equation}
\begin{equation} \label{eq:traceAdvection}
  (\widehat{\bu\otimes\bu}) \bn := \begin{cases}
 (\bu_D \otimes \bu_D )\bn + \btauA(\bu-\bu_D) & \text{on $\partial\Omega_e\cap\Ga[D]$,} \\
 (\bhu  \otimes \bhu)\bn + \btauA(\bu-\bhu) & \text{elsewhere,}  
 \end{cases}
\end{equation}
with $\btauD$ and $\btauA$ being appropriately defined stabilisation coefficients. 
Moreover,  since the incompressibility equation is only enforced in a weak sense in the present formulation,  the numerical trace of the mass flux needs to be introduced
\begin{equation} \label{eq:traceMass}
  \widehat{\bu {\cdot} \bn} := \begin{cases}
 \bu_D {\cdot} \bn + \tauP( p-\hp) & \text{on $\partial\Omega_e\cap\Ga[D]$,} \\
 \bhu   {\cdot} \bn + \tauP( p-\hp) & \text{elsewhere,}  
 \end{cases}
\end{equation}
\end{subequations}
where $\tauP$ is a newly introduced stabilisation coefficient defined to control the compressibility effect due to the relaxation of the continuity equation. A numerical study of its influence will be presented in Section~\ref{sc:InfluenceTauStokes}.

It is worth recalling that the choice of the stabilisation coefficients $\btauD$ and $\btauA$ is known to have a significant role in guaranteeing the well-posedness and optimal convergence of the method, as extensively studied in the HDG~\cite{Jay-CGL:09,Cockburn-CDG:08,Nguyen-NPC:10,HDG-Elasticity2018,giacomini2018superconvergent,Tutorial-GSH:2020,JVP_HDG-VGSH:20} and FCFV~\cite{RS-SGH:2018_FCFV1,RS-SGH:2019_FCFV2,RS-VGSH:20,MG-RS-20,vila2022non} literature.
In this work,  the stabilisation of the viscous flux is defined as
\begin{equation} \label{eq:stabilisationDiffusion}
	\btau^d := \beta \nu \Insd ,
\end{equation}
where $\beta$ is a scalar parameter ranging between 5 and 10 in all the examples presented in Section~\ref{sc:Numerical}. 
The effect of the choice of this coefficient on the robustness of the FCFV method has been extensively investigated in the literature by means of numerical experiments~\cite{RS-SGH:2018_FCFV1,RS-VGSH:20,MG-RS-20,vila2022non,Vieira-VGSH-24}.

For the convective term,  FCFV stabilisations inspired by classic Riemann solvers have been derived for compressible~\cite{vila2022non} and incompressible~\cite{Vieira-VGSH-24} flows, building on the seminal work by Toro~\cite{Toro2009} and on the unified approach devised for HDG~\cite{JVP_HDG-VGSH:20}. 
In this work, Lax-Friedrichs (LF) and Harten-Lax-van Leer (HLL) convective stabilisations are considered, employing the definitions~\cite{Vieira-VGSH-24}
\begin{subequations}\label{eq:stabConv}
\begin{align}
\btauLF &:= \max\left\{2\abs{\bhu\cdot\bn},\xi \right\}\Insd ,
\label{eq:LFConv}
\\
\btauHLL &:= \max\{2(\bhu\cdot\bn),\xi\}\Insd ,
\label{eq:HLLConv}
\end{align}
\end{subequations}
where $\xi {=} 5 \times 10^{-2}$ is a scalar cut-off parameter. A detailed derivation of these expressions is available in~\cite{Vieira-VGSH-24}.

Similarly, the trace of the boundary operator $\bB$ needs to account for the conditions on the Neumann and symmetry boundaries, namely,
\begin{align} \label{eq:hybridBC}
\bhB(\bL,\bu,\bhu,\hp) & := 
\begin{cases}
	\hp\bn + \nu \bL\bn + \btau^d(\bu-\bhu) + \bg
	& \text{ on  $\Ga[N]$,} \\
	\left\{\begin{aligned} & \tangK \cdot \bigl[ \nu\bL \bn + \btau^d(\bu-\bhu) \bigr] \\ & \bhu \cdot \bn \end{aligned} \right\}
	& \text{ on  $\Ga[S]$, for $k=1,\ldots,\nsd-1$.}
\end{cases}
\end{align}

\subsection{Integral forms of the hybrid pressure FCFV method}
\label{sc:NSweak}

The integral forms of the FCFV method with hybrid pressure are obtained after inserting the numerical fluxes~\eqref{eq:NumFlux} into the local, see~\eqref{eq:FCFV_Local}, and global, see~\eqref{eq:FCFV_Global}, problems.
It follows that $(\bL,\bu,p)$ are first determined in each cell $\Omega_e$ as a function of $(\bhu,\hp)$ by solving the following set of $\numel$ independent local problems
\begin{subequations}\label{eq:FCFV_Solve}
\begin{equation}\label{eq:FCFV_Local_Solve}
\left\{\begin{aligned} 
 -\int_{\Omega_e} \bL \, d\Omega  
 =& 
 \int_{\partial\Omega_e\cap\Ga[D]} \left(\bu_D{\otimes}\bn{+}\bn{\otimes}\bu_D {-} \frac{2}{3} (\bu_D{\cdot}\bn)\Insd \right) d\Gamma \\
 &+ \int_{\partial \Omega_e \setminus \Ga[D] } \left(\bhu   {\otimes}\bn{+}\bn{\otimes}\bhu {-} \frac{2}{3}(\bhu   {\cdot}\bn)\Insd \right) d\Gamma ,
\\[1ex]
\int_{\partial \Omega_e} \btau \bu \, d\Gamma 
 =& 
 \int_{\Omega_e}  \bs \, d\Omega
 +\int_{\partial \Omega_e\cap\Ga[D]} (\btau - (\bu_D \cdot \bn)\Insd)\bu_D \, d\Gamma \\
 &+\int_{\partial \Omega_e\setminus\Ga[D]} (\btau - (\bhu \cdot \bn)\Insd)\bhu \, d\Gamma 
 -\int_{\partial \Omega_e} \hp\bn \, d\Gamma ,
\\[1ex]
-\int_{\partial \Omega_e} \tauP p \, d\Gamma
=&
 \int_{\partial \Omega_e\cap\Ga[D]} \bu_D \cdot \bn \, d\Gamma
 +\int_{\partial \Omega_e\setminus\Ga[D]} \bhu \cdot \bn \, d\Gamma
-\int_{\partial \Omega_e} \tauP \hp \, d\Gamma ,
\end{aligned}\right.
\end{equation}
where $\btau = \btauD + \btauA$. 
Similarly,  the global problem to determine the velocity $\bhu$ and pressure $\hp$ on the mesh faces is given by
\begin{equation}\label{eq:FCFV_Global_Solve}
\left\{\begin{aligned} 
  \sum_{e=1}^{\numel} \Bigg\{ 
	\int_{\partial\Omega_e \setminus \partial\Omega} \nu\bL\bn \, d\Gamma
	+\int_{\partial\Omega_e \setminus \partial\Omega} \btau \bu \, d\Gamma
	- \int_{\partial\Omega_e \setminus \partial\Omega} \btau \,\bhu \, d\Gamma \Bigg.  \hspace{75pt} &\\
	+ \Bigg. \int_{\partial\Omega_e \cap \partial\Omega\setminus\Ga[D]} \bhB(\bL,\bu,\bhu,\hp) \, d\Gamma   
\Bigg\} &= \bm{0} ,
\\[1ex]
  \sum_{e=1}^{\numel} \Bigg\{ 
		\int_{\partial\Omega_e} \tauP p \, d\Gamma
		-\int_{\partial\Omega_e} \tauP \hp \, d\Gamma
\Bigg\} &= 0 ,
\end{aligned}\right.
\end{equation}
\end{subequations}
where the cell-by-cell expressions of $\bL$, $\bu$, and $p$ are determined from equation~\eqref{eq:FCFV_Local_Solve}.
The corresponding local and global problems for the Stokes equations are reported in Appendix~\ref{app:Stokes}.

\begin{remark} \label{rmrk:FluxCancellation}
Two observations are important to derive the first equation of~\eqref{eq:FCFV_Global_Solve}.
As in standard HDG and FCFV methods, see, e.g., ~\cite{Tutorial-GSH:2020,Vieira-VGSH-24}, the uniqueness of $\bhu$ on the faces of the internal skeleton $\Ga$ yields the cancellation of the convective part of the physical flux, hence
\begin{subequations}\label{eq:fluxSimplif}
\begin{equation}\label{eq:fluxSimplifU}
    \sum_{e=1}^{\numel}\int_{\partial\Omega_e \setminus \partial\Omega} (\btau - (\bhu \cdot \bn)\Insd)\bhu \, d\Gamma 
 = \sum_{e=1}^{\numel}\int_{\partial\Omega_e \setminus \partial\Omega} \btau \,\bhu \, d\Gamma .
\end{equation}
Following a similar reasoning,  also the second term in equation~\eqref{eq:globMom} vanishes due to the uniqueness of the hybrid pressure on $\Ga$, that is,
\begin{equation}\label{eq:fluxSimplifP}
    \int_{\partial\Omega_e \setminus \partial\Omega} \hp\bn \, d\Gamma = 0.
\end{equation}
\end{subequations}
\end{remark}

\section{Hybrid pressure FCFV discretisation}
\label{sc:Discrete}

Let $\Aset := \{1, \dotsc, \numfa^e \}$ denote the set of indices for all faces of cell $\Omega_e$.  The sets $\Iset := \{j \in \Aset \; | \; \Ga[e,j] \cap \partial \Omega = \emptyset \}$ and $\Eset := \{j \in \Aset \; | \; \Ga[e,j] \cap \partial \Omega \neq \emptyset \}$ represent the set of indices for the cell faces interior to the domain and on the boundary, respectively.  In addition,  $\Dset := \{j \in \Aset \; | \; \Ga[e,j] \cap \Ga[D] \neq \emptyset \}$,  $\Nset := \{j \in \Aset \; | \; \Ga[e,j] \cap \Ga[N] \neq \emptyset \}$, and $\Sset := \{j \in \Aset \; | \; \Ga[e,j] \cap \Ga[S] \neq \emptyset \}$ denote the sets of indices for the cell faces on the Dirichlet, Neumann, and symmetry boundaries.
To simplify the presentation,  the sets $\Bset := \Aset \setminus \Dset$ and $\Cset := \Eset \setminus \Dset$ are also introduced to identify the sets of indices for all cell faces not on a Dirichlet boundary and for all cell faces on any portion of the external boundary, excluding the Dirichlet one. 
Moreover,  recall the definition of the indicator function of a generic set $\mathcal{Z}$ as
\begin{equation}
	\chi_{\mathcal{Z}}(i) = 
	\begin{cases}
		1 & \text{ if } \ i\in\mathcal{Z}, \\
		0 & \text{ otherwise}.                                                                 
	\end{cases}
\end{equation}

The discrete form of the FCFV method with hybrid pressure is obtained by describing the cell variables $(\bL,\bu,p)$ using a piecewise constant approximation in each cell and the hybrid variables $(\bhu,\hp)$ with a constant approximation on each cell face.  The corresponding values of the cell variables in $\Omega_e$ are denoted by $\bLVe$, $\bue$, and $\pe$, whereas the hybrid variables on the $j$-th face $\Ga[e,j]$ are denoted by $\bhuj$ and $\hpj$. 
In addition, $\buDj$ and $\bgj$ represent the discretisations of the Dirichlet and Neumann data $\bu_D$ and $\bg$, respectively,  at the barycentre of face $\Ga[e,j]$, whereas $\bse$ is the discretisation of the body force $\bs$ at the centroid of $\Omega_e$. Similarly,  $\bnj$ and $\btkj$ denote the normal vector $\bn$, and the tangent vector $\tangK$ at the barycentre of face $\Ga[e,j]$.
Finally, integrals in the cells and on the faces are performed using one quadrature point.

Upon discretisation, the residuals of the local problem for each cell $\Omega_e, \ e=1,\ldots,\numel$ are given by
\begin{subequations}\label{eq:discreteRes}
\begin{equation} \label{eq:localRes}
	\left\{
	\begin{aligned}
		\bR_{e,L} :=&   
		\volE \bLVe 
		+ \sum_{j \in \Dset}\! \areaFj \Bigl(\buDj \bnj^T  +  \bnj  \buDj^T - \frac{2}{3} (\bnj^T\buDj)\Insd \Bigr) \\
		&+ \sum_{j \in \Bset}\!  \areaFj \Bigl(\bhuj \bnj^T  +  \bnj  \bhuj^T - \frac{2}{3} (\bnj^T\bhuj)\Insd \Bigr)
		,	\\[1ex]
		\bR_{e,u} :=&  
		\sum_{j \in \Aset}\! \areaFj \btau_j \bue 
		- \volE \bse
		- \sum_{j \in \Dset}\!  \areaFj  \bigl( \btau_j - (\bnj^T\buDj)\Insd \bigr) \buDj \\
		&- \sum_{j \in \Bset}\!  \areaFj  \bigl( \btau_j - (\bnj^T\bhuj)\Insd \bigr) \bhuj
		+ \sum_{j \in \Aset}\! \areaFj \hpj \bnj
		,	\\[1ex]
		\bR_{e,p} :=&  
		\sum_{j \in \Aset}\! \areaFj \tauP_j \pe
		+ \sum_{j \in \Dset}\! \areaFj \bnj^T\buDj
		+ \sum_{j \in \Bset}\!  \areaFj \bnj^T\bhuj
		- \sum_{j \in \Aset}\! \areaFj \tauP_j \hpj
		,
	\end{aligned}
	\right.
\end{equation}
where $\volE$ is the area (in 2D) or volume (in 3D) of cell $\Omega_e$ and $\areaFj$ is the length (respectively, area) of the edge (respectively, face) $\Gamma_{e,j}$ in 2D (respectively, 3D). Moreover, the stabilisation coefficients $\btau_j$ and $\tauP_j$ are all evaluated at the barycentre of the corresponding faces.

The residuals of the global problem are to be computed on the mesh faces $\Ga[e,i]$, assembling the contributions of cell $\Omega_e$, namely, 
\begin{equation} \label{eq:globalRes}
	\left\{
	\begin{aligned}
		\bR_{e,i,\hu} &:=  
		\areaFi \Bigl\{ 
		\Bigl( \nu \bLVe \bni + \btau_i \bue - \btau_i \, \bhui \Bigr)\!  \chi_{\Iset}(i) 
		+ \bhBi\chi_{\Cset}(i) \Bigr\} ,
		&&\text{for all $i \in \Bset$}
		,	\\
		\bR_{e,i,\hp} &:= 
		\areaFi \tauP_i (\pe - \hpi) ,
		&&\text{for all $i \in \Aset$}
		,
	\end{aligned}
	\right.
\end{equation}
\end{subequations}
where the expression of the discrete boundary operator for the hybrid pressure formulation on face $\Ga[e,i]$ is given by $\bhBi := \bhBi(\bLVe,\bue,\bhui,\hpi,\btau_i^d)$, with
\begin{align} \label{eq:FCFV_BC}
\bhBi :=
\begin{cases}
	\hpi\bni + \nu \bLVe\bni + \btau_i^d(\bue-\bhui) + \bgi
	& \text{ for  $i \in \Nset$,} \\
	\left\{\begin{aligned} & \btki^T \bigl[ \nu\bLVe \bni + \btau_i^d(\bue-\bhui) \bigr] \\ & \bni^T\bhui \end{aligned} \right\}
	& \text{ for  $i \in \Sset$, and for $k=1,\ldots,\nsd-1$.}
\end{cases}
\end{align}

The resulting FCFV discrete problem with hybrid pressure involves the solution of the nonlinear system of equations
\begin{equation} \label{eq:nonLinearR}
	\bR(\bLV,\buV,\bpV,\bhuV,\bhpV) = \bm{0},
\end{equation}
obtained by assembling the contributions
\begin{equation} \label{eq:nonLinearRElem}
	\bR_{e,i} := 
	\begin{Bmatrix}	
		\bR_{e,L}(\bLVe,\bhuV)
		\\
		\bR_{e,u}(\bue,\bhuV,\bhpV)
		\\
		\bR_{e,p}(\pe,\bhuV,\bhpV)
		\\
		\bR_{e,i,\hu}(\bLVe,\bue,\bhuV,\bhpV) \chi_{\Bset}(i)
		\\
		\bR_{e,i,\hp}(\pe,\bhpV)
	\end{Bmatrix},
\end{equation}
for $e = 1,\dotsc,\numel$ and for $i \in \Aset$. 

\begin{remark} \label{rmrk:vectorisationVoigt}
From a practical viewpoint, the first equation in~\eqref{eq:localRes} is vectorised by means of Voigt notation to impose the symmetry of the mixed variable pointwise.  
Similarly,  Voigt notation is employed to rewrite the terms of the first equation in~\eqref{eq:globalRes} involving $\bLVe$ to efficiently compute the contraction of the deviatoric strain rate tensor with the normal and tangent vectors by means of standard matrix-vector products.
A description of the corresponding implementation details is reported in Appendix~\ref{app:Voigt}.
\end{remark}

\subsection{Newton-Raphson method for the incompressible Navier-Stokes equations}
\label{sc:NR}

To solve the nonlinear problem~\eqref{eq:nonLinearR}, the Newton-Raphson method is employed.  
It is worth noticing that the stabilisation tensor $\btauA$ depends upon the unknown hybrid velocity,  see~\eqref{eq:stabConv}. In order to simplify the linearisation of the residuals,  at each iteration $m$ of the nonlinear solver, the stabilisation tensors are evaluated at iteration $m-1$. 
While this noticeably simplifies the implementation of the algorithm, numerical experiments not reported here for brevity have shown that it does not affect the stability and convergence of the method.  
Indeed, since the residuals at iteration $m$ are evaluated using the stabilisation with the last computed approximation, the Jacobian matrix obtained by fixing $\btauA$ at iteration $m-1$ is consistent and quadratic convergence is observed in all cases.
This behaviour confirms the results previously observed for the closely related FCFV method for laminar incompressible Navier-Stokes  equations~\cite{Vieira-VGSH-24}.
It is worth noticing that the solution of a pseudo-transient problem, performed in the aforementioned reference to guarantee convergence of the Newton-Raphson algorithm for turbulent flows, is not required in the present work. 
Nonetheless, pseudo-time marching is expected to accelerate convergence also for the hybrid pressure formulation in the presence of stronger nonlinearities, for instance in the turbulent regime, but this lies beyond the scope of this contribution.

At each iteration $m$ of the Newton-Raphson algorithm,  the linear system to be solved is
\begin{equation} \label{eq:NRiteration}
	\begin{bmatrix}
		\TUU   & \TULambda  \\
		\TLambdaU   & \TLambdaLambda\\
	\end{bmatrix}^{m}
	\begin{Bmatrix}
		\Delta \bU \\
		\Delta \bLambda 
	\end{Bmatrix}^{\!m}
	= - 
	\begin{Bmatrix}
		\bR_U \\
		\bR_{\Lambda}
	\end{Bmatrix}^{\!m} ,
\end{equation}
where the residuals of the local and global problems respectively are
\begin{subequations}\label{eq:NRiterationDefs}
\begin{equation} \label{eq:NRiterationRes}
		\bR_U := 
		\begin{Bmatrix}
			\bR_L \\
			\bR_u \\
			\bR_p
		\end{Bmatrix},
		\qquad
		\bR_{\Lambda} := 
		\begin{Bmatrix}
			\bR_{\hu} \\
			\bR_{\hp} 
		\end{Bmatrix},
\end{equation}
the block matrices are defined as
\begin{equation}\label{eq:NRiterationMats}
\begin{gathered}
		\TUU := 
		\begin{bmatrix}
			\TLL   & \bm{0} & \bm{0}   \\
			\bm{0}   & \Tuu   & \bm{0}    		 \\
			\bm{0} & \bm{0} & \Tpp         	 
		\end{bmatrix},
		\qquad
		\TULambda := 
		\begin{bmatrix}
			\TLhu  & \bm{0}   \\
			\Tuhu  & \Tuhp  		  \\
			\Tphu & \Tphp
		\end{bmatrix},
		\\[1ex]
		\TLambdaU := 
		\begin{bmatrix}
			\ThuL  & \Thuu  & \bm{0}   \\
			\bm{0} & \bm{0} & \Thpp 		 
		\end{bmatrix},
		\qquad
		\TLambdaLambda := 
		\begin{bmatrix}
			\Thuhu & \Thuhp  \\
			\bm{0} & \Thphp 
		\end{bmatrix},
\end{gathered}
\end{equation}
with $\bT_{ab}$ being the tangent matrix obtained from the contribution $(\bT_{ab})_{e,i} := \partial \bR_{e,i,a}/\partial b$, whereas the solution increments from iteration $m$ to iteration $m+1$ are denoted by $\Delta \circledcirc^{m} = \circledcirc^{m+1} - \circledcirc^{m}$, with
\begin{equation} \label{eq:NRiterationInc}
		\bU := 
		\begin{Bmatrix}
			\bLV \\
			\buV \\
			\bpV 
		\end{Bmatrix},
		\qquad
		\bLambda := 
		\begin{Bmatrix}
			\bhuV \\
			\bhpV 
		\end{Bmatrix} .
\end{equation}
\end{subequations}

The structure of the resulting system~\eqref{eq:NRiteration} allows for the cell unknowns in $\bU$ to be eliminated by computing the Schur complement of the problem, thus leading to a smaller global system involving only the face unknowns $\bLambda$, namely, 
\begin{equation} \label{eq:NRiterationReduced}
	\mathbb{K}^{m} \Delta \bLambda^{m} = \mathbb{F}^{m},
\end{equation}
with $\mathbb{K} := \TLambdaLambda - \TLambdaU \TUU^{-1} \TULambda$ and $\mathbb{F} := -\bR_{\Lambda} + \TLambdaU \TUU^{-1}\bR_U$.
Note that the inverse of matrix $\TUU$ can be easily computed analytically. 
Indeed, stemming from the choice of applying the divergence theorem as discussed in Section~\ref{sc:NSweakLoc}, $\TUU$ features a block-diagonal structure.  

\begin{remark}
Contrary to the FCFV formulation in~\cite{RS-SGH:2018_FCFV1,Vieira-VGSH-24},  the matrix $\mathbb{K}$ in~\eqref{eq:NRiterationReduced} does not feature a saddle point structure.
This follows from the introduction of the hybrid pressure unknown in the formulation and, consequently, to the presence of a non-zero block $\Thphp$ in the matrix $\TLambdaLambda$.
Section~\ref{sc:Computational} further discusses the structure and the properties of the resulting matrix $\mathbb{K}$, whereas its detailed expression for the linear case (see equation~\eqref{eq:linearSysStokes}) and the derivation of the hybrid pressure formulation for the Stokes equations are reported in Appendix~\ref{app:Stokes}.
\end{remark}

Upon solving~\eqref{eq:NRiterationReduced}, the cell unknowns $\bU$ can be retrieved by computing, independently in each cell,
\begin{equation} \label{eq:NRiterationLocal}
	\TUU^{m} \Delta \bU^{m} =  -\bR_U^{m} - \TULambda^{m} \Delta \bLambda^{m}.
\end{equation}
It is worth noticing that, given the block-diagonal structure of $\TUU$, the local variables $\bLV$, $\buV$,  and $\bpV$ are decoupled from one another and the computation of equation~\eqref{eq:NRiterationLocal} is computationally inexpensive.

\begin{remark}
The computation of the reduced matrix $\mathbb{K}$ and vector $\mathbb{F}$ in equation~\eqref{eq:NRiterationReduced} corresponds to the hybridisation step described in Section~\ref{sc:NSweakLoc} and is the algebraic counterpart of the solution of the local problem~\eqref{eq:FCFV_Local_Solve}.
This strategy is equivalent to the well known static condensation and hybridisation procedures commonly performed for high-order finite element approximations~\cite{Guyan-65,Fraeijs-65,Cockburn:16,HDGlab-GSH-20}.
\end{remark}

\section{Computational aspects}
\label{sc:Computational}

The implementation of the hybrid pressure FCFV solver features three stages: (i) constructing matrix $\mathbb{K}$ and vector $\mathbb{F}$ at each iteration of the Newton-Raphson algorithm; (ii) solving equation~\eqref{eq:NRiterationReduced} to determine the face variables; (iii) retriving the cell variables by means of~\eqref{eq:NRiterationLocal}.
As discussed in Section~\ref{sc:NR}, the first and last steps are computationally inexpensive since (i) can be performed analytically and (iii) exploits the block-diagonal structure of matrix $\TUU$ to decouple the computation.
Hence, the main cost of the hybrid pressure FCFV method lies in the solution of the global problem~\eqref{eq:NRiterationReduced} in step (ii).
The simulations of the Navier-Stokes cases in Section~\ref{sc:Numerical} are obtained with a Fortran 90 implementation of the method, relying on PETSc~\cite{PETSC} for direct linear solvers.

In this context, the introduction of the face unknowns for pressure represents a major difference with respect to the FCFV method~\cite{RS-SGH:2018_FCFV1,Vieira-VGSH-24}, resulting in larger systems to be solved.
Table~\ref{tab:dofCompare} presents an estimate of the expected number of degrees of freedom of the hybrid pressure formulation with respect to the FCFV approach~\cite{RS-SGH:2018_FCFV1}, neglecting the boundary unknowns. On the one hand, the hybrid pressure formulation accounts for one degree of freedom for each component of the velocity and an additional one for pressure on each face yielding a total of $(\nsd+1)\numfa$ unknowns,  $\numfa$ being the number of mesh faces. On the other hand, the FCFV approach also features one face unknown per component of velocity, whereas pressure is discretised in each of the $\numel$ cells, leading to a system of dimension $\nsd\numfa+\numel$.
Hence, using simplicial meshes,  less than $15\%$ extra degrees of freedom are required to construct the hybrid pressure system with respect to the FCFV one.
\begin{table}[!htb]
\caption{Comparison of the expected number of degrees of freedom for the hybrid pressure formulation and the FCFV method using meshes of different cell types.}
\centering
\begin{tabular}{|l|c|c|c|c|c|}\hline
 \multirow{2}{*}{Cell type} & \multirow{2}{*}{Vertices}  &  \multirow{2}{*}{Cells} &  \multirow{2}{*}{Faces} & Hybrid pressure & FCFV \\
  & & & & $\numdof{}$ & $\numdof{}$ \\ 
 \hline
 Quadrilaterals  & $n$ & $n$      & $2n$     & $6n$      & $5n$ \\
 Triangles          & $n$ & $2n$    & $3n$     & $9n$     & $8n$ \\ 
 Tetrahedra       & $n$ & $5n$    & $10n$   & $40n$    & $35n$ \\ 
 Hexahedra       & $n$ & $n$      & $3n$     & $12n$     & $10n$ \\ 
 Prisms              & $n$ & $2n$    & $5n$     & $20n$    & $17n$ \\ 
 Pyramids         & $n$ & $8n/5$ & $4n$     & $16n$     & $68n/5$ \\ 
 \hline
\end{tabular}	
\label{tab:dofCompare}
\end{table}

\begin{figure}[!htb]
	\centering
	\subfigure[]{\includegraphics[width=0.4\textwidth]{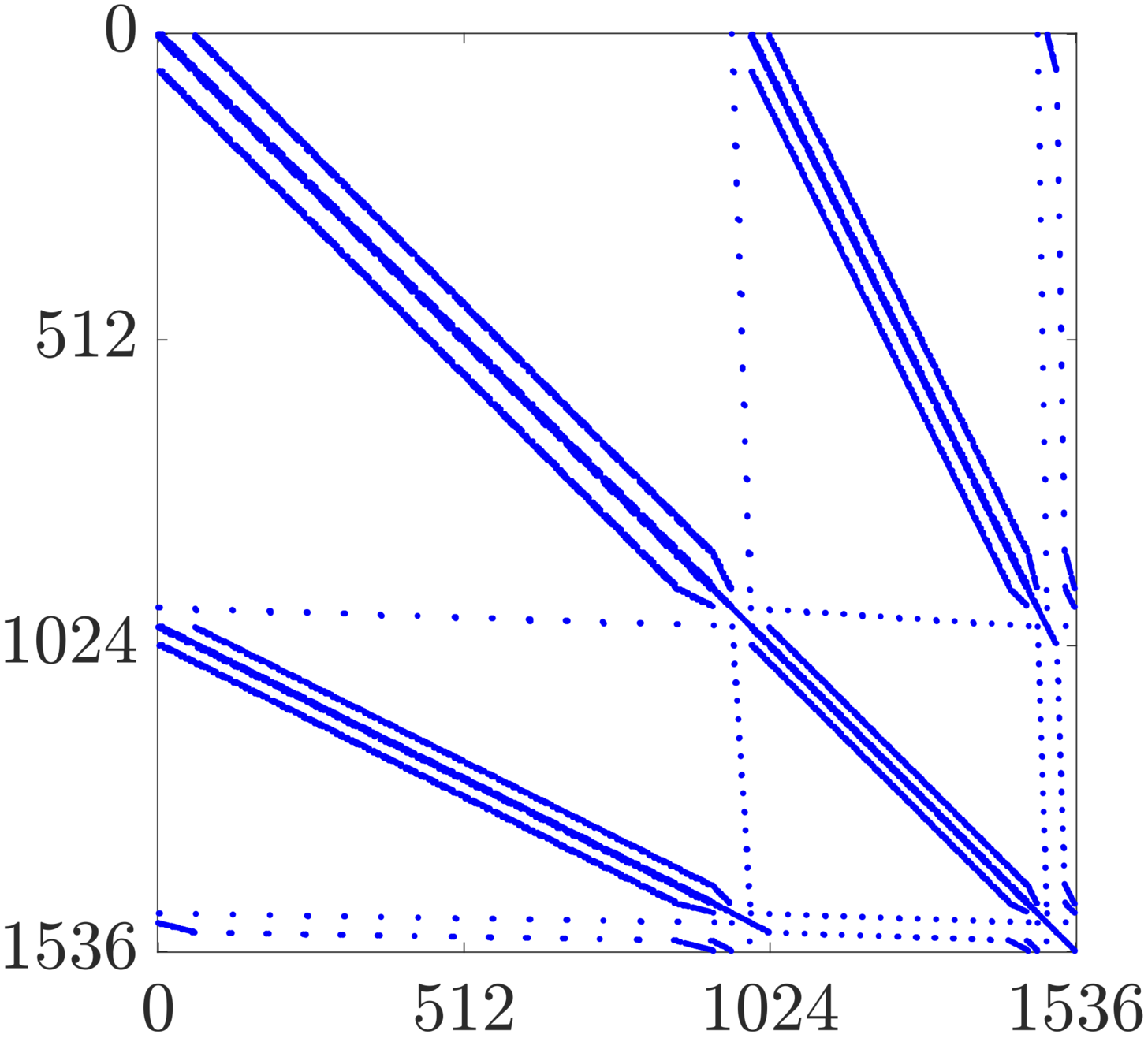}\label{fig:spyStokesHP}}
	\hspace{5pt}
	\subfigure[]{\includegraphics[width=0.4\textwidth]{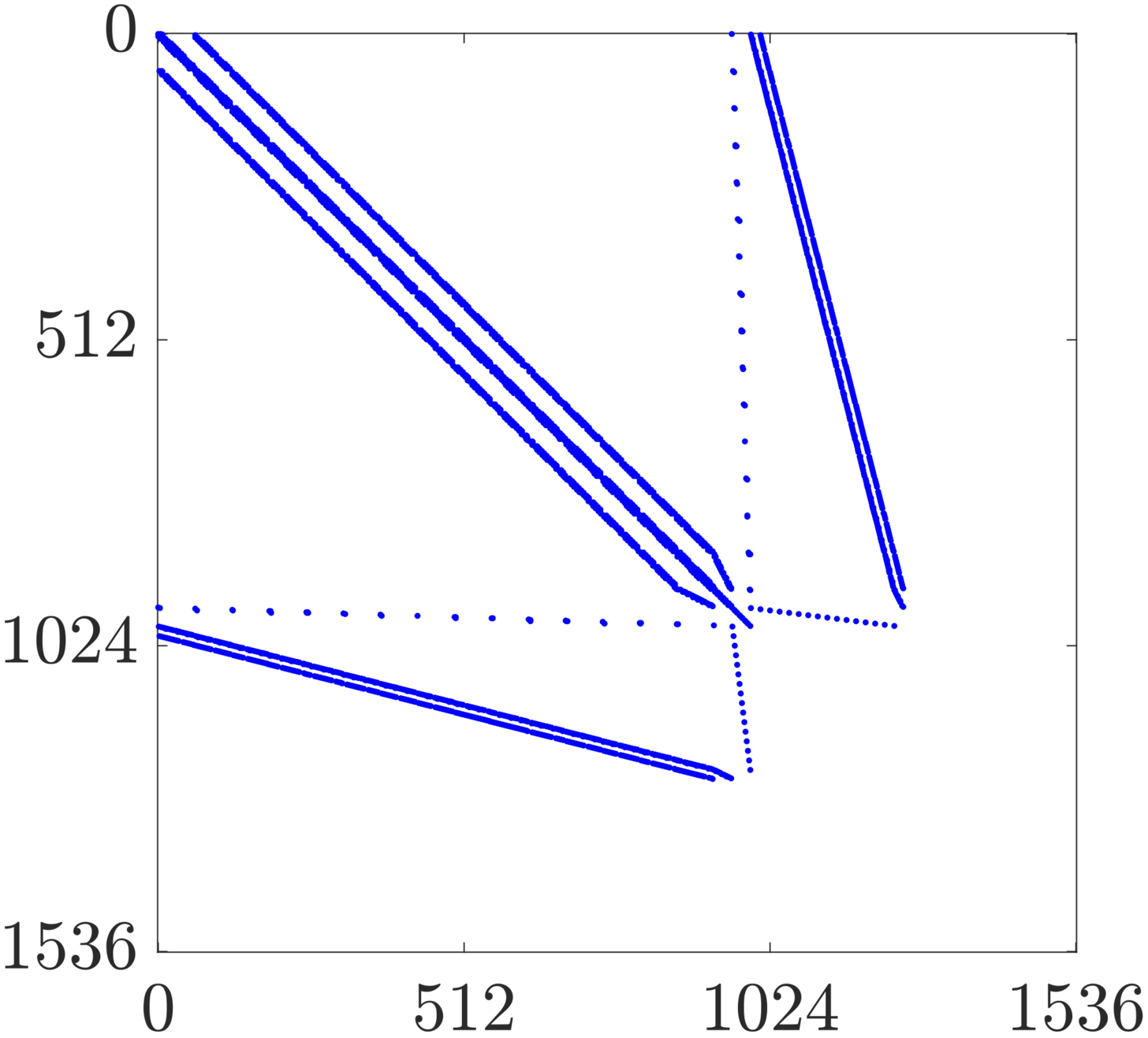}\label{fig:spyStokesFCFV}}
	
	\caption{Sparsity pattern of the matrix of the global system of a Stokes flow computed using (a) hybrid pressure formulation and (b) FCFV method. }
	\label{fig:spyStokes}
\end{figure}

The introduction of the face unknowns for pressure is also responsible for a significant change in the definition of the matrix of the global problem. Figure~\ref{fig:spyStokes} reports the sparsity pattern of the global system obtained using the hybrid pressure formulation and the FCFV method. 
For the sake of simplicity, a linear case is considered by neglecting the convection term in the Navier-Stokes equations. The resulting problem under analysis is a two-dimensional Stokes flow in a unit square domain, discretised by means of a structured mesh consisting of $16 \times 16$ quadrilateral cells.
The sparsity pattern is presented without any renumbering of the degrees of freedom in order to showcase the block structure of the problem. Of course, before solving the system, a reordering procedure should be applied to reduce the bandwidth of the matrix.

Several aspects of the structure of the hybrid pressure matrix are worth commenting.
Whilst the matrix blocks associated with the face velocities (top left) are comparable in the two approaches, both the off-diagonal terms and the pressure block present significant differences.

The off-diagonal blocks feature a more complex structure than the FCFV ones, as visible in Figure~\ref{fig:spyStokesHP}.
Indeed,  the dependence of both cell unknowns $\bu$ and $p$ on the face unknowns $\bhu$ and $\hp$, see~\eqref{eq:FCFV_Local_Solve}, is responsible for introducing a stronger coupling between velocity and pressure in the hybrid pressure formulation,  instead of $\bu$ depending exclusively on $\bhu$ (see equation~34b of the FCFV formulation~\cite{RS-SGH:2018_FCFV1}) and the cell pressure being used to enforce incompressibility cell-by-cell (see equation~37b of the FCFV formulation~\cite{RS-SGH:2018_FCFV1}).
Moreover,  the employment of a pointwise symmetric mixed variable (i.e., $\bL$ in~\eqref{eq:CellEq} instead of the scaled gradient of velocity in equation~25  of the FCFV formulation~\cite{RS-SGH:2018_FCFV1}) also introduces additional coupling between hybrid velocity and hybrid pressure via the inter-cell continuity of the momentum flux in equation~\eqref{eq:FCFV_Global_Solve}.

This additional richness of information is also testified by the larger number of non-zero entries in the resulting system for the hybrid pressure formulation: whilst the $1,536 \times 1,536$ matrix in Figure~\ref{fig:spyStokesHP} features $30,448$ non-zero entries, the FCFV matrix depicted in Figure~\ref{fig:spyStokesFCFV} has $10,728$ non-null terms over a total of $1,248 \times 1,248$ entries, whereas the number of degrees of freedom of pressure only increased from $256$ to $544$.

Moreover,  it is worth remarking that the hybrid pressure formulation no longer features the saddle point structure present in the FCFV matrix, as commented in Section~\ref{sc:NR}. This is a direct consequence of introducing the new variable $\hp$, that yields the matrix block $\Thphp$ in $\TLambdaLambda$~\eqref{eq:NRiterationMats} and the non-zero bottom-right block visible in Figure~\ref{fig:spyStokesHP}.

\begin{figure}[!htb]
	\centering
	\subfigure[]{\includegraphics[width=0.23\textwidth]{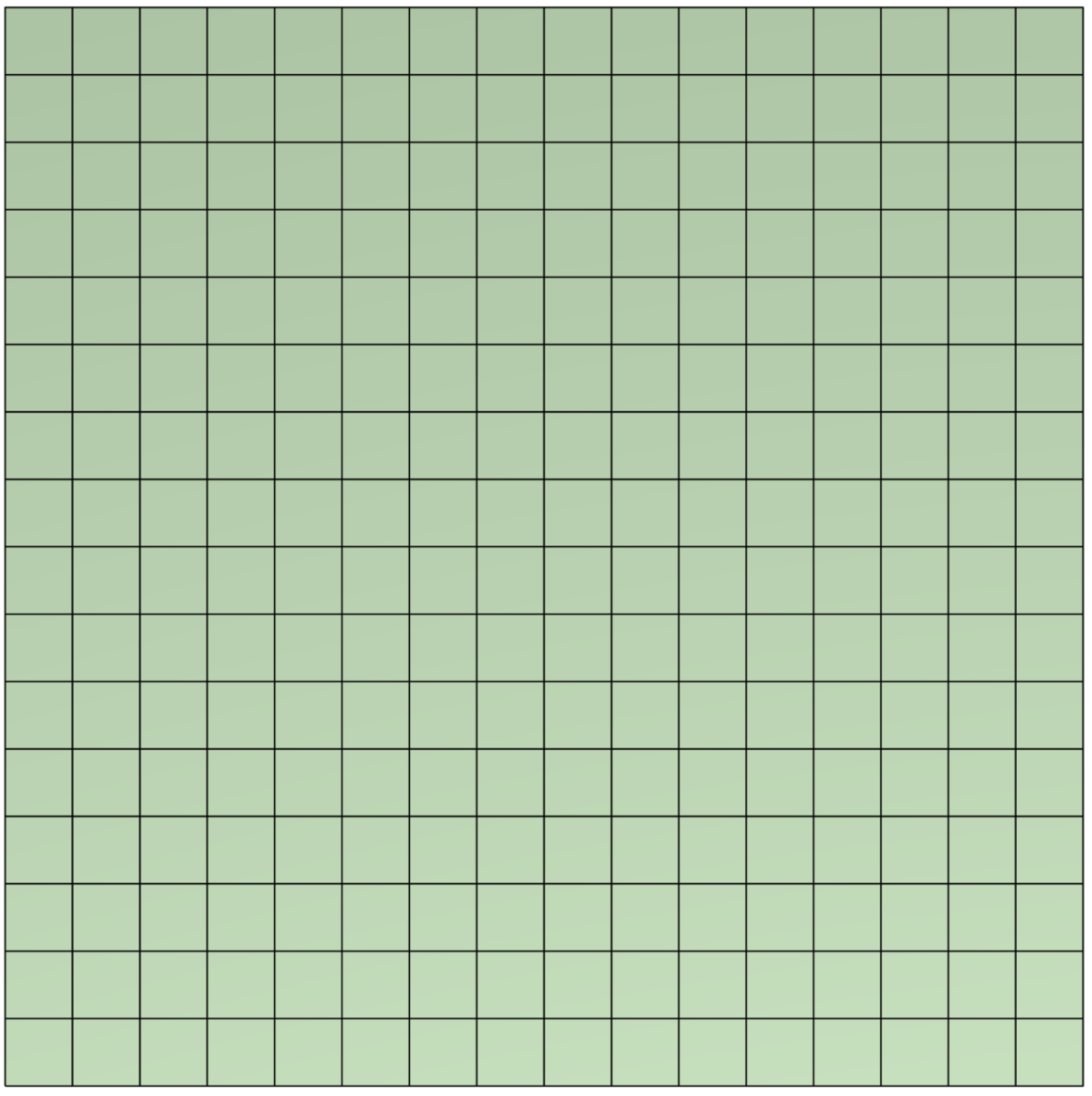}\label{fig:meshes2DStokesQUA}}
	\hspace{5pt}
	\subfigure[]{\includegraphics[width=0.23\textwidth]{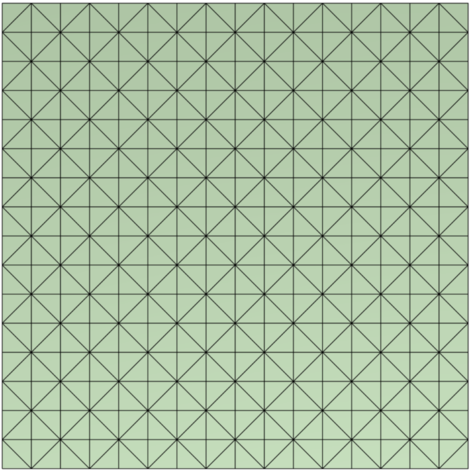}\label{fig:meshes2DStokesTRI}}
	\hspace{5pt}
	\subfigure[]{\includegraphics[width=0.23\textwidth]{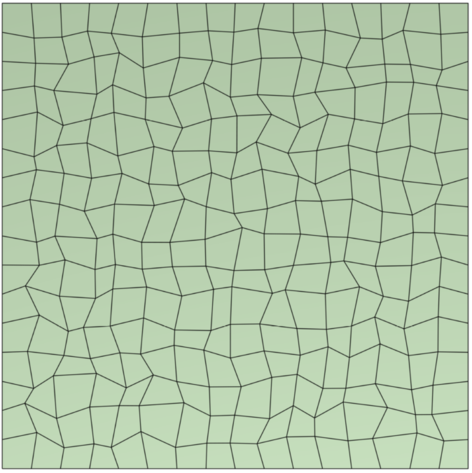}\label{fig:meshes2DStokesQUAd}}
	\hspace{5pt}	
	\subfigure[]{\includegraphics[width=0.23\textwidth]{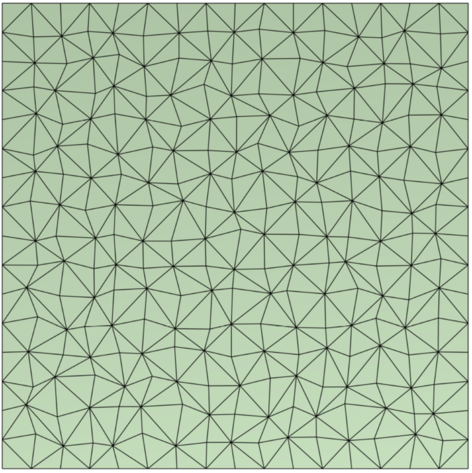}\label{fig:meshes2DStokesTRId}}
	\caption{Computational meshes of quadrilateral and triangular cells used to study the spectrum of the global matrix.  (a,b) Regular meshes. (c,d) Distorted meshes.}
	\label{fig:meshes2DStokes}
\end{figure}

A numerical study of the properties of the spectrum of the global matrix is now presented using the discretisation of a Stokes flow defined on a unit square domain.
Four meshes are considered, two featuring quadrilateral cells and two consisting of triangles: the quadrilateral meshes consist of $16 \times 16$ cells, whereas the meshes of triangular cells are composed of $16 \times 16 \times 2$ triangles.  Each set of meshes accounts for \emph{regular} meshes with uniform cell size (Figures~\ref{fig:meshes2DStokesQUA} and~\ref{fig:meshes2DStokesTRI}) and \emph{distorted} meshes in which the coordinates of the internal mesh nodes have been perturbed by means of a random value at most equal to the $30\%$ of the shortest edge emanating from such a node (Figures~\ref{fig:meshes2DStokesQUAd} and~\ref{fig:meshes2DStokesTRId}). 

Table~\ref{tab:eigenStokes2D} reports the minimum and maximum eigenvalues of the matrix of the global problem obtained from the hybrid pressure and FCFV discretisations.
\begin{table}[!htb]
\caption{Minimum and maximum eigenvalues of the global matrix for a Stokes flow using hybrid pressure formulation and FCFV method.}
\centering
\begin{tabular}{|l|c|c|c|c|}\hline
 \multirow{2}{*}{Mesh type} & $\mathcal{R}e(\lambda_{\min})$  & $\mathcal{R}e(\lambda_{\max})$ &  $\mathcal{I}m(\lambda_{\min})$ & $\mathcal{I}m(\lambda_{\max})$ \\
 \cline{2-5}
 & \multicolumn{4}{c|}{Hybrid pressure} \\
 \hline
 Regular quadrilaterals      & $-11.7$     & $-0.16 \times 10^{-3}$      & $-0.34 \times 10^{-5}$     & $0.34 \times 10^{-5}$ \\
 Distorted quadrilaterals   & $-13.9$    & $-0.16 \times 10^{-3}$      & $-0.11 \times 10^{-5}$     & $0.11 \times 10^{-5}$ \\
 Regular triangles              & $-18.9$     & $-0.17 \times 10^{-3}$      & $-0.65 \times 10^{-6}$     & $0.65 \times 10^{-6}$ \\
 Distorted triangles           & $-30.0$     & $-0.17 \times 10^{-3}$      & $-0.24 \times 10^{-4}$     & $0.24 \times 10^{-4}$ \\
 \hline
 \multicolumn{1}{c|}{} & \multicolumn{4}{c|}{FCFV} \\
 \hline
 Regular quadrilaterals      & $-7.7$     & $0.40 \times 10^{-2}$      & $0.00$     & $0.00$ \\
 Distorted quadrilaterals   & $-10.9$   & $0.49 \times 10^{-2}$      & $0.00$     & $0.00$ \\
 Regular triangles              & $-13.5$   & $0.32 \times 10^{-2}$      & $0.00$     & $0.00$ \\
 Distorted triangles           & $-22.7$   & $0.50 \times 10^{-2}$      & $0.00$     & $0.00$ \\
 \hline
\end{tabular}	
\label{tab:eigenStokes2D}
\end{table}

Contrary to the FCFV matrix whose eigenvalues are all real numbers,  the results show that the hybrid pressure matrix exhibits a certain number of complex eigenvalues.  It is worth noticing that,  in all analysed cases, these eigenvalues represent less than $7\%$ of the total number of eigenvalues. Moreover,  the magnitude of the corresponding imaginary parts always ranges between $10^{-6}$ and $10^{-4}$, thus significantly below the accuracy of the mesh under analysis,  the corresponding cell size being $0.62 \times 10^{-1}$.
In addition, whilst the FCFV matrix features a saddle point structure, thus having both positive and negative eigenvalues, the matrix arising from the hybrid pressure formulation only presents negative eigenvalues.  
Note also that the maximum eigenvalue for the hybrid pressure approach always remains bounded below zero, being almost insensitive to cell type and cell distortion.
Hence, the performed numerical studies suggest that the matrix obtained from the hybrid pressure discretisation is negative definite, but further theoretical analyses are required to confirm this conclusion.
Nonethless, these observations give promising insights on the possibility to develop efficient iterative schemes and preconditioners for the hybrid pressure FCFV method, which lie outside the scope of the present work.

\section{Numerical results}
\label{sc:Numerical}

Four examples are presented in this section to numerically evaluate the properties of the hybrid pressure formulation. 
First, in Section~\ref{sc:ConvergenceStokes}, it is shown that the method can be successfully employed to approximate the linear case of Stokes equations, achieving optimal convergence of order $1$ for velocity, pressure and deviatoric stress tensor, independently of cell type (quadrilaterals/triangles) and cell distortion. Moreover, the capability of the method to guarantee mass conservation although incompressibility is only enforced in a weak sense, is verified numerically. 
Then, the method is used to simulate a suite of steady-state incompressible Navier-Stokes flows.  The Couette flow (Section~\ref{sc:ConvergenceNS}) is employed to assess the convergence properties of the method. The lid-driven cavity flow is simulated for different Reynolds numbers in Section~\ref{sc:CavityNS}, demonstrating the capability of the method to accurately approximate convection-dominated flows, using different Riemann solvers. Finally, the three-dimensional flow past a sphere confirms the robustness and the accuracy of the approach, also with unstructured meshes of tetrahedra (Section~\ref{sc:SphereNS}).
Specific attention is devoted to the comparison of the results with reference solutions, when available in the literature, and with the previously published FCFV method~\cite{Vieira-VGSH-24},  showcasing the superior performance of the hybrid pressure formulation in the presence of convective effects.

\subsection{A synthetic Stokes flow}
\label{sc:ConvergenceStokes}

In this section, the optimal convergence properties of the method are numerically verified for the Stokes equations using a two-dimensional synthetic problem with viscosity $\nu=1$.
The computational domain is $\Omega=[0,1]^2$ and the boundary $\partial\Omega$ is composed of two disjoint portions,  $\Ga[N] := \{(x_1,x_2) \in \RR^2 \; | \; x_2=0\}$ and $\Ga[D] := \partial \Omega \setminus \Ga[N]$.
Starting from an analytical solution of the flow equations~\cite{Cheung-CCKQ-15},  velocity $\bu=(u_1,u_2)^T$ and pressure are defined as
\begin{equation}\label{eq:analSolS}
\bu(x_1,x_2) =
\begin{Bmatrix}
(1-\cos(2\pi x_1))\sin(2\pi x_2) \\
- \sin(2\pi x_1) (1-\cos(2\pi x_2))
\end{Bmatrix}
, \qquad
p(x_1,x_2) = \cos(\pi x_1)+\cos(\pi x_2) ,
\end{equation}
and they are used to determine the analytical expressions of body force $\bs$, Dirichlet datum $\bu_D$ and Neumann condition $\bg$.

Five structured meshes of quadrilateral cells are employed to estimate the convergence rate of the method. The meshes are uniformly refined  and the $i$-th mesh consists of $(8 \times 2^i)\times(8 \times 2^i)$ quadrilateral cells.
In addition, the effect of cell distortion on the accuracy of the computations is assessed, by displacing the position of each internal mesh node by a random perturbation of up to $30\%$ of the length of the shortest edge sharing the node. A detailed description of the strategy used to distort the meshes is available in~\cite{RS-SGH:2018_FCFV1}. The meshes previously displayed in Figure~\ref{fig:meshes2DStokesQUA} and~\ref{fig:meshes2DStokesQUAd} respectively represent the first level of refinement for regular and distorted meshes of quadrilateral cells.

\begin{figure}[!htb]
	\centering
	\subfigure[]{\includegraphics[width=0.32\textwidth]{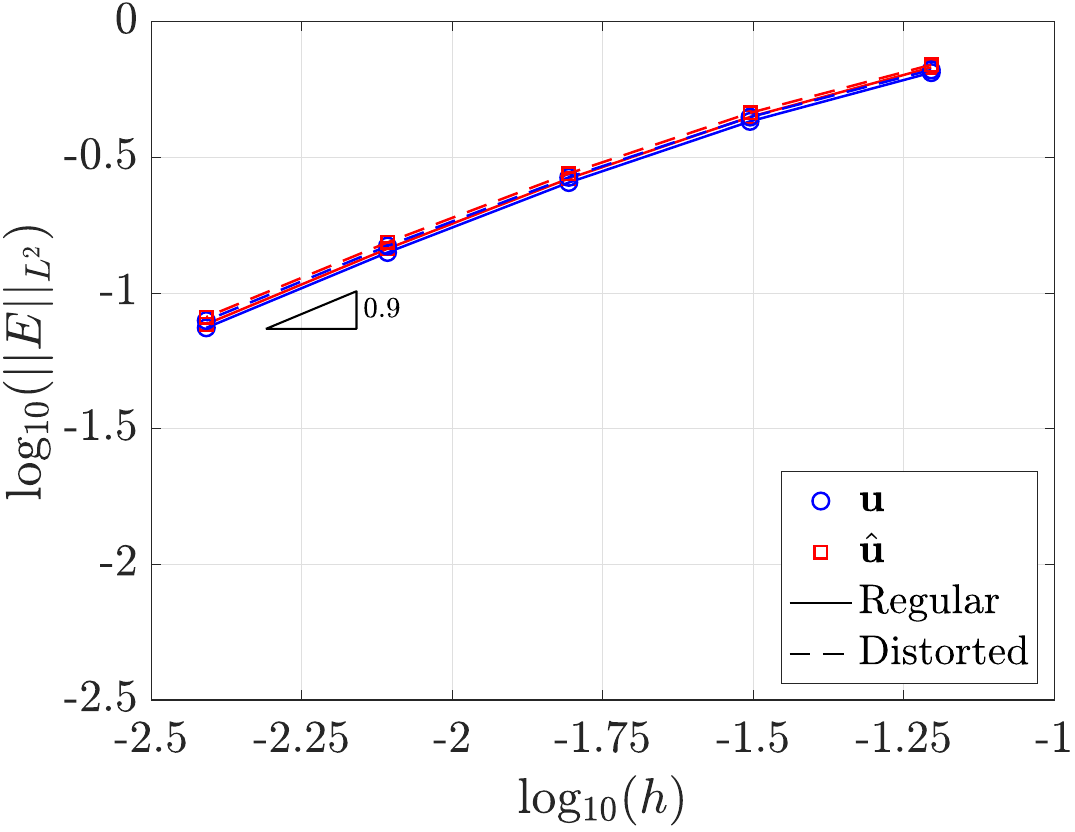}}
	\hspace{2pt}
	\subfigure[]{\includegraphics[width=0.32\textwidth]{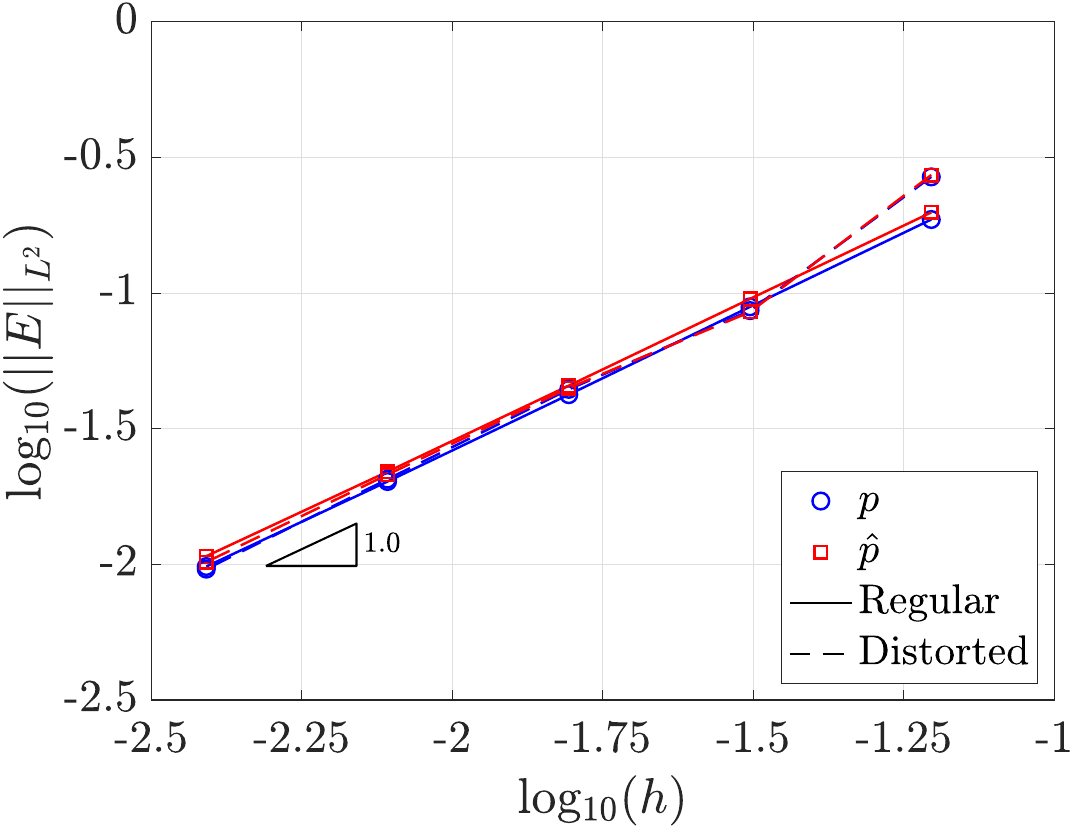}}
	\hspace{2pt}	
	\subfigure[]{\includegraphics[width=0.32\textwidth]{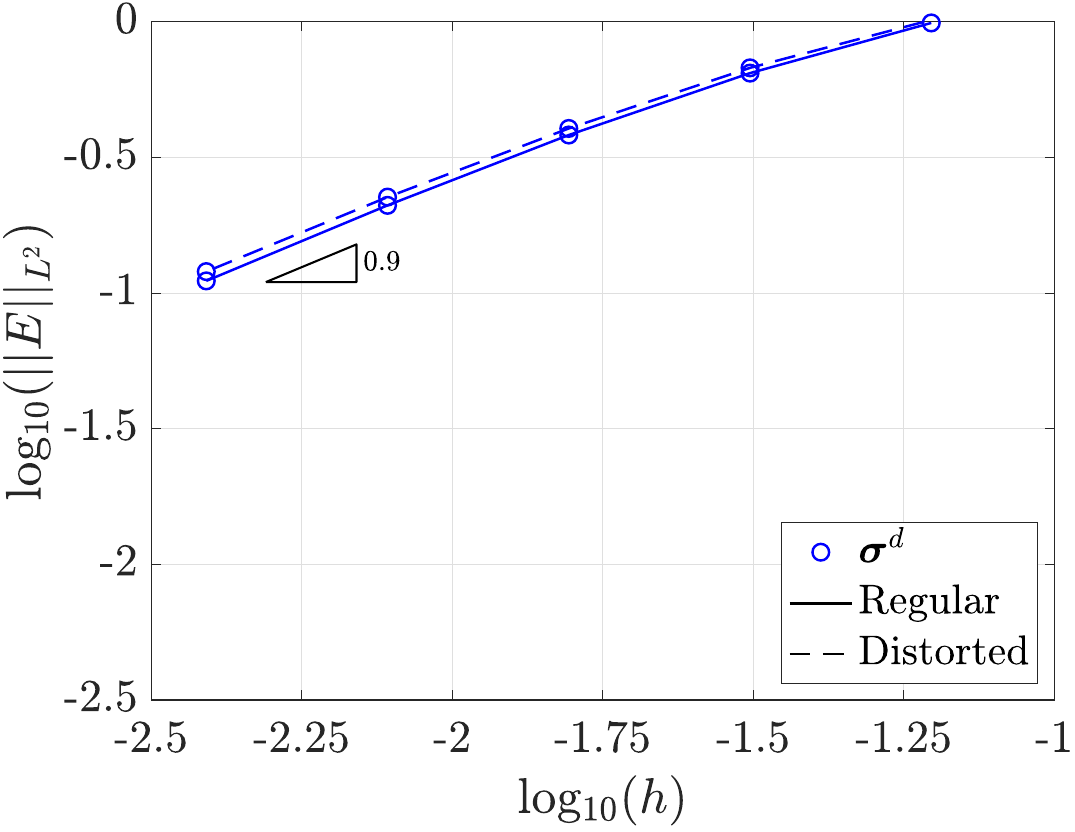}}	
		
	\caption{Synthetic Stokes flow - Mesh convergence of the error of (a) velocity and hybrid velocity, (b) pressure and hybrid pressure, and (c) deviatoric stress tensor, measured in the $\eltwo$ norm as a function of the cell size $h$ on meshes of regular and distorted quadrilateral cells. }
	\label{fig:convSqua}
\end{figure}

Setting $\tauP=10^{-1}$, Figure~\ref{fig:convSqua} displays the relative error, measured in the $\eltwo$ norm as a function of the characteristic cell size $h$,  for cell and face velocity (left),  cell and face pressure (centre), and deviatoric stress tensor (right). 
The results show that, upon mesh refinement,  the FCFV method with hybrid pressure achieves the optimal convergence rate of order one for velocity, pressure and stress tensor.
As already observed for the FCFV method in~\cite{RS-SGH:2018_FCFV1,vila2022non,vila2023benchmarking,Vieira-VGSH-24},  the convergence rates are preserved even in the presence of distorted cells. Moreover,  the error levels achieved on regular and distorted meshes remain practically identical for all variables, showing that the hybrid pressure formulation is less sensitive to cell distortion than traditional CCFV and VCFV methods~\cite{diskin2010comparison,diskin2011comparison}.
\begin{figure}[!htb]
	\centering
	\subfigure[]{\includegraphics[width=0.32\textwidth]{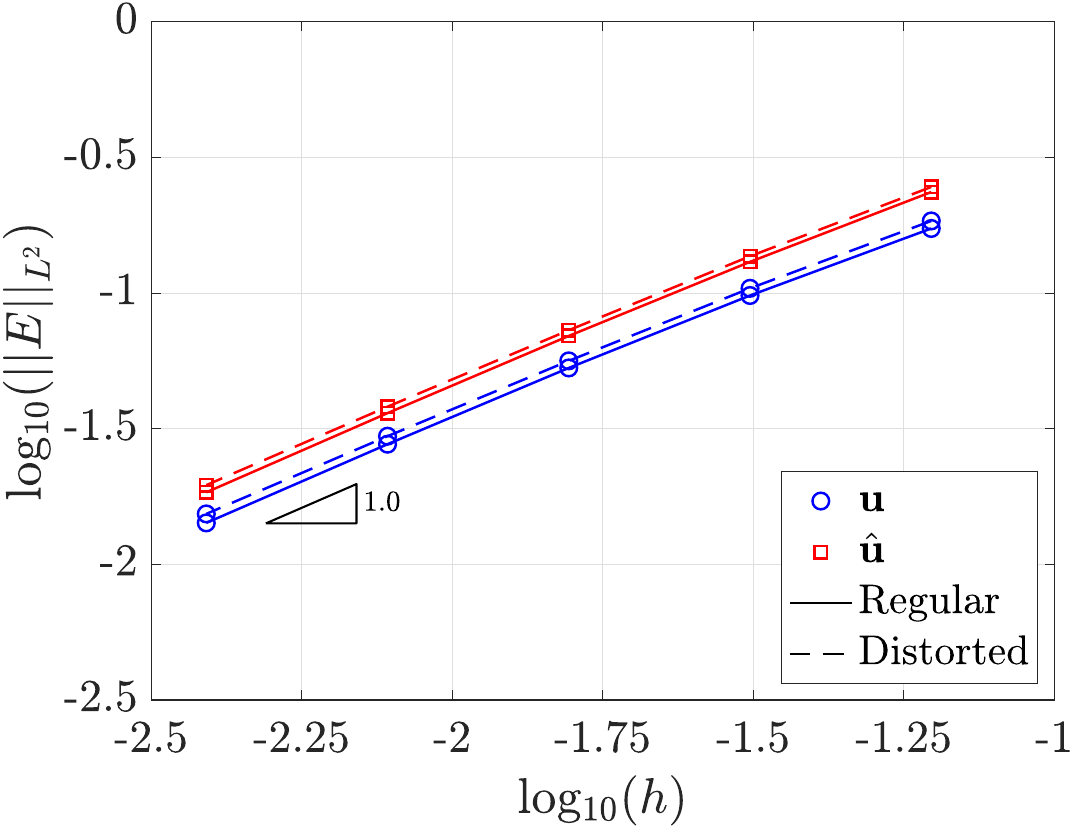}}
	\hspace{2pt}	
	\subfigure[]{\includegraphics[width=0.32\textwidth]{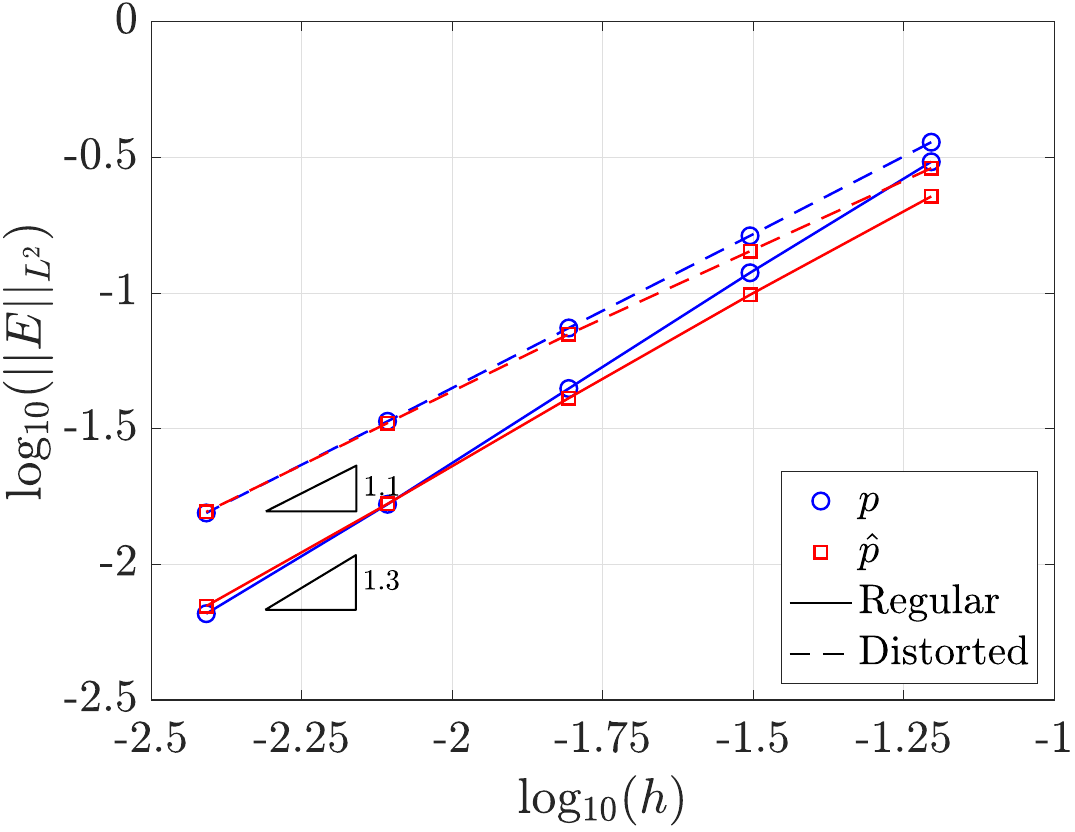}}
	\hspace{2pt}	
	\subfigure[]{\includegraphics[width=0.32\textwidth]{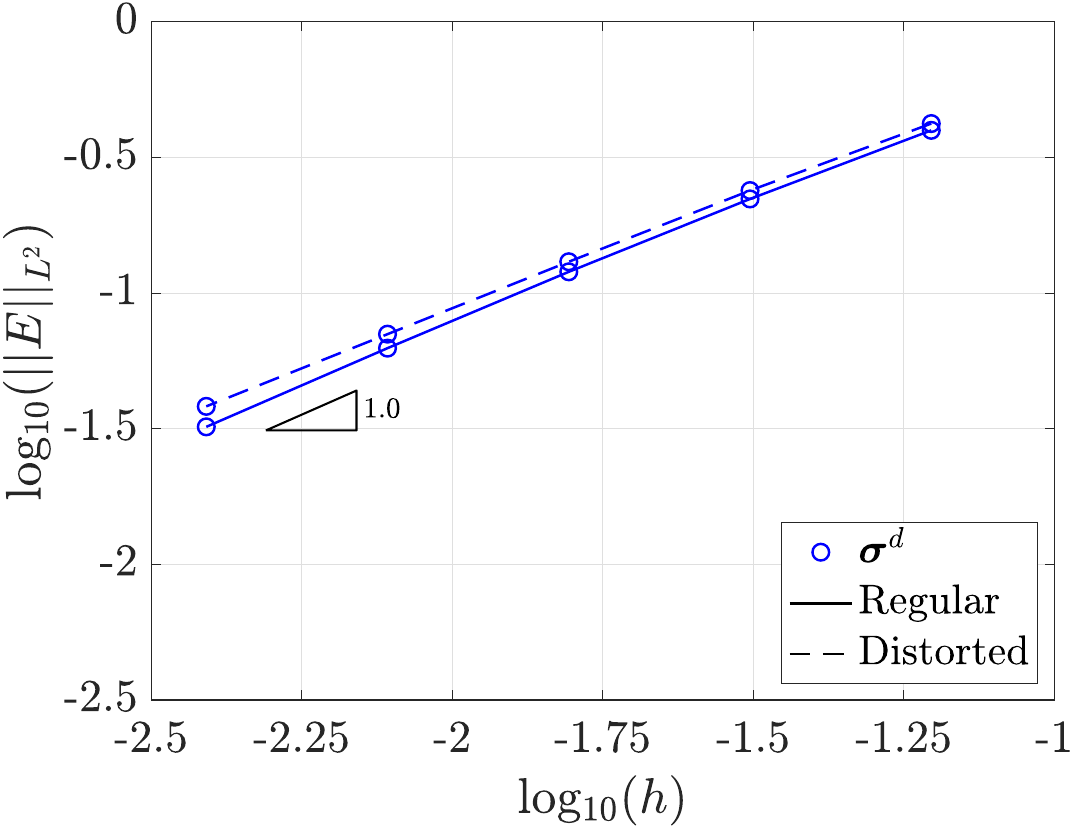}}	
	
	\caption{Synthetic Stokes flow - Mesh convergence of the error of (a) velocity and hybrid velocity, (b) pressure and hybrid pressure, and (c) deviatoric stress tensor, measured in the $\eltwo$ norm as a function of the cell size $h$ on meshes of regular and distorted triangular cells. }
	\label{fig:convStri}
\end{figure}

The test is repeated using a set of meshes with triangular cells,  obtained by splitting each quadrilateral in the above meshes by means of one diagonal. The resulting $i$-th mesh consists of $(8 \times 2^i)\times(8 \times 2^i)\times 2$ triangles. The first level of refinement of these computational meshes has been presented in Figure~\ref{fig:meshes2DStokesTRI} and~\ref{fig:meshes2DStokesTRId} for regular and distorted cells.
The corresponding convergence of the error for cell and face velocity,  cell and face pressure, and deviatoric stress tensor is reported in Figure~\ref{fig:convStri}.
The results confirm the optimal convergence rate of order one and almost identical errors for cell and face velocity and deviatoric stress tensor, on both regular and distorted meshes.
Cell and face pressure converge optimally with order one on distorted meshes, whereas, in this problem, they experience superconvergence on regular meshes, achieving order $1.3$.

\subsubsection{Influence of the stabilisation parameter on mass conservation}
\label{sc:InfluenceTauStokes}

\begin{figure}[!htb]
	\centering
	\subfigure[]{\includegraphics[width=0.32\textwidth]{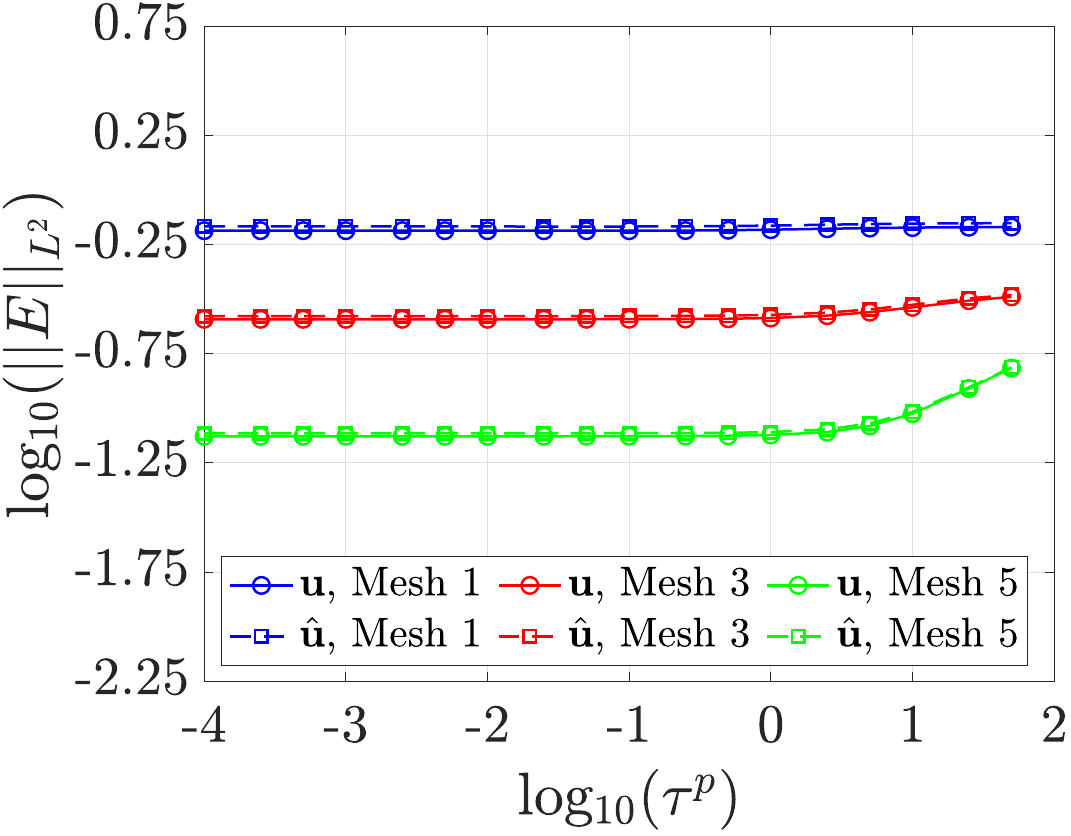}}
	\hspace{2pt}	
	\subfigure[]{\includegraphics[width=0.32\textwidth]{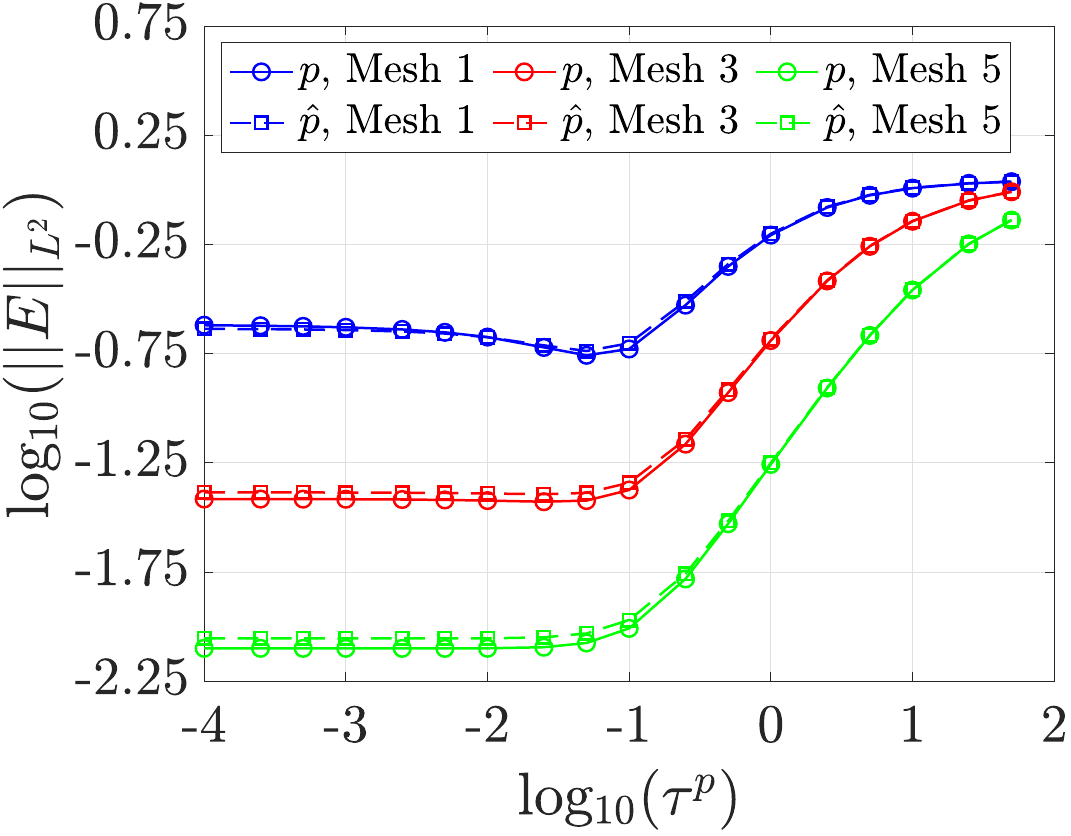}}
	\hspace{2pt}	
	\subfigure[]{\includegraphics[width=0.32\textwidth]{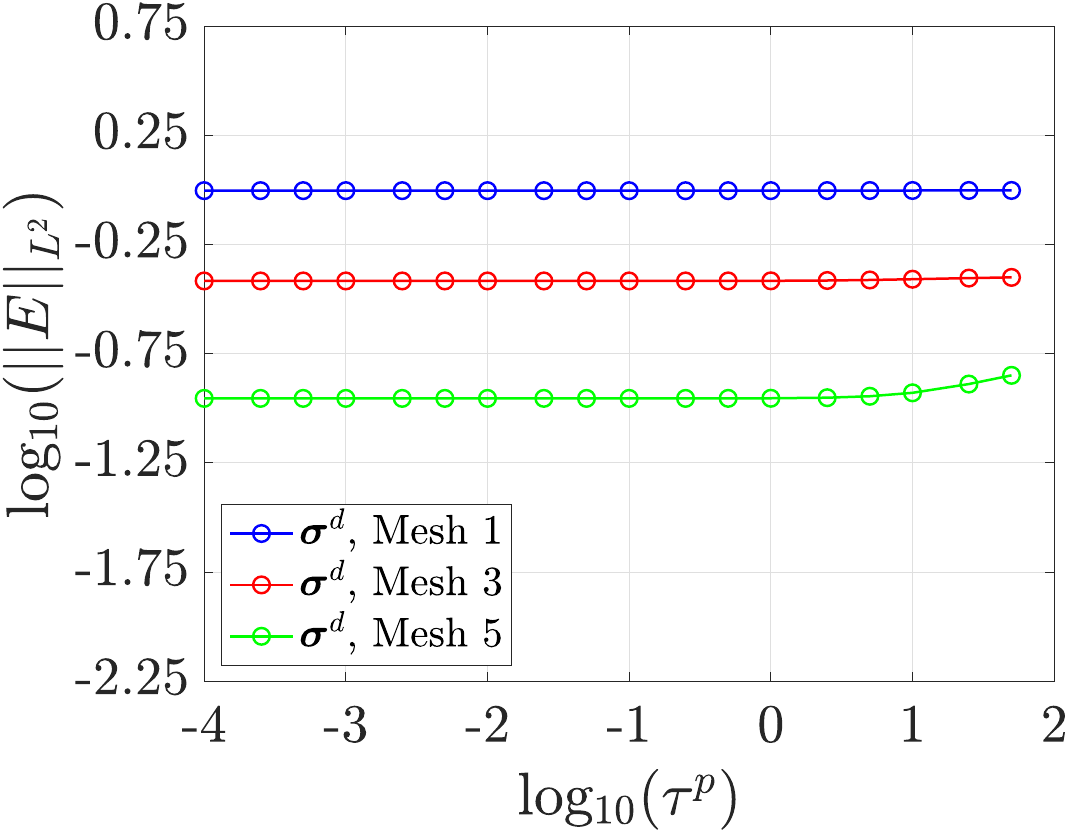}}	
	
	\subfigure[]{\includegraphics[width=0.32\textwidth]{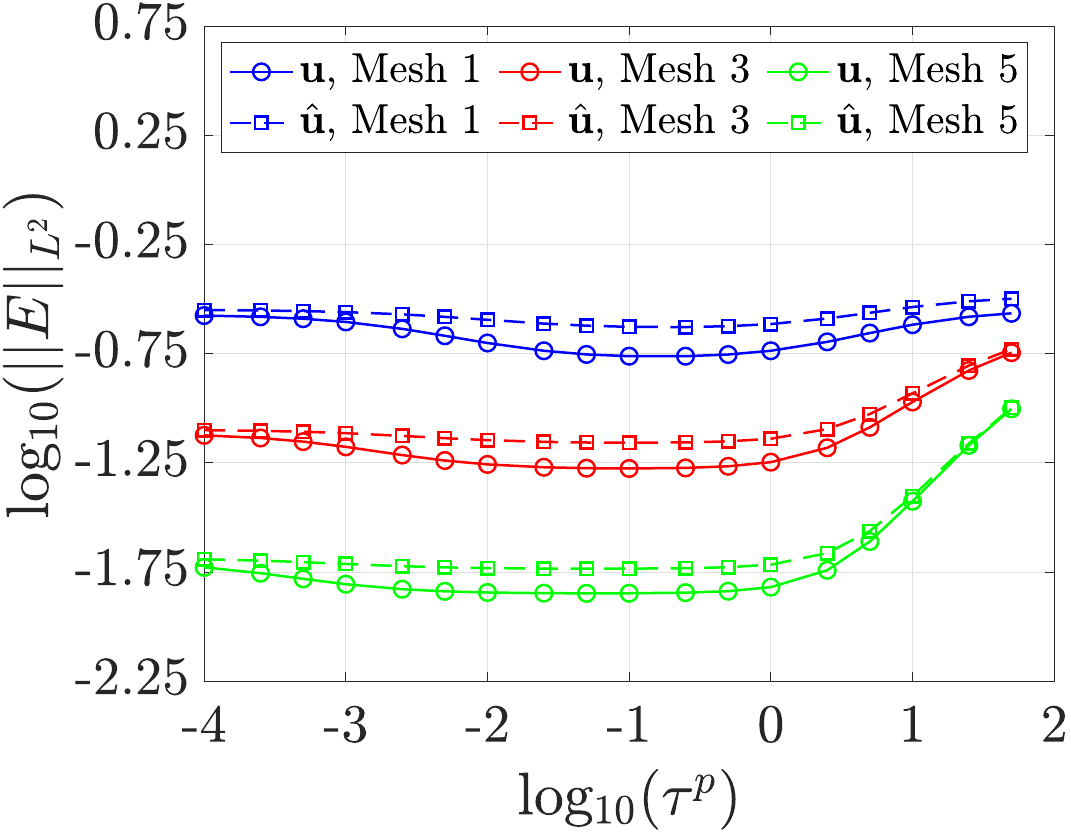}}
	\hspace{2pt}	
	\subfigure[]{\includegraphics[width=0.32\textwidth]{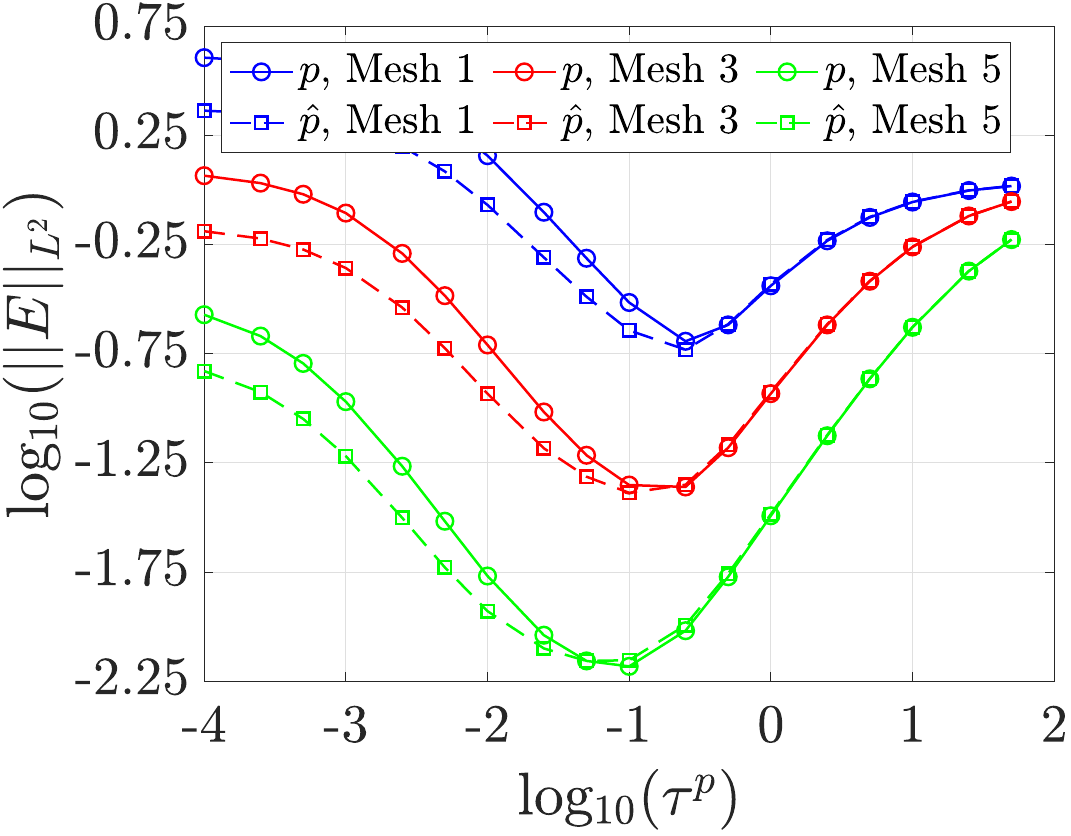}}
	\hspace{2pt}	
	\subfigure[]{\includegraphics[width=0.32\textwidth]{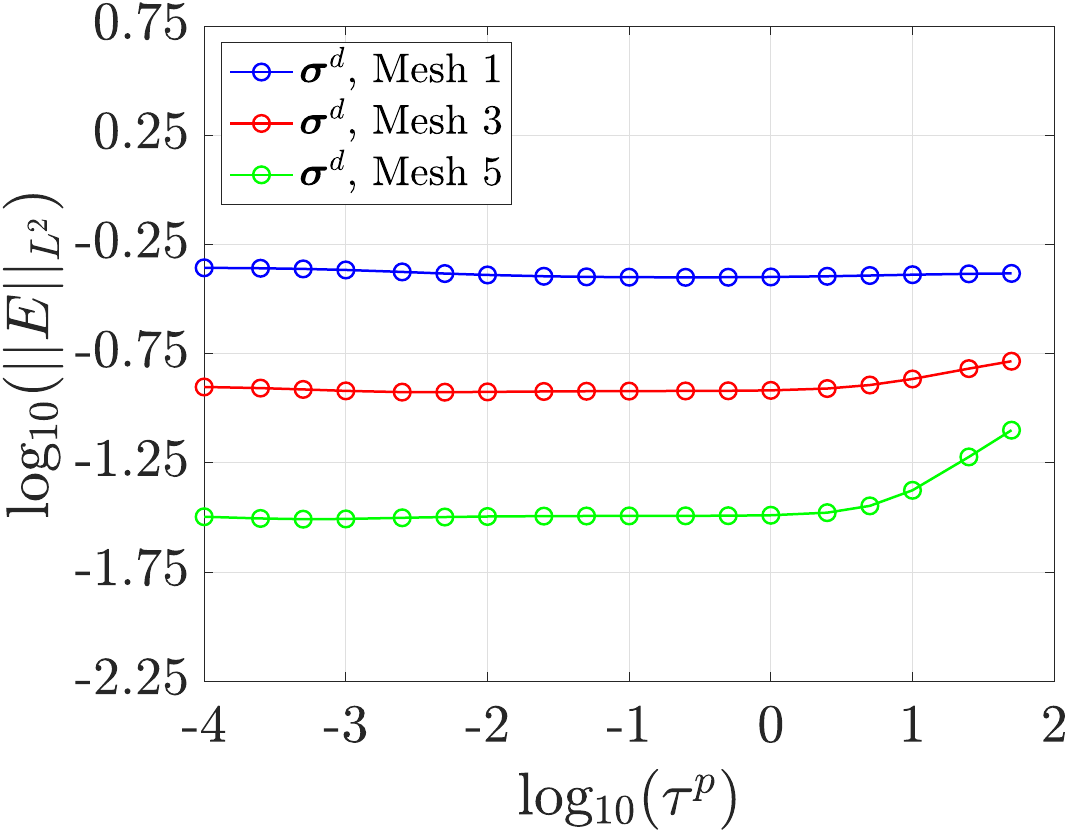}}	
	
	\caption{Synthetic Stokes flow - Sensitivity of the error for (a,d) velocity and hybrid velocity,  (b,e) pressure and hybrid pressure, and (c,f) deviatoric stress tensor,  measured in the $\eltwo$ norm as a function of the stabilisation coefficient $\tauP$ for different levels of mesh refinement.  Top row: quadrilateral meshes. Bottom row: triangular meshes.}
	\label{fig:sensTauUPL}
\end{figure}

As previously mentioned (see Section~\ref{sc:Fluxes}), the hybrid pressure formulation requires the introduction of new stabilisation coefficient in the definition of the mass flux.
A sensitivity study is presented to showcase the robustness of the method to the choice of this parameter. 
Figure~\ref{fig:sensTauUPL} displays the $\eltwo$ norm of the error for cell and face velocity (left), cell and face pressure (centre), and deviatoric stress tensor (right) as a function of $\tauP$, for different levels of refinement of the meshes of quadrilateral and triangular cells employed in the previous convergence study.

The accuracy of the approximation of $\bu$, $\bhu$, and $\bsigmaD$ appears to be almost independent of the choice of $\tauP$, for both quadrilateral and triangular cells and for any mesh refinement considered.
On the contrary, the accuracy of the cell and face pressure is affected by the choice of $\tauP$, but a common trend can be identified, for different meshes and cell types. Given the results in Figure~\ref{fig:sensTauUPL},  all the simulations in the present work are performed using $\tauP=10^{-1}$.

It is worth recalling that, according to the definition of the mass flux~\eqref{eq:traceMass}, $\tauP$ has a direct influence on the fulfilment of the incompressibility equation. In particular, $\tauP$ can be interpreted as a measure of the compressibility effect introduced by enforcing mass conservation weakly.
To assess the accuracy of the hybrid pressure formulation to impose the incompressiblity constraint, let us introduce the cell mass flux
\begin{equation}\label{eq:cellFlux}
\fluxEl := \int_{\partial\Omega_e \setminus \Ga[D]} \bhu \cdot \bn\, d\Ga + \int_{\partial\Omega_e \cap \Ga[D]} \bu_D \cdot \bn\, d\Ga .
\end{equation}

\begin{figure}[!htb]
	\centering
	\subfigure[]{\includegraphics[width=0.4\textwidth]{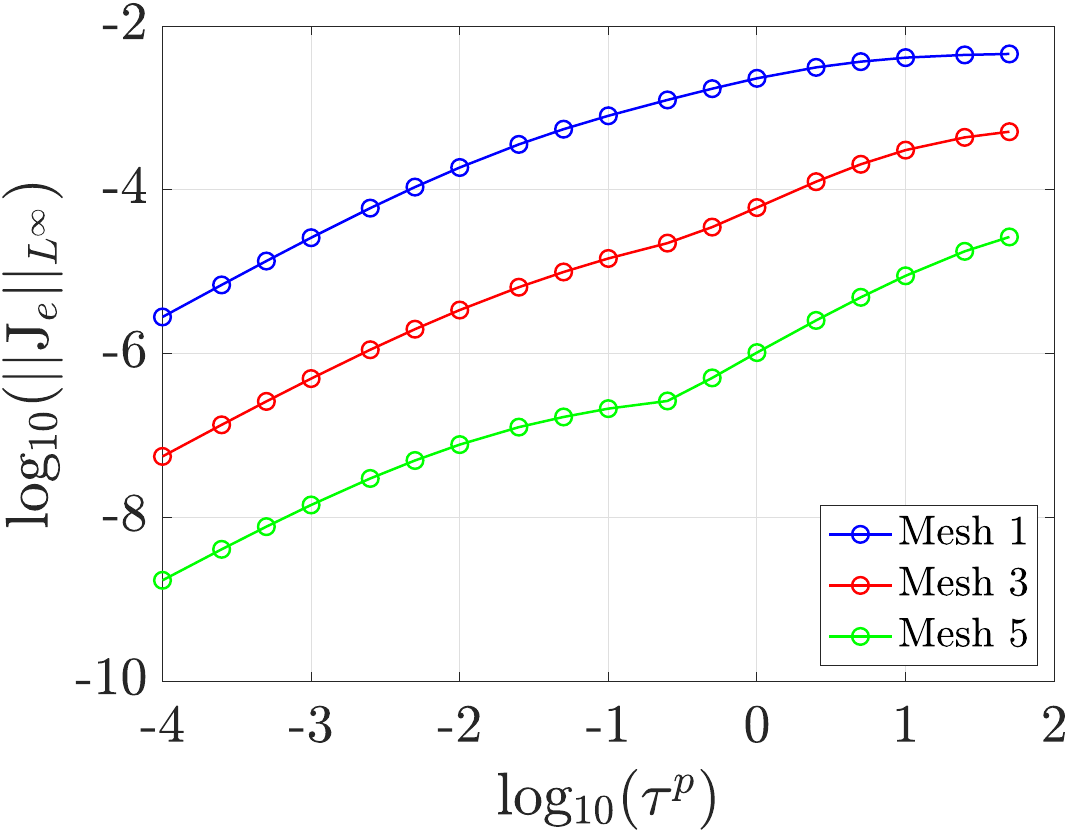}}
	\hspace{5pt}	
	\subfigure[]{\includegraphics[width=0.4\textwidth]{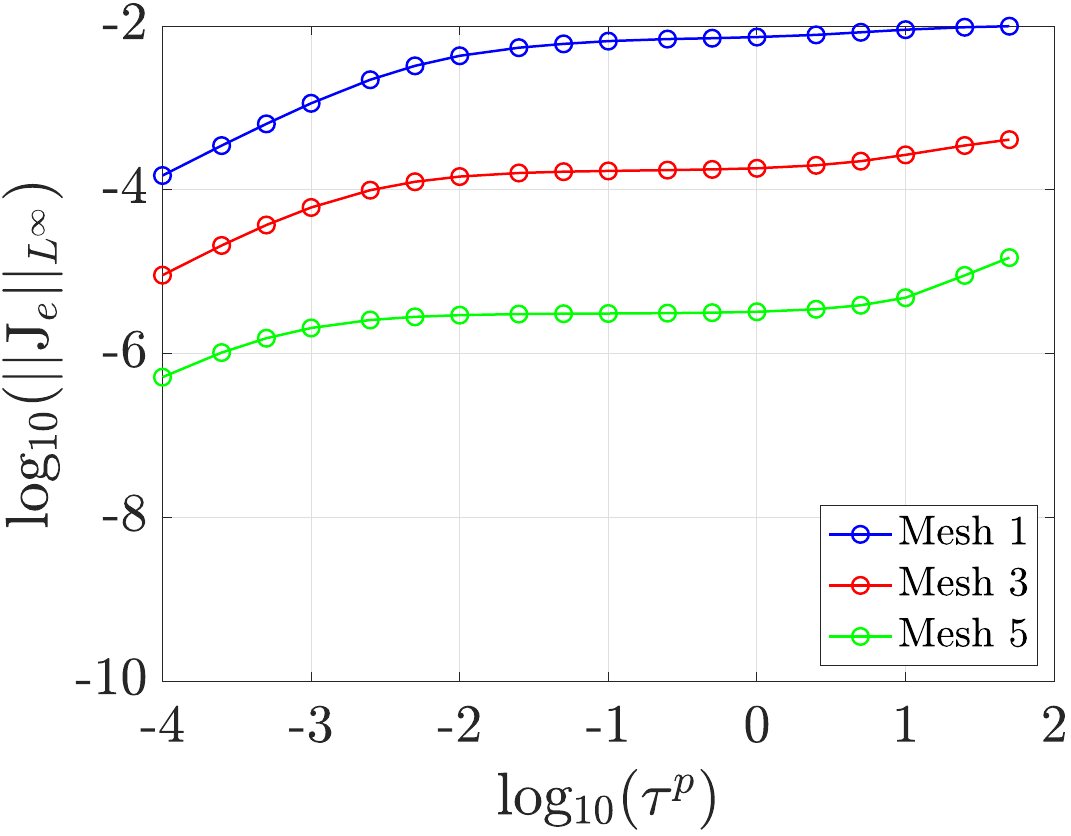}}
	
	\caption{Synthetic Stokes flow - Sensitivity of the cell mass flux $\fluxEl$, measured in the $\elinf$ norm, as a function of the stabilisation coefficient $\tauP$ on different levels of mesh refinment of (a) quadrilateral and (b) triangular cells.}
	\label{fig:sensTauDiv}
\end{figure}

Figure~\ref{fig:sensTauDiv} displays the maximum norm of the cell mass flux $\fluxEl$, for different levels of mesh refinement and cell types, when the value of the stabilisation coefficient $\tauP$ is modified.
Note that while the exact value of $\fluxEl$ is zero for each cell, this quantity serves as a measure of the error in the approximation of the cell's mass flux.
Of course, as $\tauP \rightarrow 0$,  the approximation of the velocity tends to be perfectly divergence-free and $\fluxEl \rightarrow 0$. On the contrary, for large values of $\tauP$, the incompressibility constraint is not well enforced but the error decreases as the mesh is refined.
Nonetheless,  the maximum mass flux error,  reported in Table~\ref{tab:errMass} for $\tauP = 10^{-1}$ and different  meshes of quadrilateral and triangular cells,  is always at least one order of magnitude smaller than the corresponding best approximation error of order $h$ that can be achieved by a given mesh.  
Hence,  although the hybrid pressure formulation only enforces incompressibility in a weak sense, the error introduced in the mass conservation of a given cell is negligible with respect to the spatial discretisation error. Moreover, at the global level, upon summing the contributions $\fluxEl$ of all cells $\Omega_e, \ e=1,\ldots,\numel$, the mass flux error achieves values of order $10^{-12}$ for all meshes, and the mass conservation equation is thus verified at machine precision.
Similar results, not reported here for brevity, were obtained using distorted meshes, with maximum cell mass flux errors ranging from $0.92 \times 10^{-2}$ to $0.13 \times 10^{-4}$.
\begin{table}[!htb]
\caption{Synthetic Stokes flow - Maximum cell mass flux $\fluxEl$ with $\tauP = 10^{-1}$ for different levels of refinement of the regular meshes.}
\centering
\begin{tabular}{|l|c|c|c|c|c|}\hline
 Mesh                  &  $1$                               &   $2$                             &   $3$                              & $4$                              & $5$                          \\ \hline
 $h$                     &  $0.62 \times 10^{-1}$ &   $0.31 \times 10^{-1}$ &   $0.16 \times 10^{-1}$ &   $0.78 \times 10^{-2}$ &   $0.39 \times 10^{-2}$ \\ \hline
 Quadrilaterals    &  $0.80 \times 10^{-3}$ &   $0.12 \times 10^{-3}$ &   $0.14 \times 10^{-4}$ &   $0.16 \times 10^{-5}$ &   $0.21 \times 10^{-6}$  \\
 Triangles            &  $0.65 \times 10^{-2}$ &   $0.11 \times 10^{-2}$ &   $0.17 \times 10^{-3}$ &   $0.23 \times 10^{-4}$ &   $0.31 \times 10^{-5}$  \\ \hline
\end{tabular}	
\label{tab:errMass}
\end{table}

Additional numerical tests have been performed for the Stokes flow, including the two-dimensional cavity flow and the flow past a sphere in 3D.  The results confirmed the previously observed optimal convergence properties of the hybrid pressure formulation, its accuracy using different cell types (quadrilaterals, triangles, and tetrahedra), and its robustness in the presence of unstructured, distorted, and stretched meshes.
It is worth noticing that all results were comparable with the ones provided by the previously published FCFV method~\cite{RS-SGH:2018_FCFV1}, and they are not reported here for brevity.  In particular,  it was observed that the advantages of one approach over the other were mainly problem-dependent when dealing with the Stokes equations.
Hence,  the remaining of this work will focus on Navier-Stokes flows to showcase the superiority of the hybrid pressure formulation in the presence of convection phenomena.

\subsection{Navier-Stokes Couette flow}
\label{sc:ConvergenceNS}

In this section, the numerical convergence analysis of the hybrid pressure formulation is performed for a two-dimensional incompressible Navier-Stokes flow with analytical solution.
The test case under analysis is the co-axial Couette flow, defined on the computational domain $\Omega := \{ (x_1,x_2) \in \RR^2 \ | \ R_i \leq r(x_1,x_2) \leq R_o \}$,  with $r := (x_1^2+x_2^2)^{1/2}$.
The setup of the problem represents an annulus centred in $(0,0)$, with inner and outer radii $R_i=1$ and $R_o=2$, respectively, rotating with imposed angular velocities $\omega_i = 0$ and $\omega_o = 0.5$, respectively. 

The analytical velocity $\bu = (u_r,u_\phi)^T$ and pressure, expressed in polar coordinates, are given by
\begin{equation}\label{eq:Couette_Statment}
\bu(r) =
\begin{Bmatrix}
0 \\
C_1 r + C_2 \displaystyle\frac{1}{r}
\end{Bmatrix}
, \qquad
p(r)	=  C_1^2 \frac{r^2}{2} + 2 C_1 C_2 \log(r) - \frac{C_2^2}{2 r^2} + C,
\end{equation}
with $C_1 := (\omega_o R_o^2 - \omega_i R_i^2) /(R_o^2 - R_i^2)$ and $C_2 := (\omega_i - \omega_o)R_i^2 R_o^2/(R_o^2 - R_i^2)$, whereas $C$ denotes a constant selected to guarantee that the pressure at $R_o$ is equal to 1. 
Note that the solution pair $(\bu,p)$ does not depend upon viscosity and $\nu$ is set to 1.  Considering $|R_o - R_i|$ as characteristic length of the problem and setting the magnitude of $\bu$ at $r=R_o$ as characteristic velocity,  the corresponding value of the Reynolds number is 1.

\begin{figure}[!htb]
	\centering
	\subfigure[]{\includegraphics[width=0.23\textwidth]{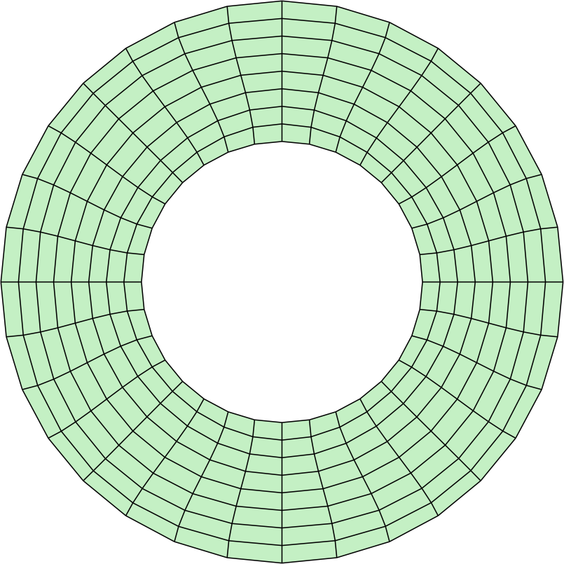}}
	\hspace{5pt}
	\subfigure[]{\includegraphics[width=0.23\textwidth]{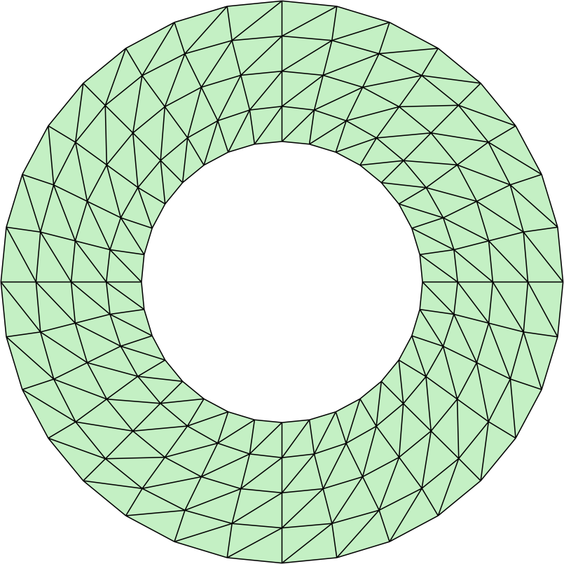}}
	\hspace{5pt}
	\subfigure[]{\includegraphics[width=0.23\textwidth]{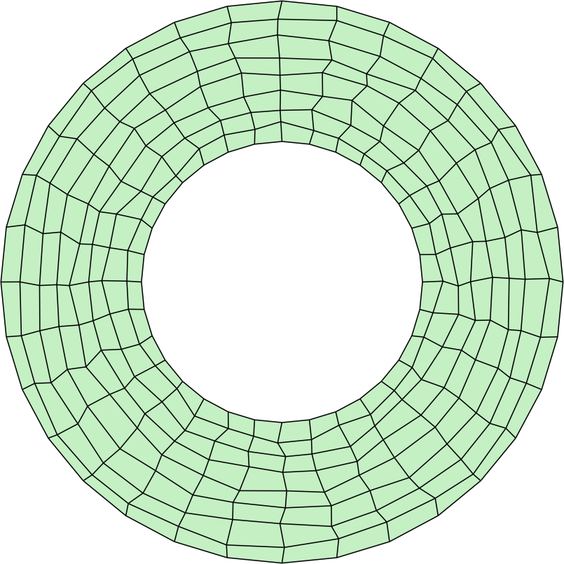}}
	\hspace{5pt}
	\subfigure[]{\includegraphics[width=0.23\textwidth]{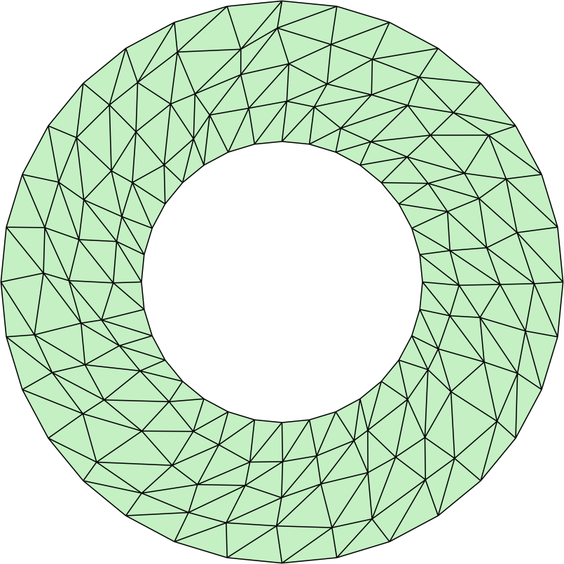}}
	\caption{Couette Navier-Stokes flow - First level of refinement of the meshes of quadrilateral and triangular cells.  (a,b) Regular meshes. (c,d) Distorted meshes.}
	\label{fig:CouetteMeshes}
\end{figure}

Setting $\xi = 5 \times 10^{-2}$ in the definition of the HLL convective stabilisation~\eqref{eq:HLLConv},  the Newton-Raphson solver is executed with a tolerance of $10^{-10}$. 
A mesh convergence study is performed using structured meshes of regular and distorted, quadrilateral and triangular cells. 
Figure~\ref{fig:CouetteMeshes} displays the first level of mesh refinement, featuring $32$ subdivisions in the tangential direction and $8$ in the radial one. In general,  the $i$-th level of mesh refinement is characterised by $(8 \times 2^i) \times (8 \times 2^i)$ cells,  as described in~\cite{Vieira-VGSH-24}.
Similarly to the experiments presented in Section~\ref{sc:ConvergenceStokes}, the effect of cell distortion is evaluated introducing a random perturbation of the internal mesh nodes.

\begin{figure}[!htb]
	\centering
	\subfigure[]{\includegraphics[width=0.32\textwidth]{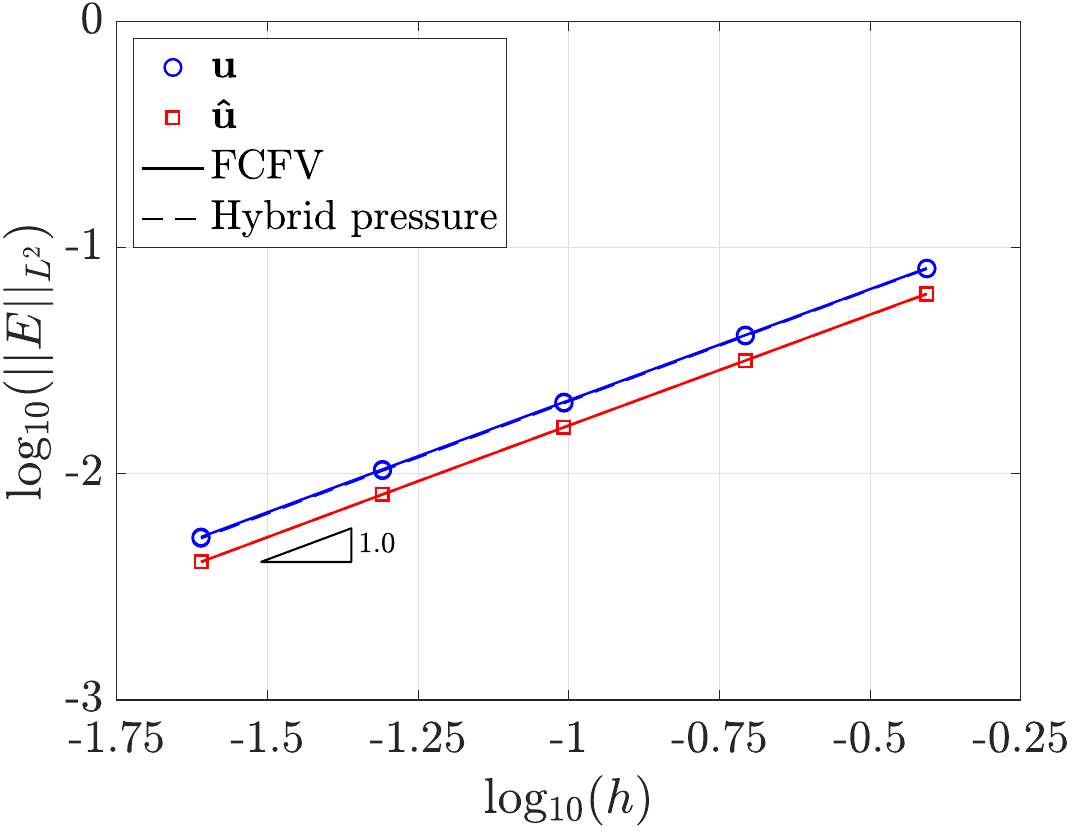}}
	\hspace{2pt}	
	\subfigure[]{\includegraphics[width=0.32\textwidth]{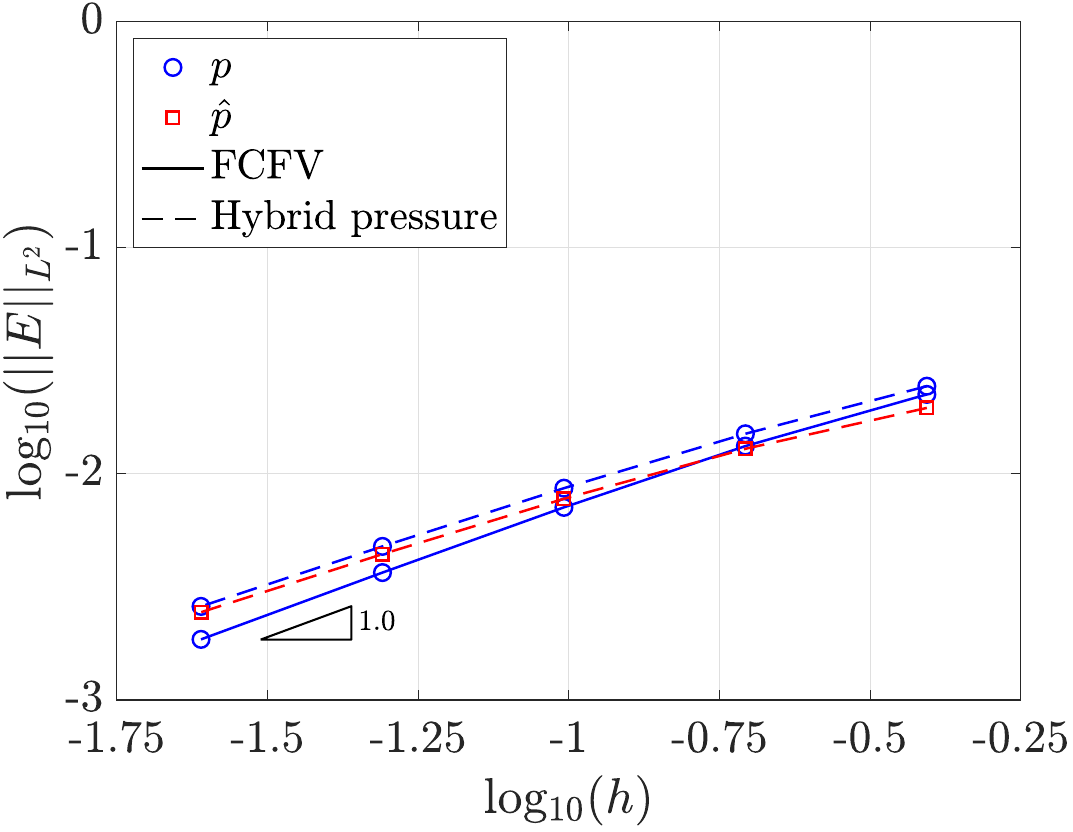}\label{fig:convNSquaPreg}}
	\hspace{2pt}	
	\subfigure[]{\includegraphics[width=0.32\textwidth]{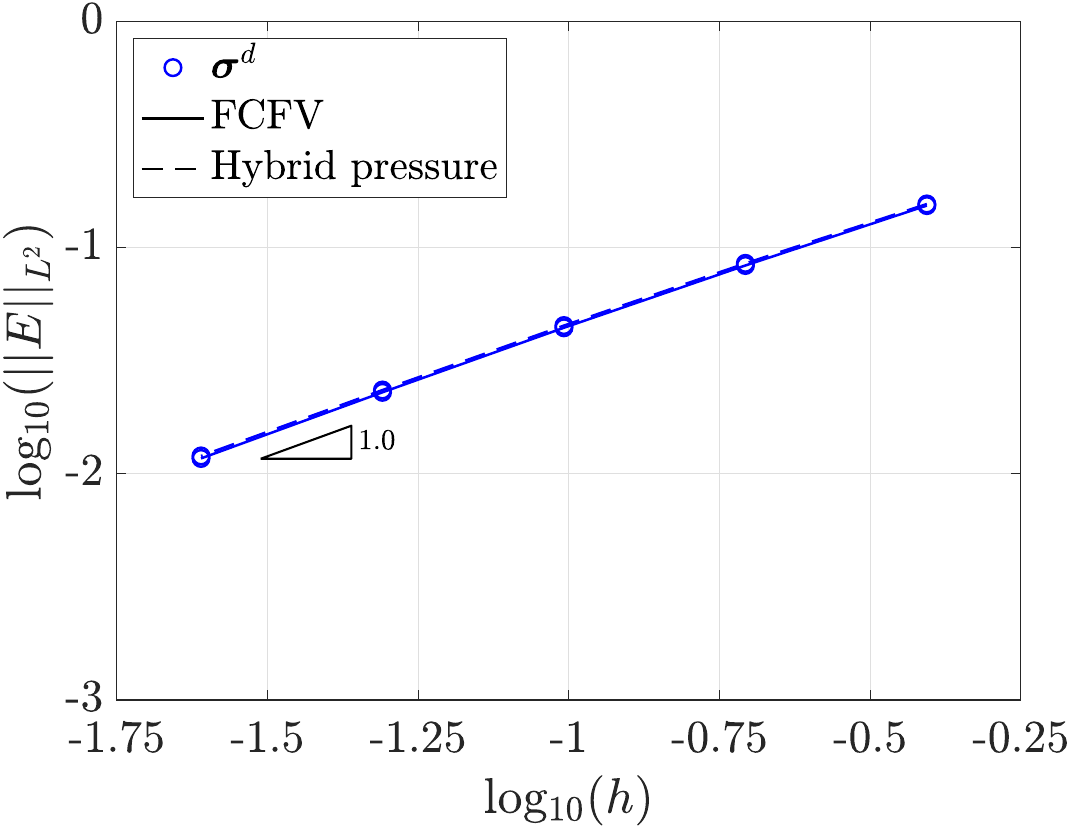}}	
	
	\subfigure[]{\includegraphics[width=0.32\textwidth]{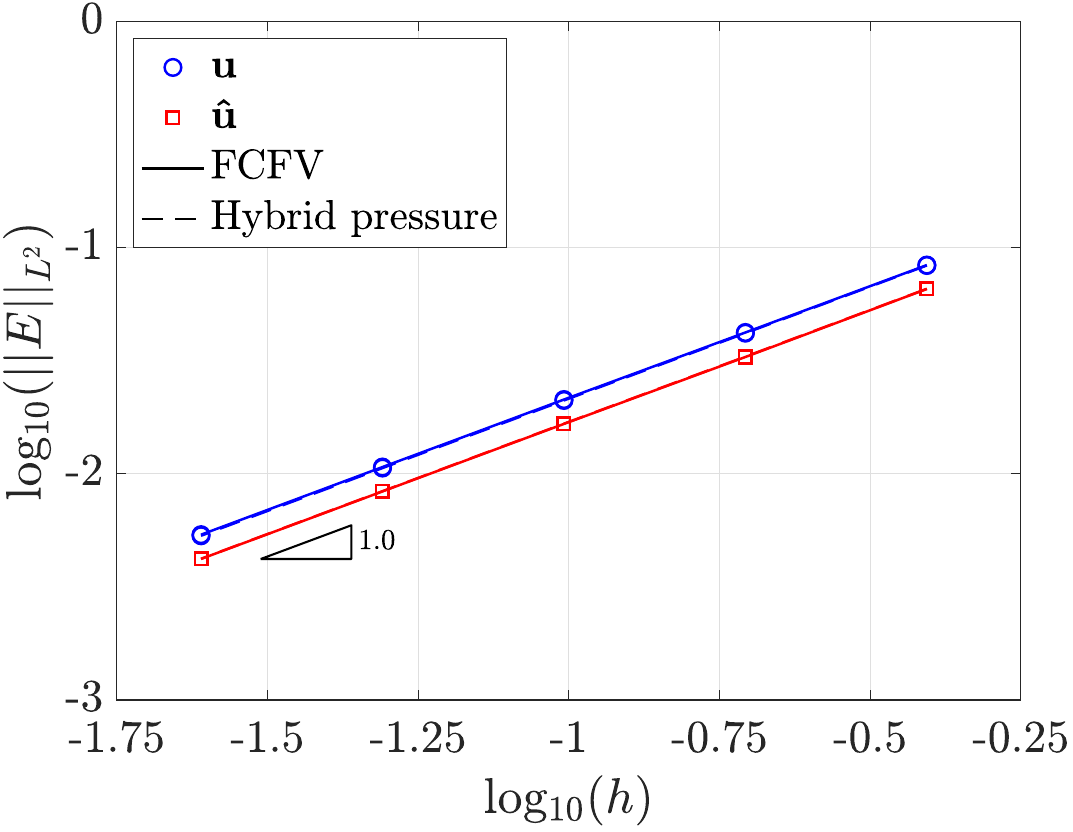}}
	\hspace{2pt}	
	\subfigure[]{\includegraphics[width=0.32\textwidth]{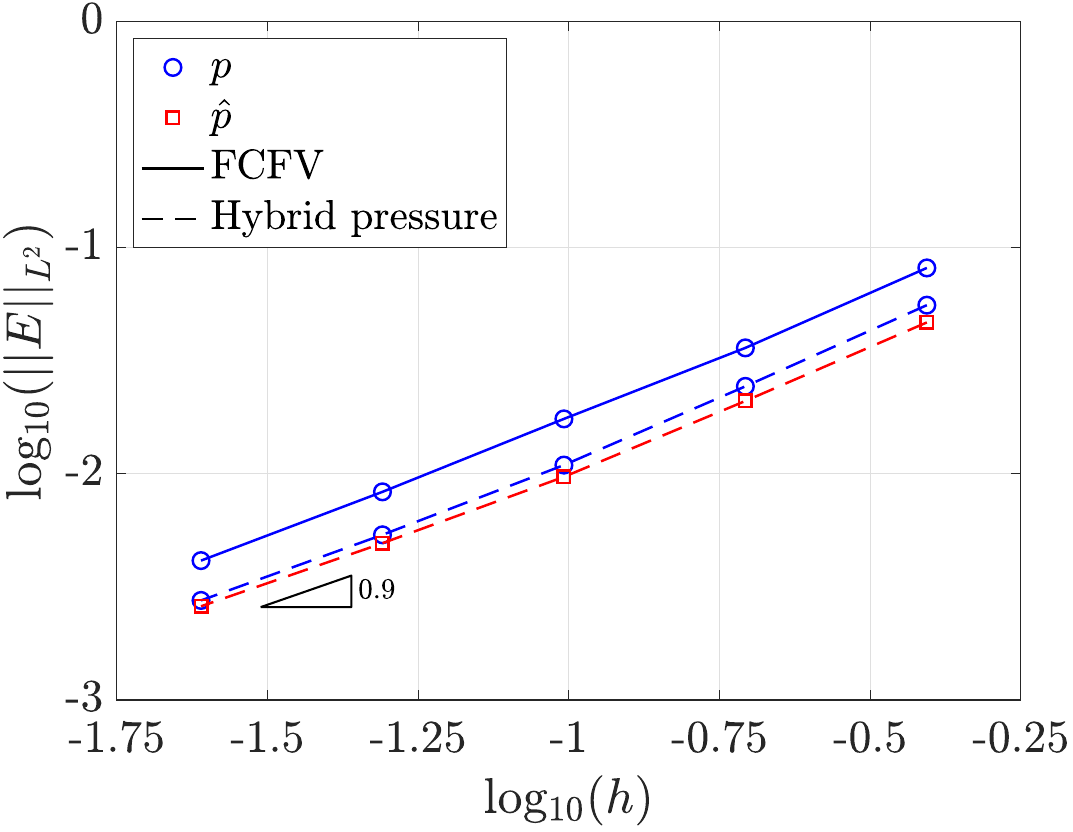}\label{fig:convNSquaPdist}}
	\hspace{2pt}	
	\subfigure[]{\includegraphics[width=0.32\textwidth]{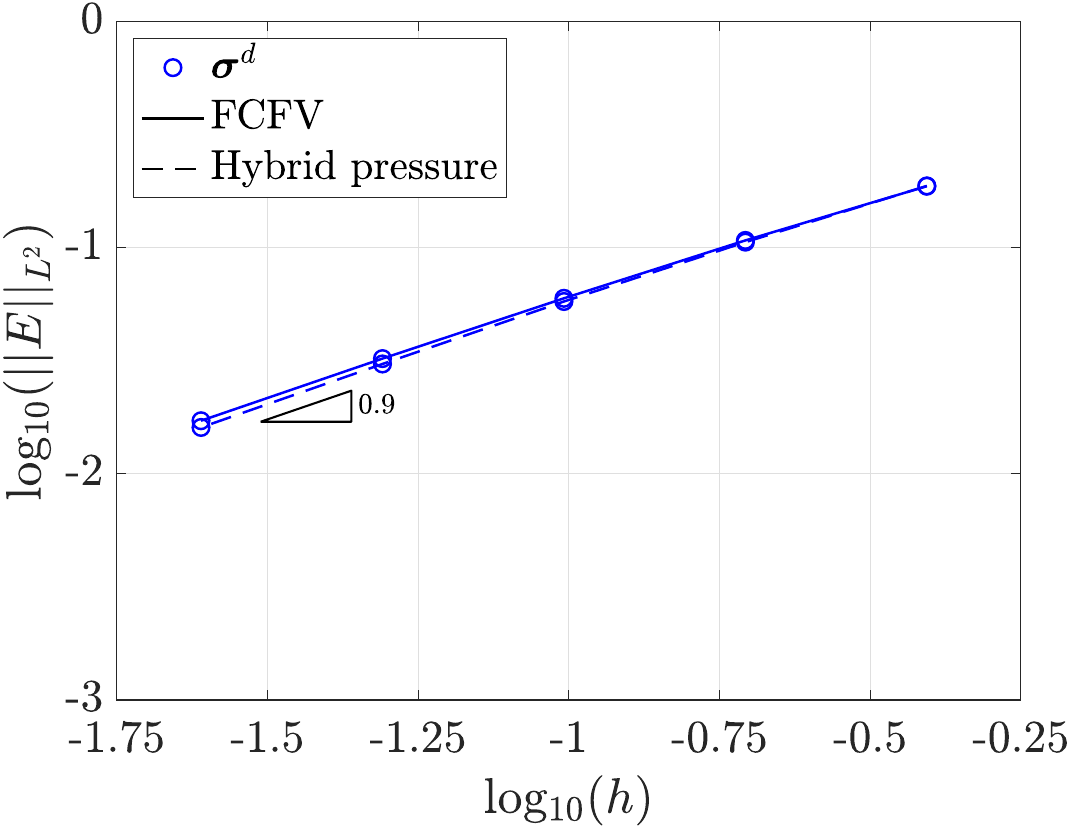}}	
	
	\caption{Couette Navier-Stokes flow - Mesh convergence of the error of (a,d) velocity and hybrid velocity, (b,d) pressure and hybrid pressure, and (c,f) deviatoric stress tensor, measured in the $\eltwo$ norm as a function of the cell size $h$ on meshes of quadrilateral cells. Top row: regular meshes. Bottom row: distorted meshes.}
	\label{fig:convNSqua}
\end{figure}

\begin{figure}[!htb]
	\centering
	\subfigure[]{\includegraphics[width=0.32\textwidth]{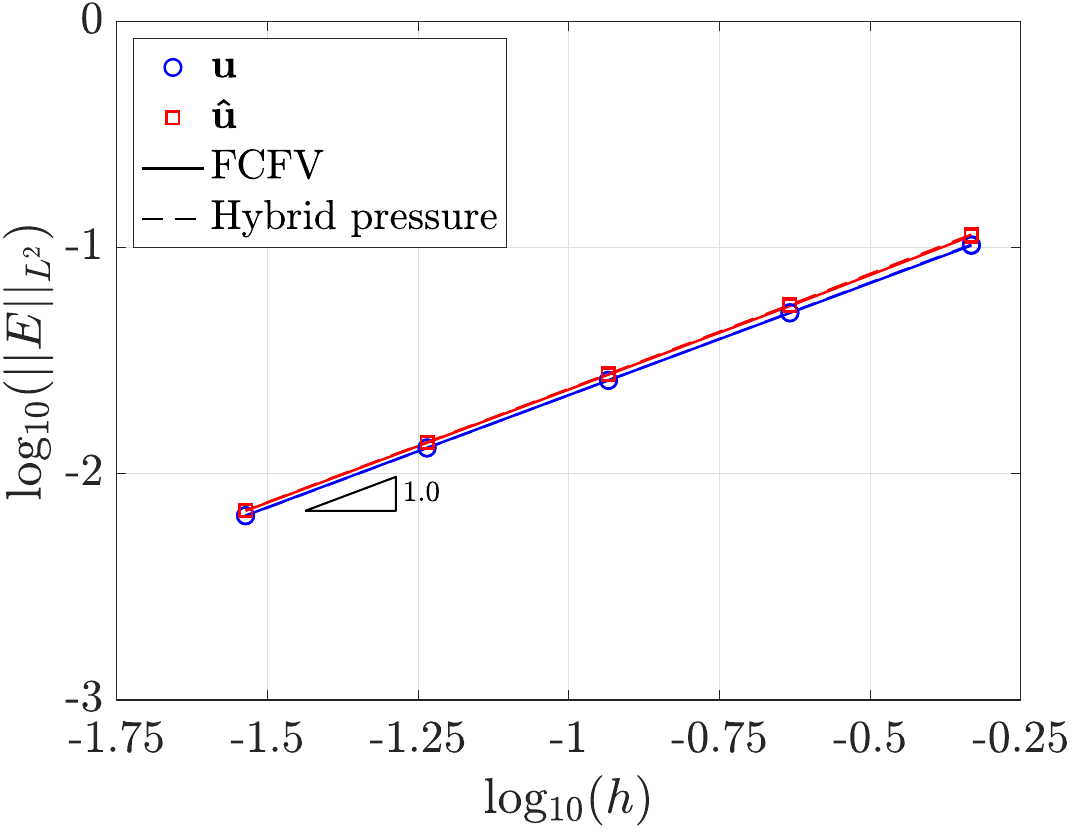}}
	\hspace{2pt}	
	\subfigure[]{\includegraphics[width=0.32\textwidth]{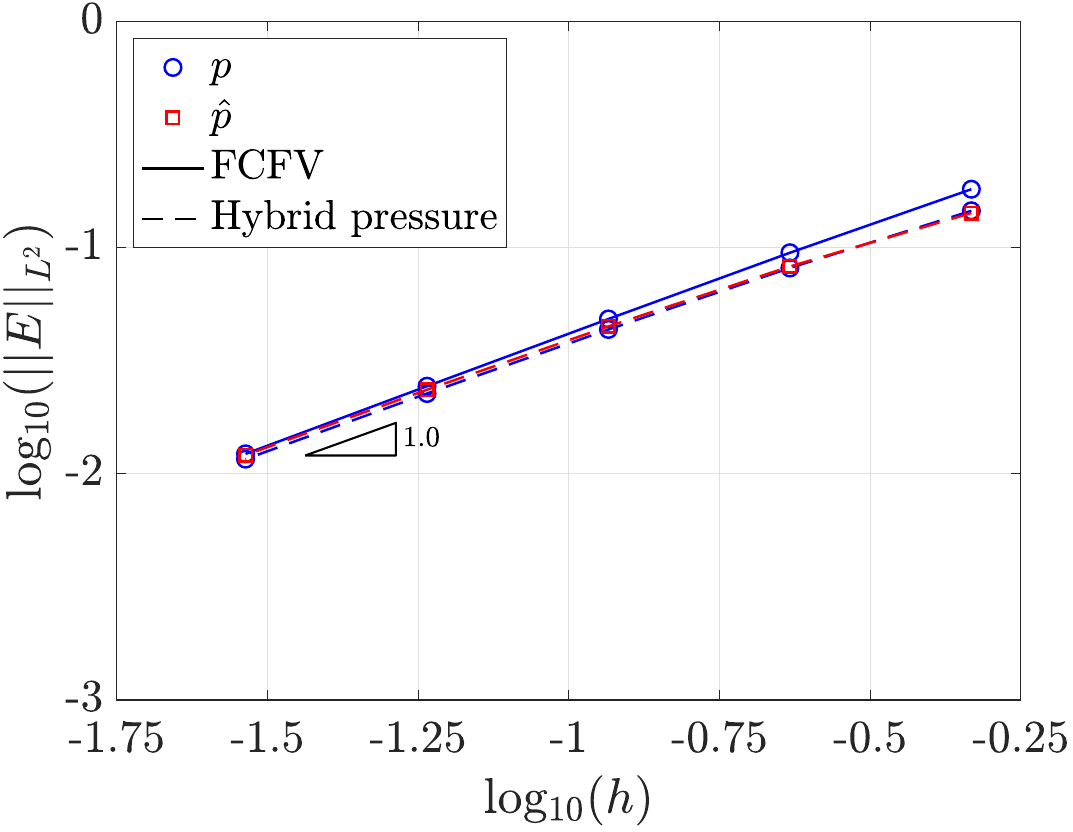}}
	\hspace{2pt}	
	\subfigure[]{\includegraphics[width=0.32\textwidth]{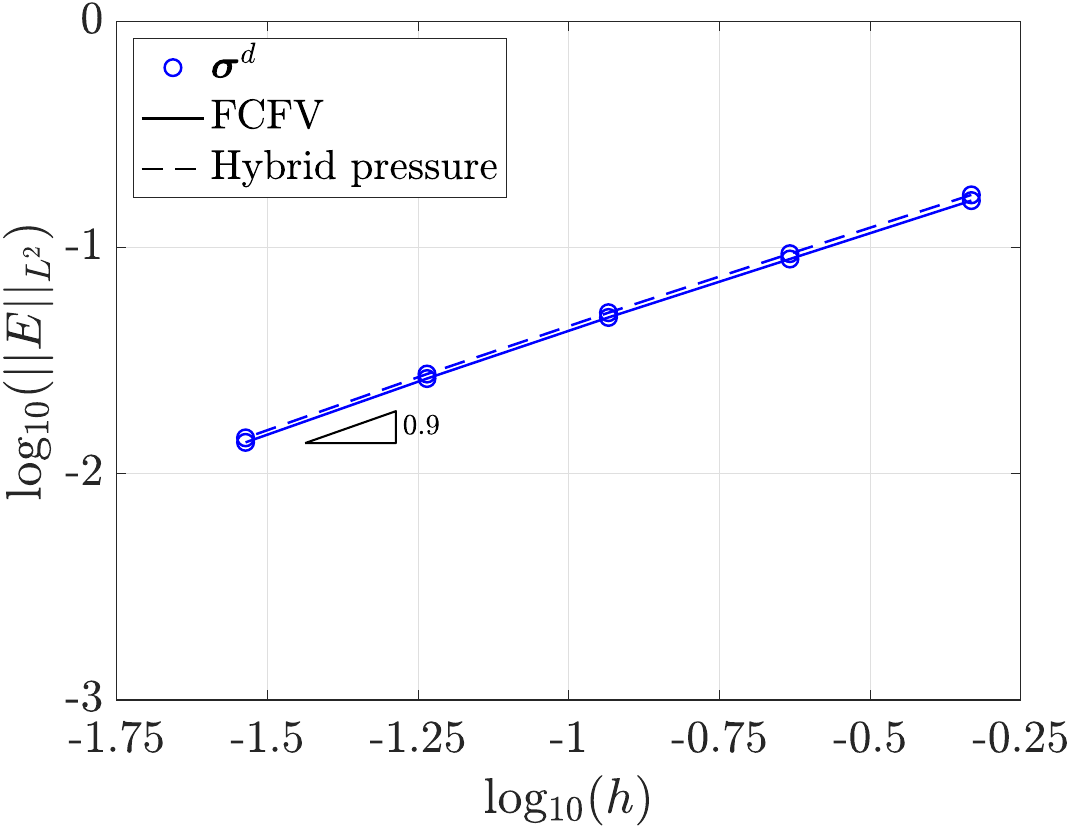}}	
	
	\subfigure[]{\includegraphics[width=0.32\textwidth]{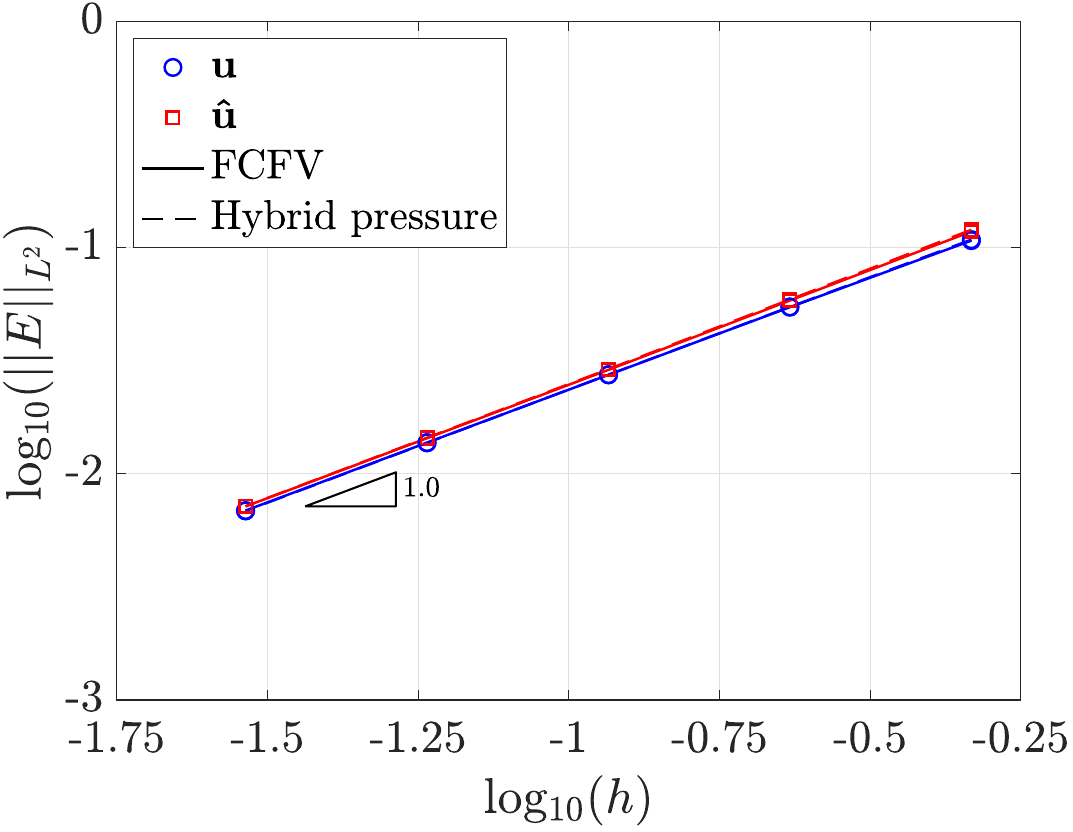}}
	\hspace{2pt}	
	\subfigure[]{\includegraphics[width=0.32\textwidth]{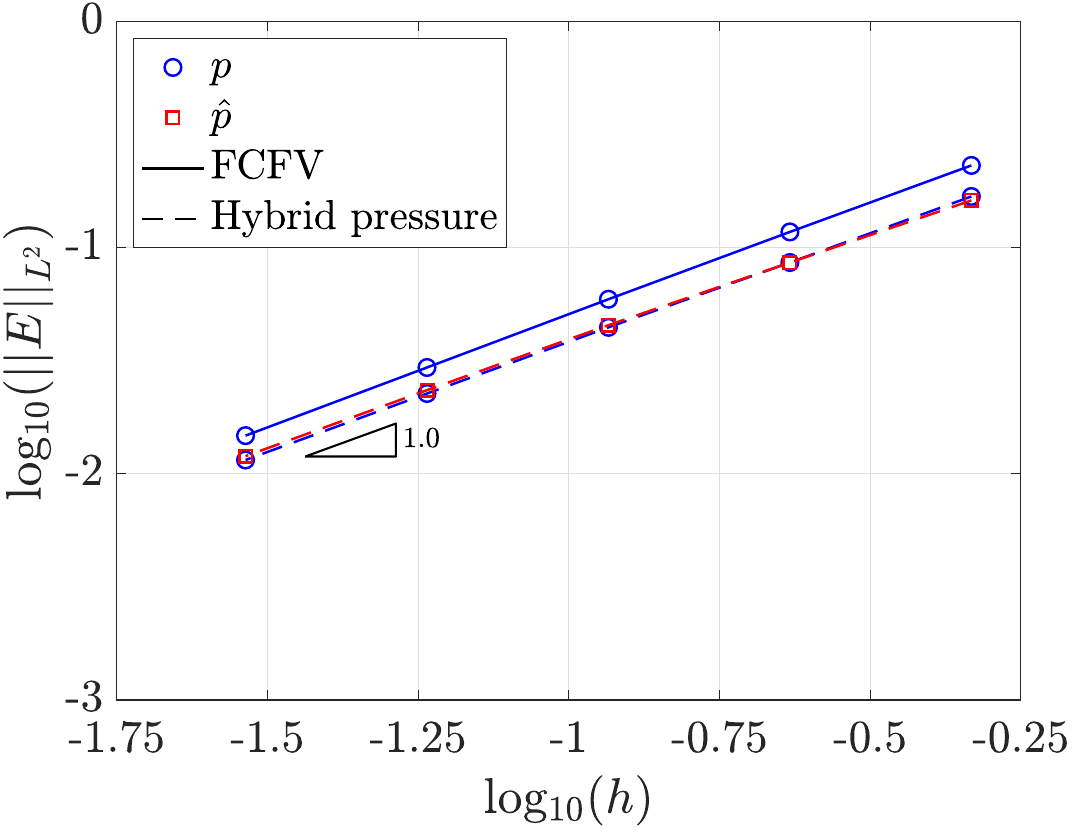}}
	\hspace{2pt}	
	\subfigure[]{\includegraphics[width=0.32\textwidth]{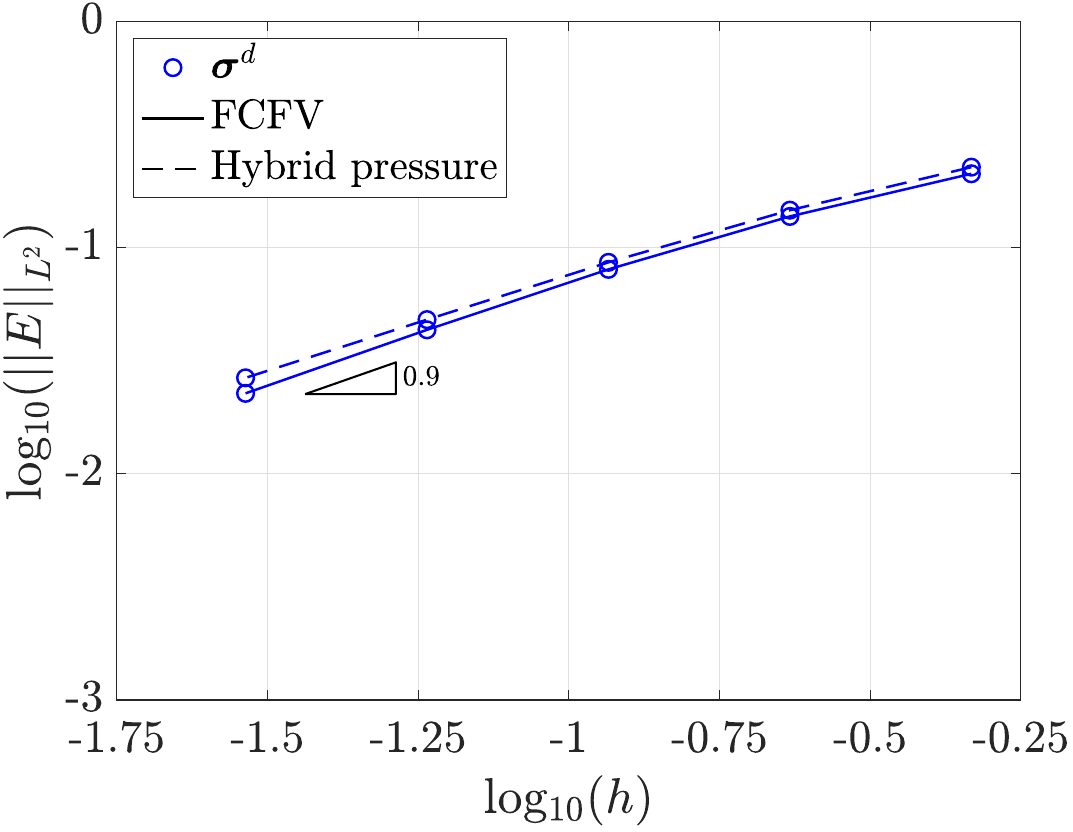}}	
	
	\caption{Couette Navier-Stokes flow - Mesh convergence of the error of (a,d) velocity and hybrid velocity, (b,d) pressure and hybrid pressure, and (c,f) deviatoric stress tensor, measured in the $\eltwo$ norm as a function of the cell size $h$ on meshes of triangular cells. Top row: regular meshes. Bottom row: distorted meshes.}
	\label{fig:convNStri}
\end{figure}

Figures~\ref{fig:convNSqua} and~\ref{fig:convNStri} display the optimal convergence of the error for cell and face velocity (left), cell and face pressure (centre), and deviatoric stress tensor (right), upon mesh refinement, for quadrilateral and triangular cells. 
The behaviour of the $\eltwo$ norm of the relative error confirms the results previously observed in the case of the Stokes flow, with all variables optimally converging with order one, independently of cell type or distortion.
It is worth noticing that, given the extremely low value of the Reynolds number,  the flow under analysis is similar to a Stokes flow and no significant accuracy gain is expected in this problem. Indeed,  the comparison of the hybrid pressure formulation with the FCFV method in~\cite{Vieira-VGSH-24} shows that the approximations of velocity and deviatoric stress tensor are almost superimposed. Concerning the discretisation of $p$,  the hybrid pressure formulation shows slightly increased robustness, achieving approximately the same accuracy using regular and distorted meshes. On the contrary, the FCFV method attains optimal convergence for pressure in both types of grid but it suffers a slight loss of accuracy due to the distortion of the cells (see Figures~\ref{fig:convNSquaPreg} and~\ref{fig:convNSquaPdist}).

\begin{figure}[!htb]
	\centering
	\subfigure[Quadrilaterals - Mesh 1]{\includegraphics[width=0.23\textwidth]{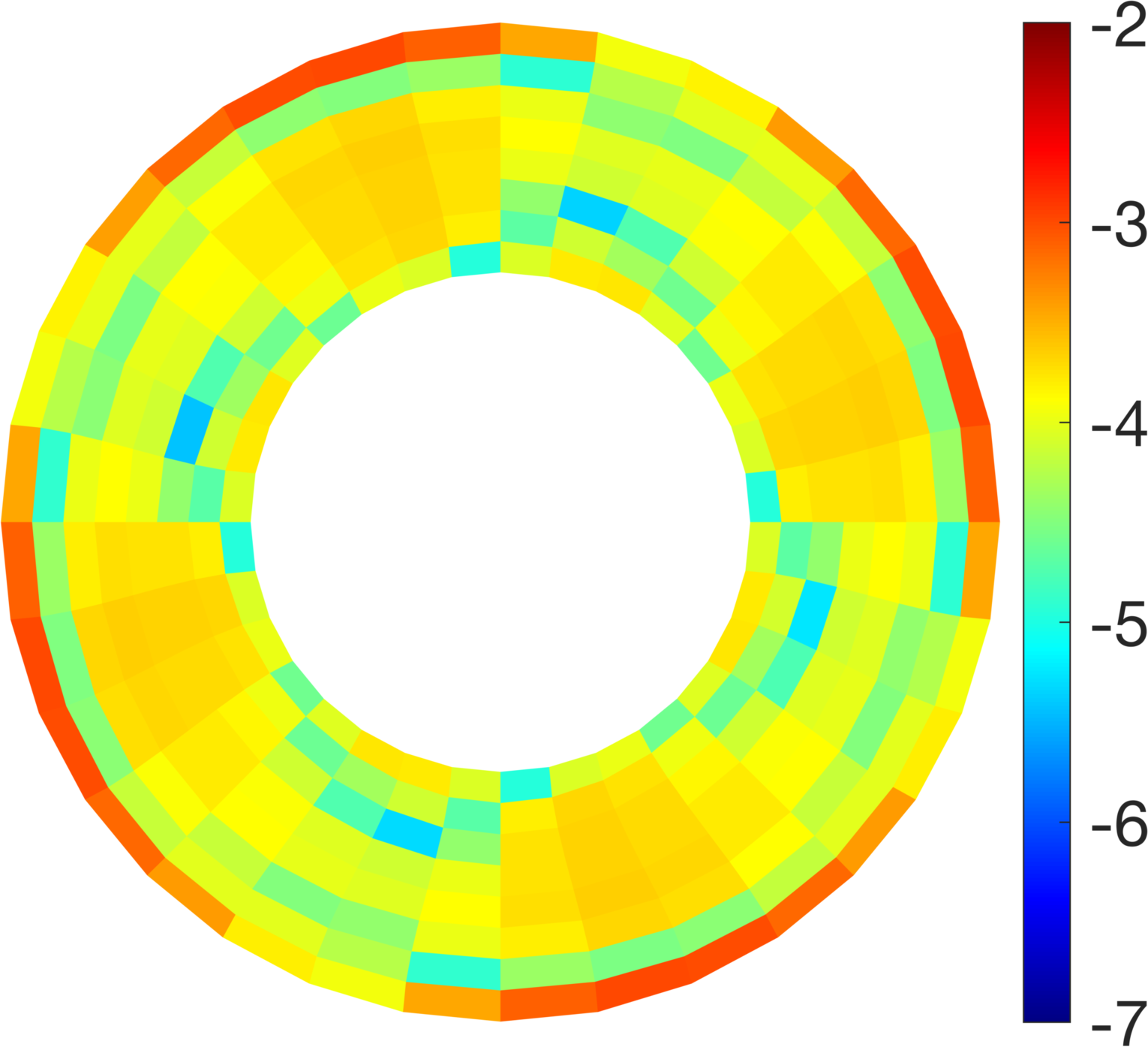}}
	\hspace{5pt}	
	\subfigure[Quadrilaterals - Mesh 4]{\includegraphics[width=0.23\textwidth]{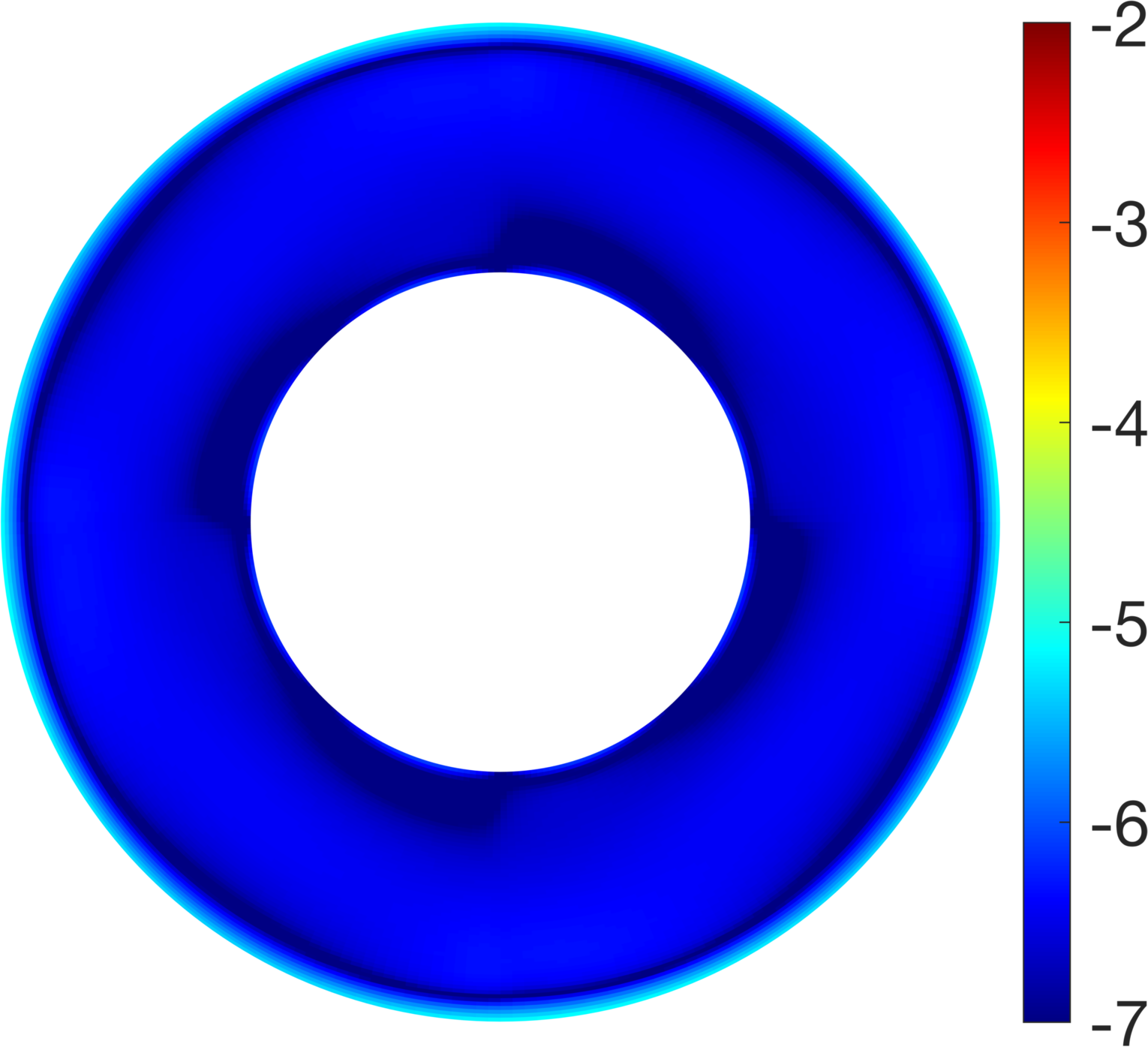}}
	\hspace{5pt}	
	\subfigure[Triangles - Mesh 1]{\includegraphics[width=0.23\textwidth]{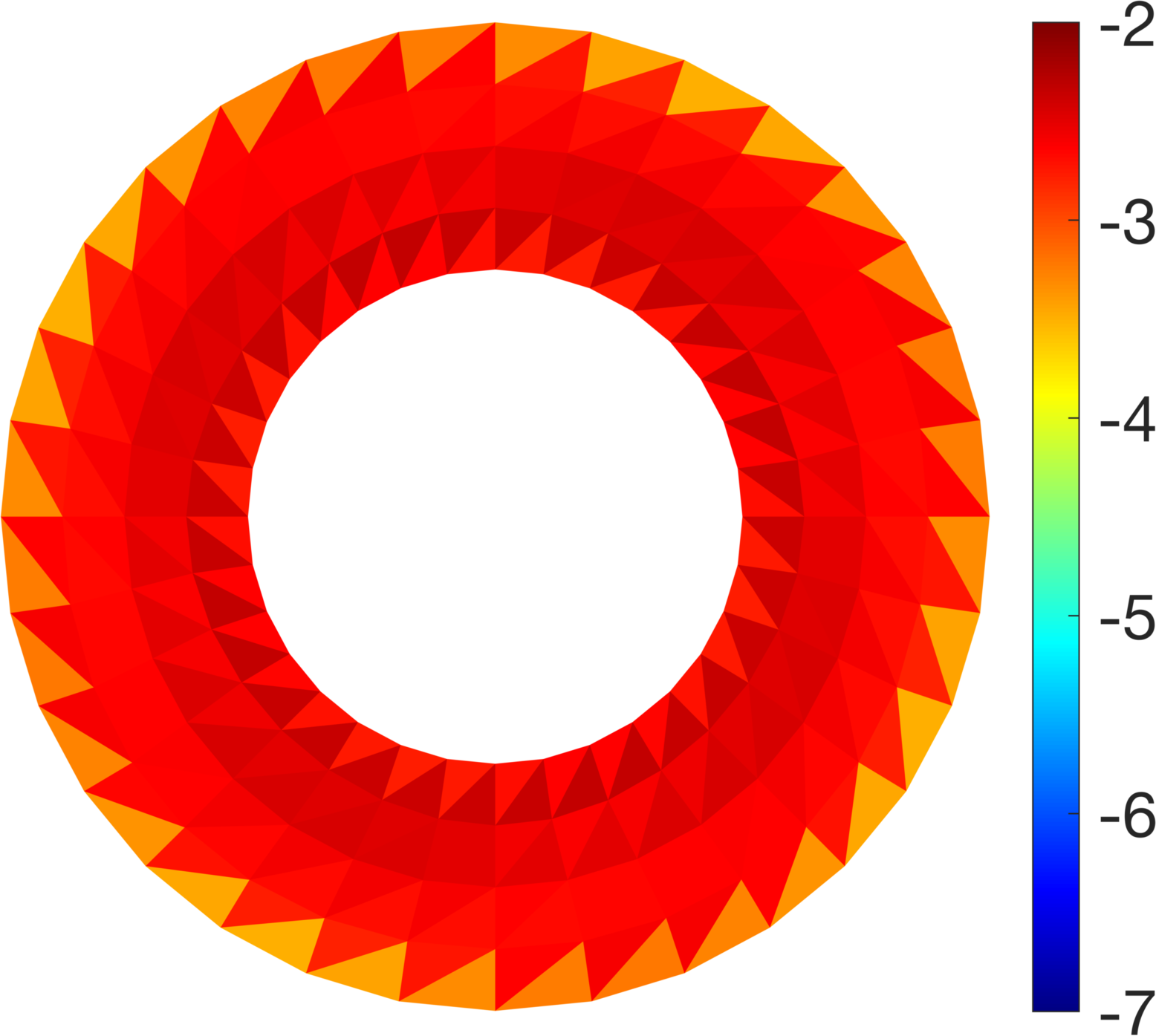}}
	\hspace{5pt}	
	\subfigure[Triangles - Mesh 4]{\includegraphics[width=0.23\textwidth]{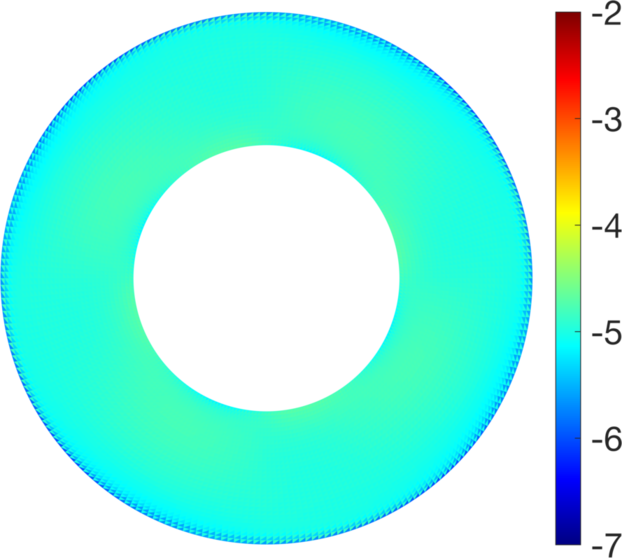}}	
		
	\subfigure[Quadrilaterals - Mesh 1]{\includegraphics[width=0.23\textwidth]{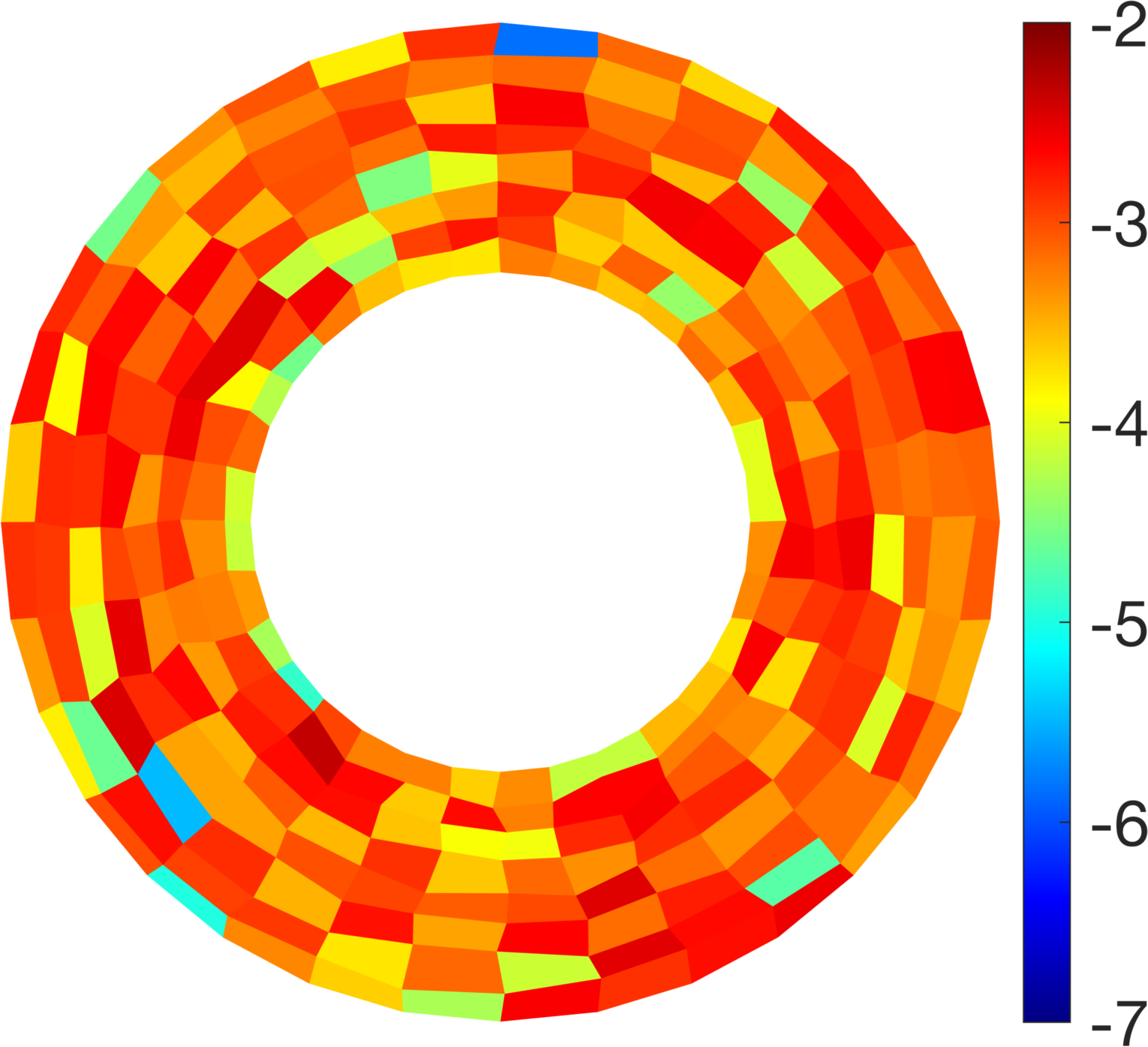}}
	\hspace{5pt}	
	\subfigure[Quadrilaterals - Mesh 4]{\includegraphics[width=0.23\textwidth]{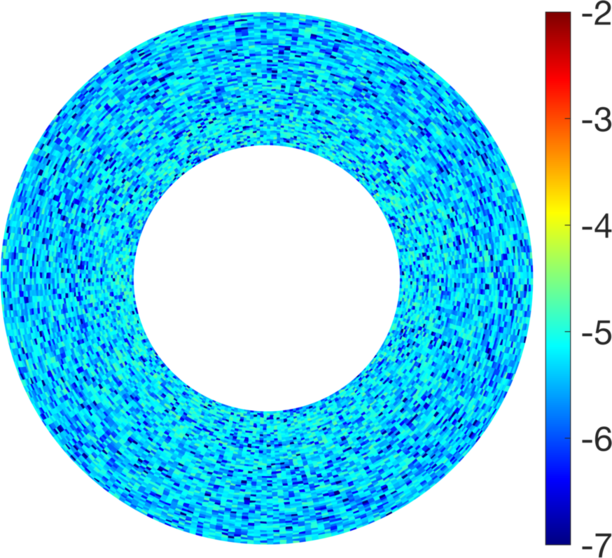}}
	\hspace{5pt}	
	\subfigure[Triangles - Mesh 1]{\includegraphics[width=0.23\textwidth]{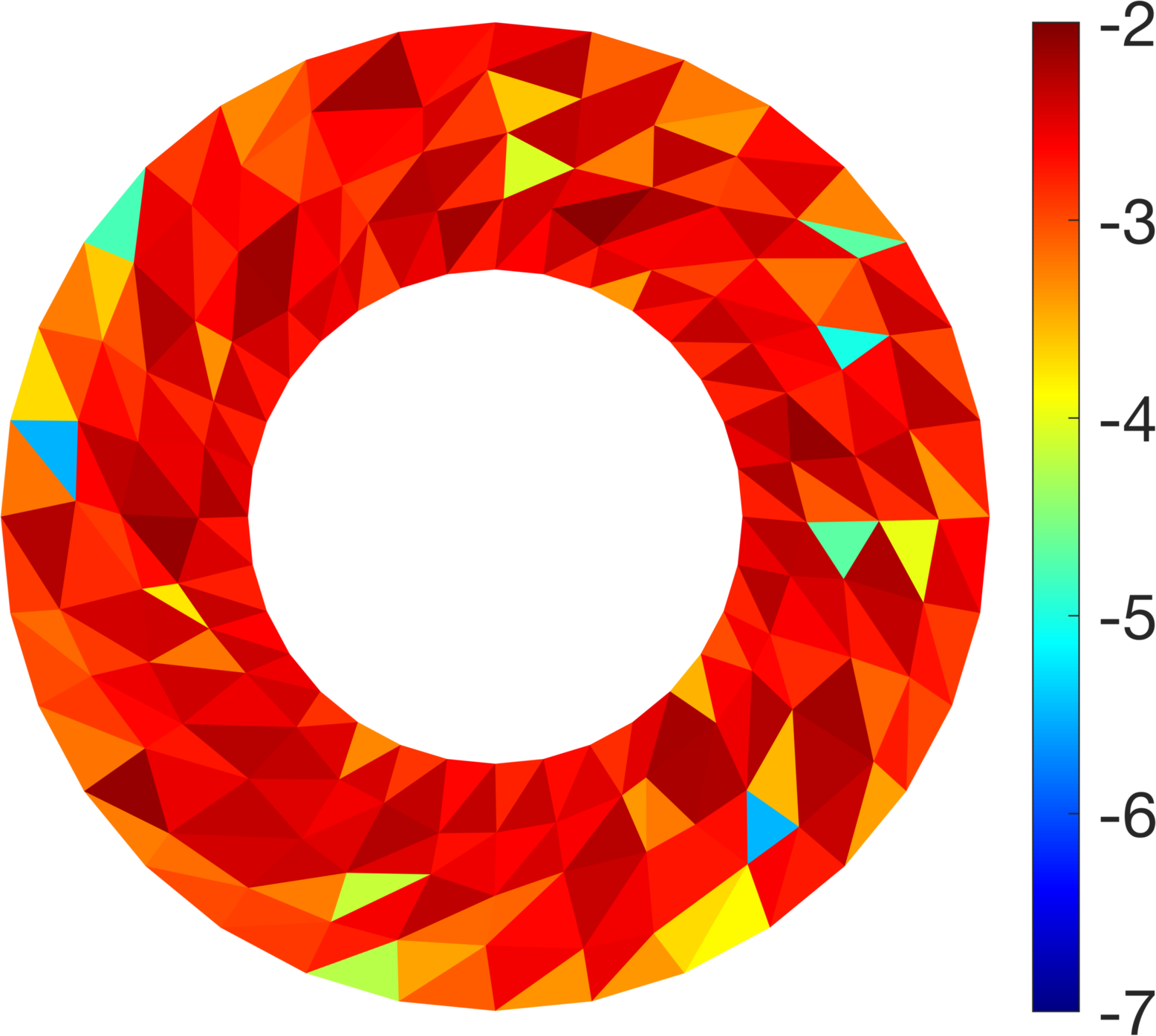}}
	\hspace{5pt}	
	\subfigure[Triangles - Mesh 4]{\includegraphics[width=0.23\textwidth]{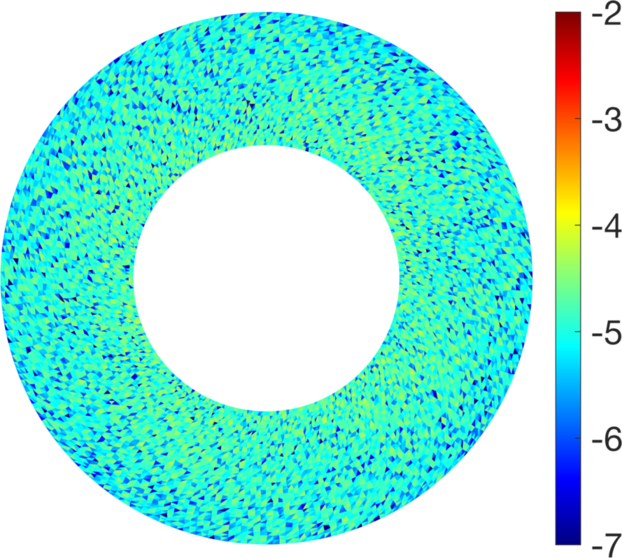}}
	
	\caption{Couette Navier-Stokes flow - Cell mass flux $\log_{10}|\fluxEl|$ for the first and fourth levels of mesh refinement using quadrilateral and triangular cells. Top row: regular meshes. Bottom row: distorted meshes.}
	\label{fig:divNS}
\end{figure}

Concerning mass conservation,  Figure~\ref{fig:divNS} reports the absolute value of the cell mass flux in logarithmic scale,  that is, $\log_{10}|\fluxEl|$.
The results show that the maximum error in the cell-by-cell mass conservation ranges from $10^{-3}$ to $0.4 \times 10^{-2}$ for the first mesh of quadrilateral and triangular cells, respectively. Note that this is significantly below the corresponding best approximation error than can be achieved with these meshes, being $h \simeq 0.4$ for both quadrilaterals and triangles.
Moreover,  this error consistently decreases when the mesh is refined, achieving values of the order of $10^{-6}$ on the fourth mesh, while the corresponding mesh size is of order $0.2 \times 10^{-1}$ for both meshes.
It is worth noticing that the distortion of the cells (bottom row of the figure) affects only marginally the error in the cell-wise incompressibility, increasing the maximum error by a moderate factor,  always significantly below the best approximation error of the mesh.
Finally,  the total mass flux obtained summing the contributions $\fluxEl$ for all mesh cells achieves machine precision independently of the level of mesh refinement, thus confirming the ability of the hybrid pressure formulation to globally preserve the mass in the domain.

\subsection{Navier-Stokes flow in a lid-driven cavity}
\label{sc:CavityNS}

The third example presents the flow in a lid-driven cavity at $Re=1,000$ and $Re=3,200$.
Given the domain $\Omega = [0,1]^2$, a constant horizontal velocity field $\bu_D = [1,0]^T$ is imposed on the top lid, whereas no-slip conditions are enforced on the remaining portions of the boundary.
Given the height of the domain as characteristic length and the magnitude of $\bu_D$ as characteristic velocity, the viscosity $\nu$ is set to $1/Re$.

The goal of this test is to evaluate the influence of different choices of convective stabilisation on the hybrid pressure approximation and compare its performance with the FCFV approach in~\cite{Vieira-VGSH-24}.
More precisely, Lax-Friedrichs (LF) and Harten-Lax-van Leer (HLL) Riemann solvers, see~\eqref{eq:stabConv},  are employed.
The nonlinear problem is solved using the Newton-Raphson algorithm, with a tolerance of $10^{-10}$. 

\begin{figure}[!htb]
	\centering
	\subfigure[]{\includegraphics[width=0.3\textwidth]{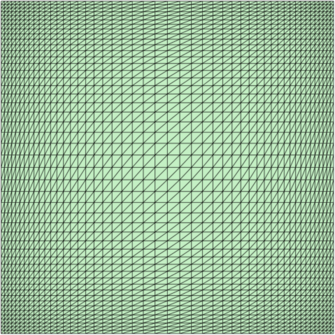}}
	\hspace{5pt}
	\subfigure[]{\includegraphics[width=0.3\textwidth]{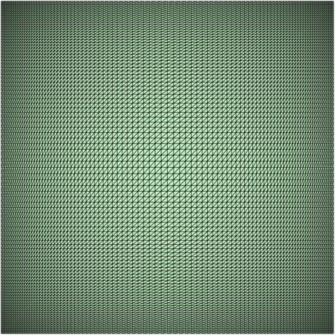}}
	
	\caption{Cavity flow - (a) First and (b) second level of mesh refinement.}
	\label{fig:meshCavityNS}
\end{figure}

The outcome of the hybrid pressure formulation is compared to the FCFV method in~\cite{Vieira-VGSH-24}, a Taylor-Hood (Q2Q1) finite element result with streamline-upwind Petrov-Galerkin (SUPG) stabilisation~\cite{Donea-Huerta-2003} and the reference solution by Ghia and co-workers~\cite{Ghia-82}.  
On the one hand, the computational meshes (see Figure~\ref{fig:meshCavityNS}) for the hybrid pressure and FCFV formulations are characterised by a local refinement near the boundaries, with the height of the first layer of boundary cells being $10^{-2}/i$ for the $i$-th level of mesh refinement and the growth ratio ranging from $1.06$ in the coarse mesh to $1.003$ in the finest mesh.  The resulting $i$-th mesh features $(24 \times 2^i) \times (24 \times 2^i) \times 2$ triangular cells. 
On the other hand, to compute the reference solution,  a mesh of $700 \times 700$ quadrilateral cells is employed. This mesh sets the height of the first cell to $0.96 \times 10^{-4}$ and features a local refinement near the four physical walls.

\begin{figure}[!htb]
	\centering
	\subfigure[]{\includegraphics[width=0.4\textwidth]{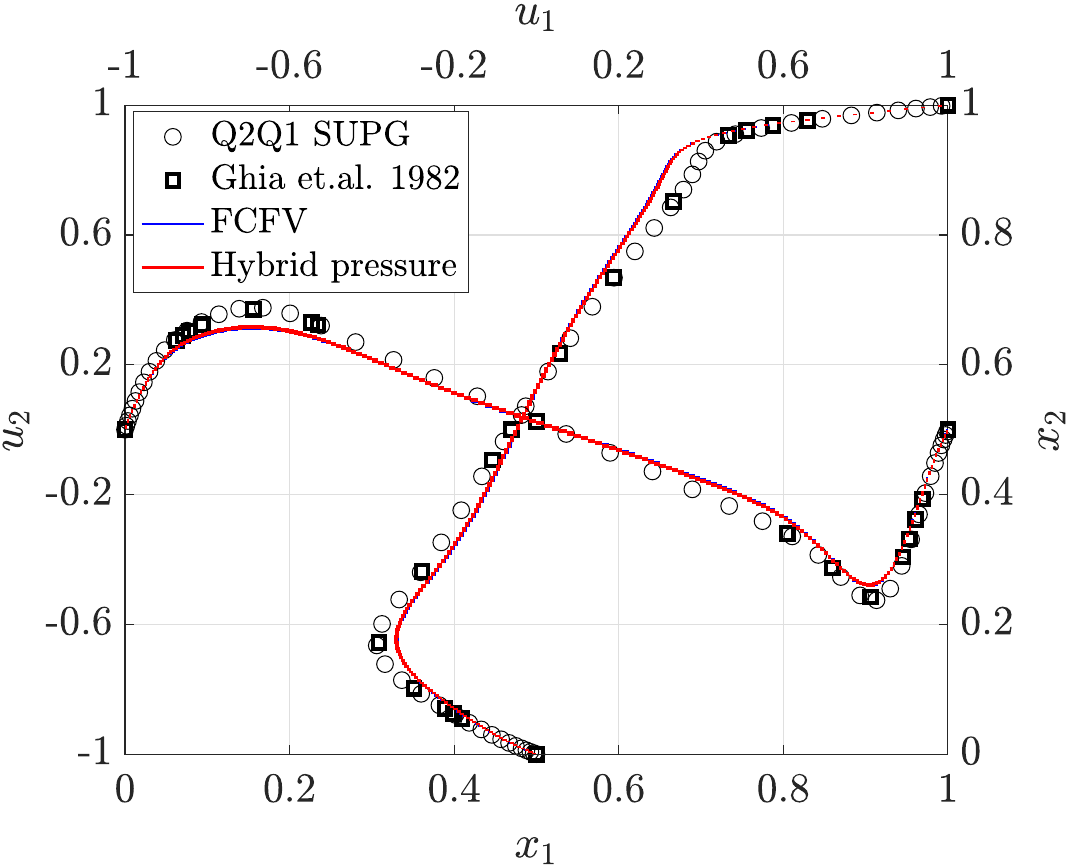}\label{fig:cavity1000Vel-LF3}}
	\hspace{5pt}
	\subfigure[]{\includegraphics[width=0.4\textwidth]{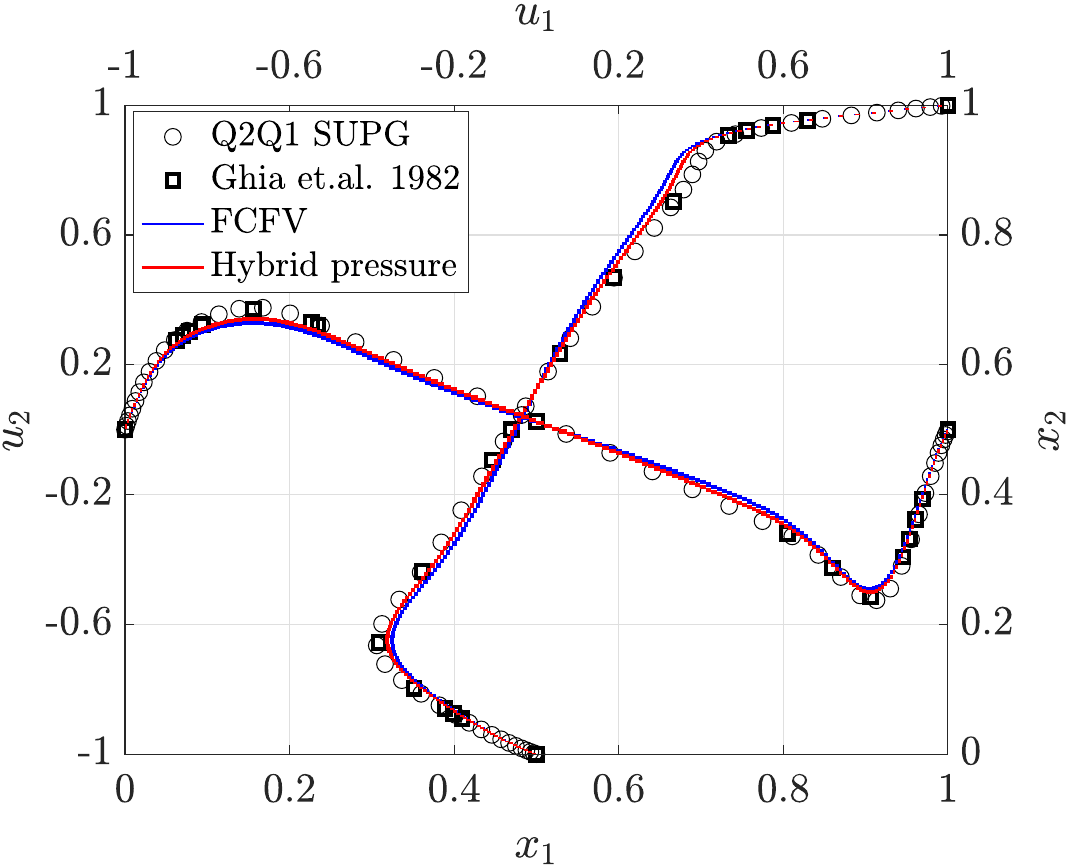}\label{fig:cavity1000Vel-HLL3}}
		
	\subfigure[]{\includegraphics[width=0.4\textwidth]{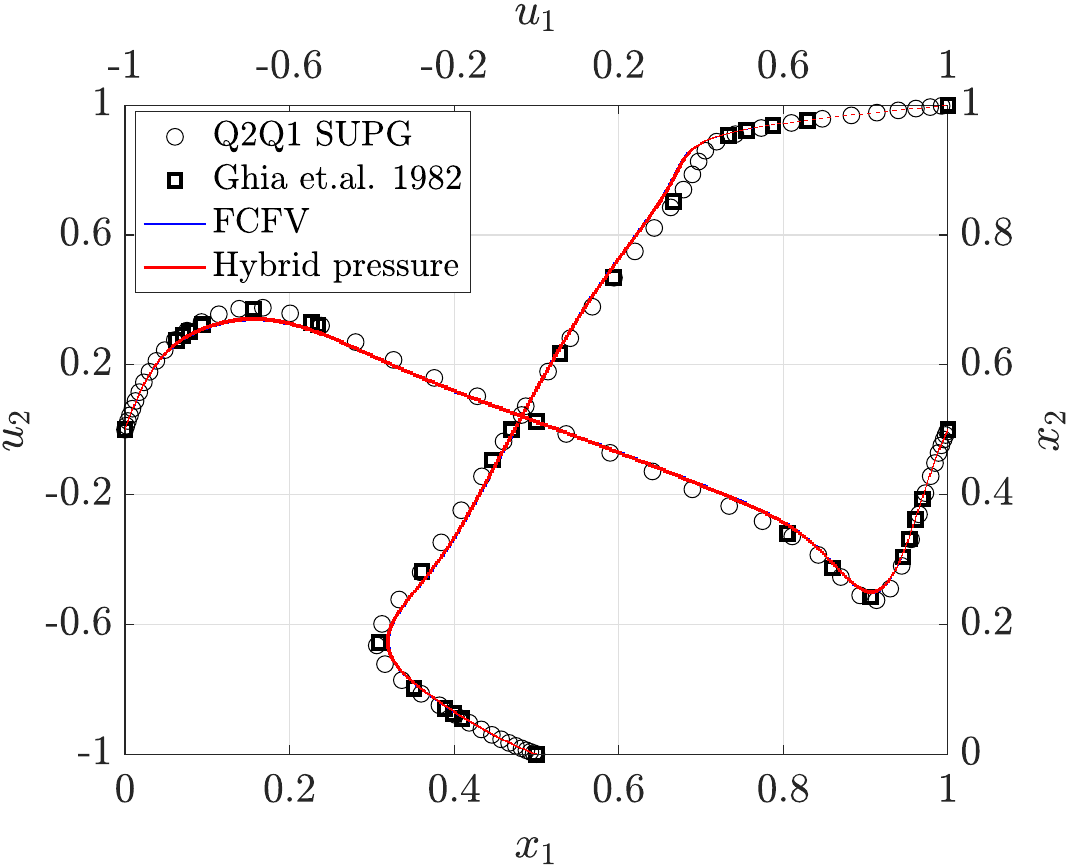}\label{fig:cavity1000Vel-LF4}}
	\hspace{5pt}
	\subfigure[]{\includegraphics[width=0.4\textwidth]{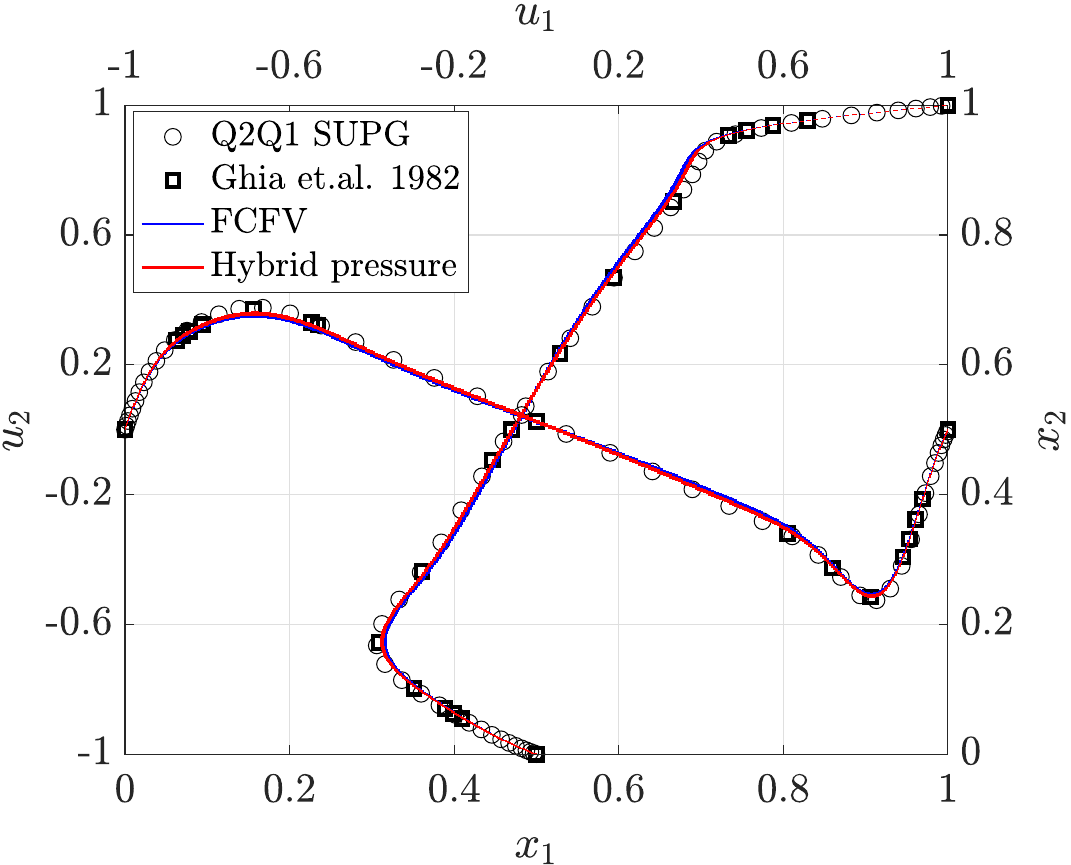}\label{fig:cavity1000Vel-HLL4}}
				
	\caption{Cavity flow at $Re=1,000$ - Profiles of velocity along the centrelines using (a,c) Lax-Friedrichs, and (b,d) Harten-Lax-van Leer Riemann solvers on different meshes. Top row: third mesh. Bottom row: fourth mesh.}
	\label{fig:cavity1000Vel}
\end{figure}

\begin{figure}[!htb]
	\centering
	\subfigure[]{\includegraphics[width=0.4\textwidth]{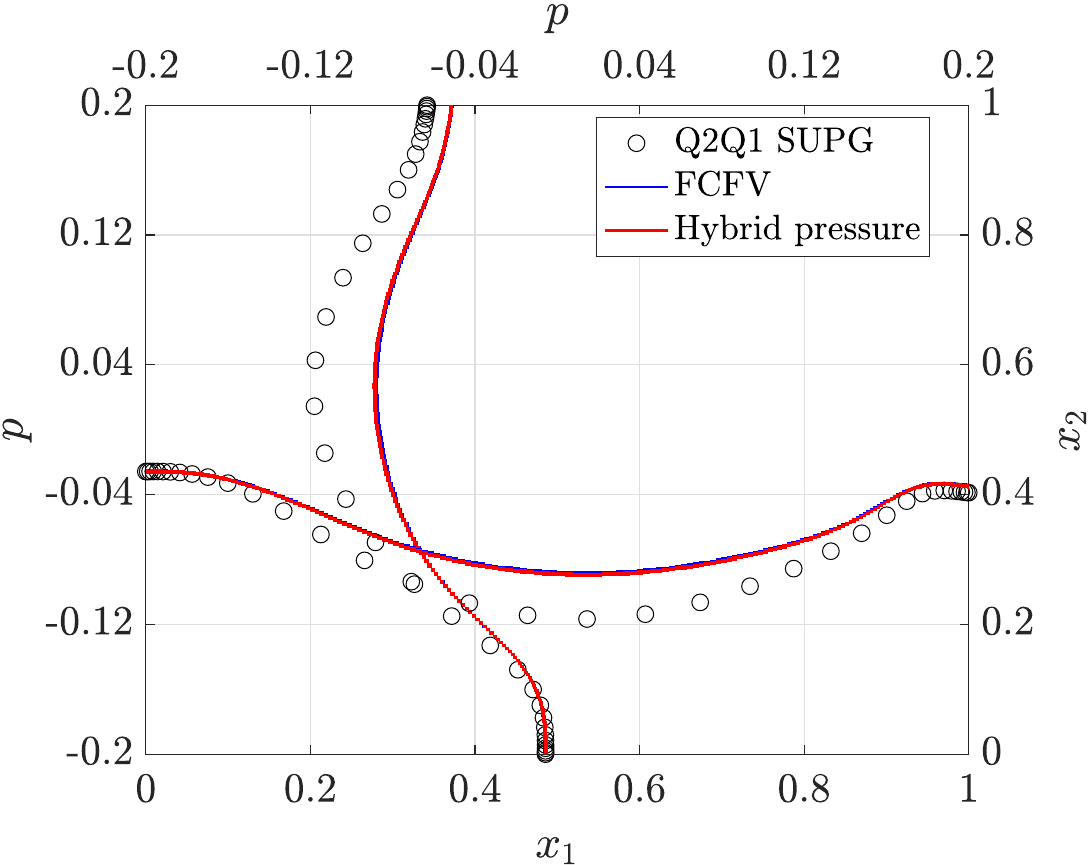}\label{fig:cavity1000Pres-LF3}}
	\hspace{5pt}
	\subfigure[]{\includegraphics[width=0.4\textwidth]{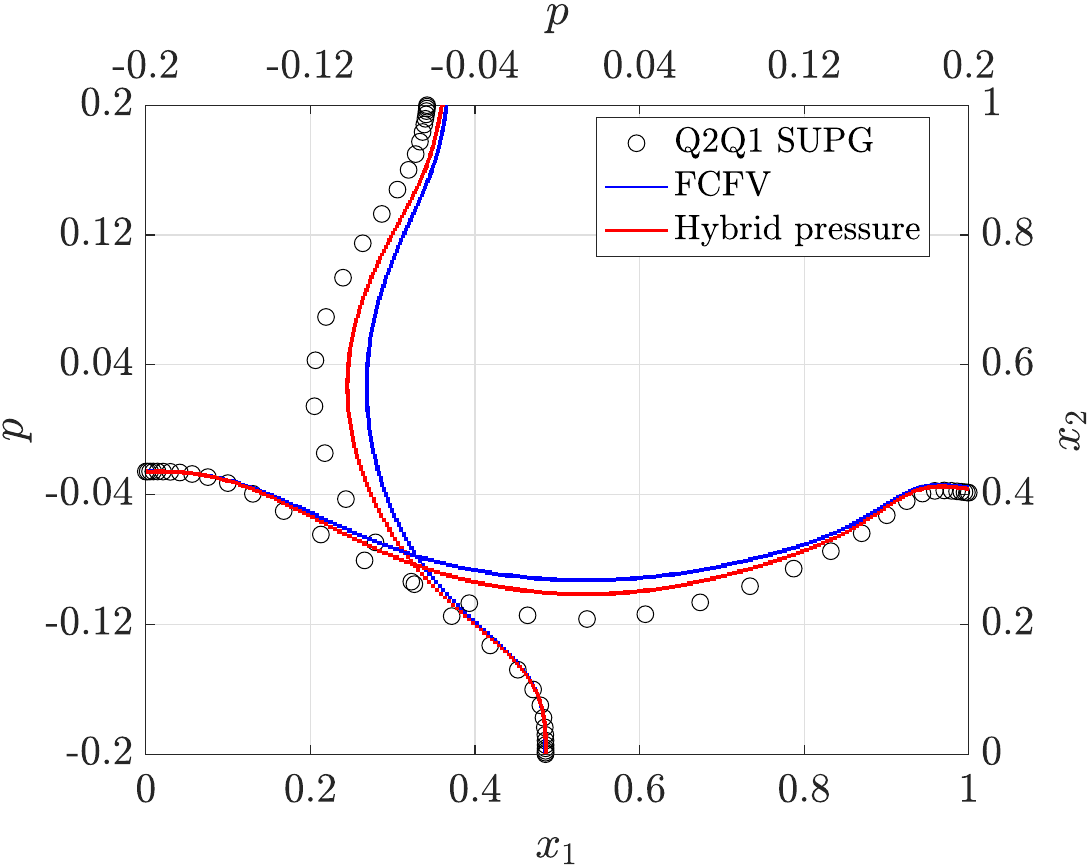}\label{fig:cavity1000Pres-HLL3}}
		
	\subfigure[]{\includegraphics[width=0.4\textwidth]{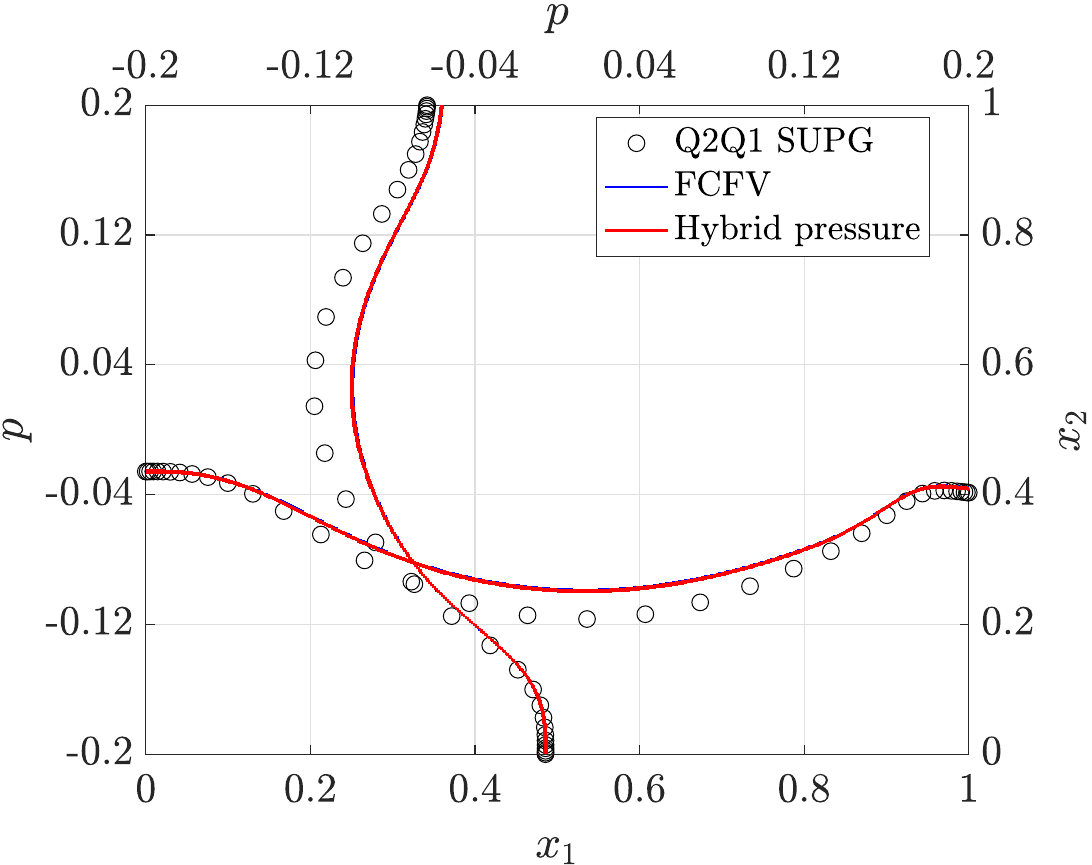}\label{fig:cavity1000Pres-LF4}}
	\hspace{5pt}
	\subfigure[]{\includegraphics[width=0.4\textwidth]{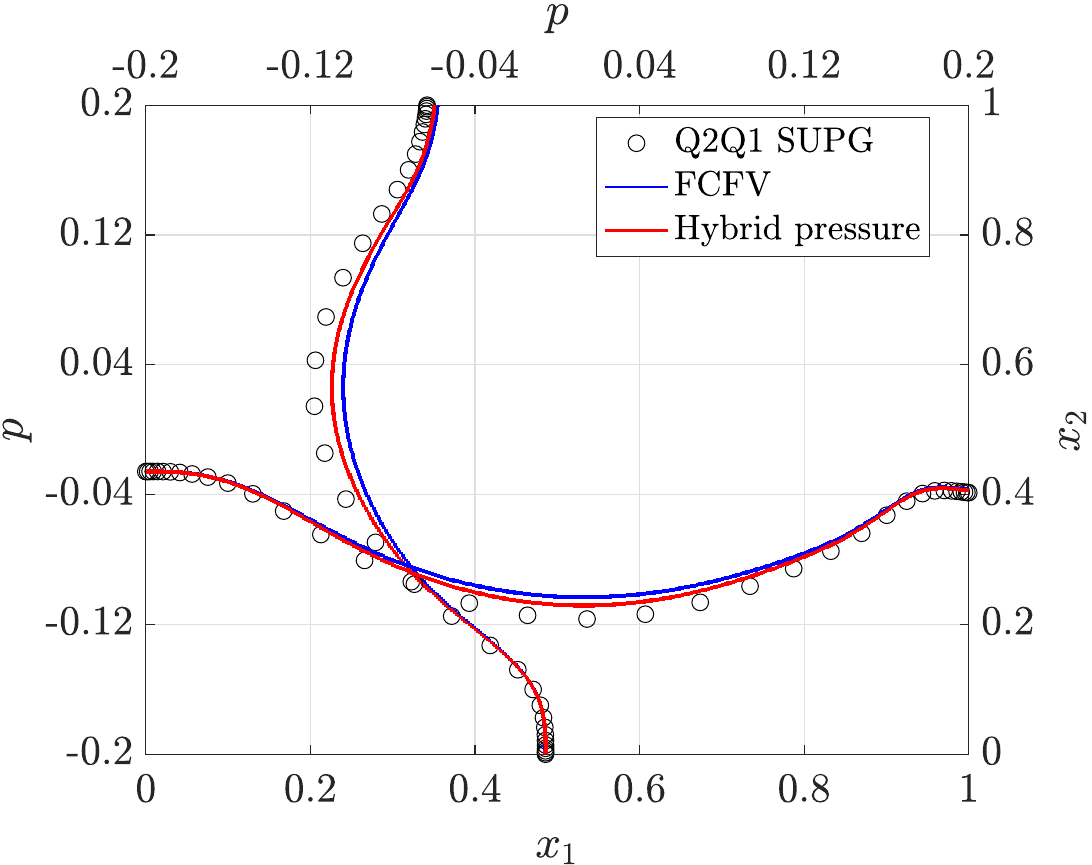}\label{fig:cavity1000Pres-HLL4}}	
	
	\caption{Cavity flow at $Re=1,000$ - Profiles of pressure along the centrelines using (a,c) Lax-Friedrichs, and (b,d) Harten-Lax-van Leer Riemann solvers on different meshes. Top row: third mesh. Bottom row: fourth mesh.}
	\label{fig:cavity1000Pres}
\end{figure}

The quantities of interest in this problem are the profiles of velocity and pressure at the centrelines of the domain.
Velocity profiles for $Re=1,000$ are reported in Figure~\ref{fig:cavity1000Vel} for the third and fourth level of mesh refinement using LF and HLL stabilisations. For the sake of brevity, only meshes without distortion are presented.
Using LF Riemann solver, the results of hybrid pressure and FCFV formulations are perfectly superposed,  showing excellent agreement with the Taylor-Hood and Ghia \etal\ references (Figure~\ref{fig:cavity1000Vel-LF3}).  As expected, the accuracy furtherly increases when the mesh is refined (Figure~\ref{fig:cavity1000Vel-LF4}). Alternatively,  a more accurate approximation can be achieved on the coarse mesh by employing the HLL Riemann solver (Figure~\ref{fig:cavity1000Vel-HLL3}), confirming the superior performance of this stabilisation, already observed for the FCFV method in~\cite{Vieira-VGSH-24}. In addition,  it is worth noticing that the hybrid pressure formulation with HLL stabilisation also provides slightly more accurate results with respect to the FCFV approximation with HLL on both meshes.
Nonetheless, it is in the approximation of $p$ that the hybrid pressure formulation with HLL stabilisation clearly shows its superior performance with respect to the FCFV scheme and the LF stabilisation. Figure~\ref{fig:cavity1000Pres} displays the profiles of pressure at the centrelines for $Re=1,000$. Whilst using LF Riemann solver the hybrid pressure formulation provides results comparable to the FCFV method and the fourth mesh is required to achieve an accurate description of pressure (Figure~\ref{fig:cavity1000Pres-LF4}), the HLL stabilisation significantly improves the quality of the approximation, even on the coarser mesh (Figure~\ref{fig:cavity1000Pres-HLL3}). Moreover, the hybrid pressure formulation clearly outperforms the FCFV method using the HLL Riemann solver: using the third mesh, the hybrid pressure formulation with HLL stabilisation provides results with accuracy comparable to the ones achieved by the FCFV scheme with HLL on the fourth level of mesh refinement, see~Figure~\ref{fig:cavity1000Pres-HLL4}, while reducing the number of cells from $294,912$ to $73,728$. 

\begin{figure}[!htb]
	\centering
	\subfigure[]{\includegraphics[width=0.48\textwidth]{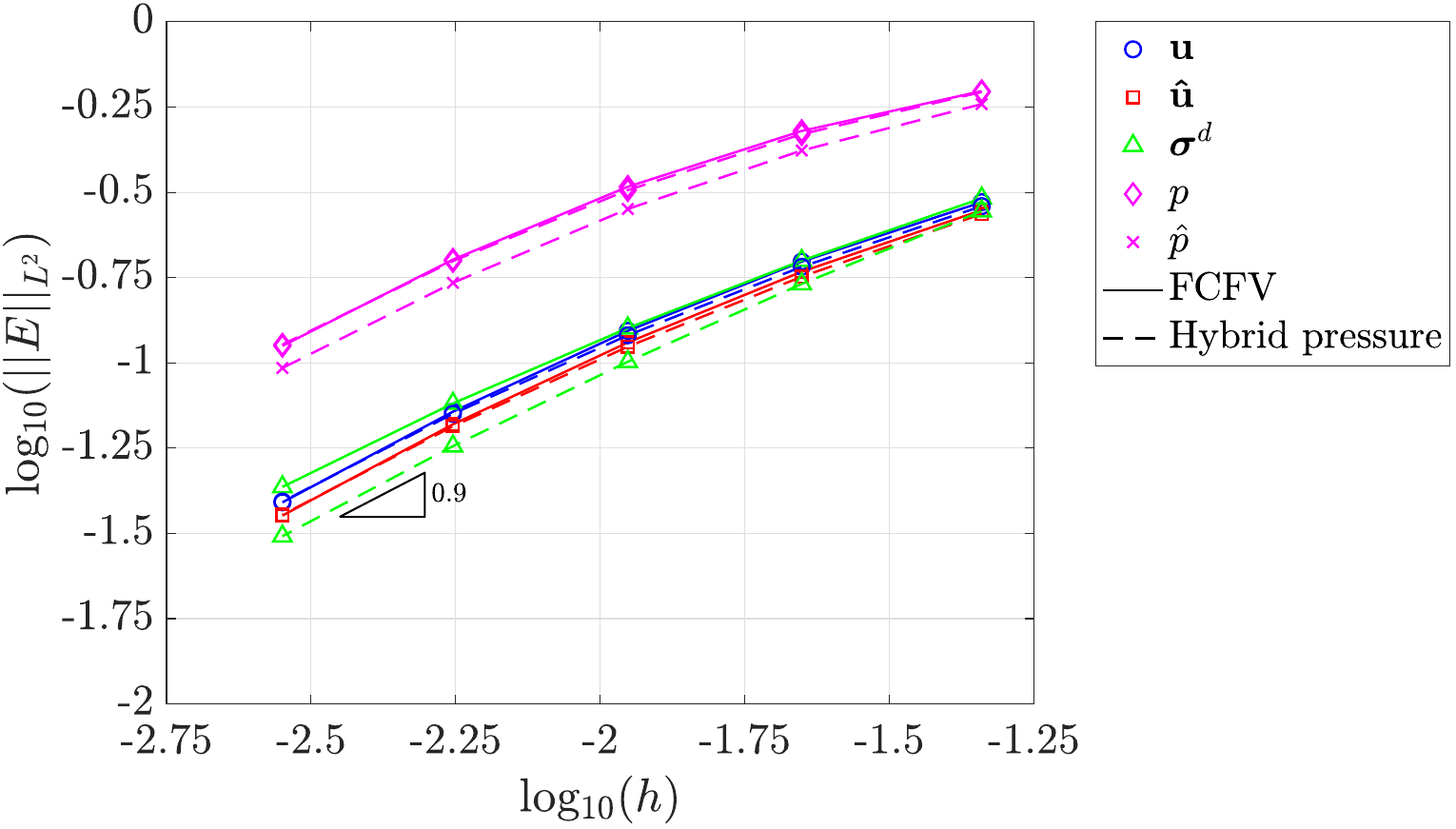}}
	\hspace{2pt}
	\subfigure[]{\includegraphics[width=0.48\textwidth]{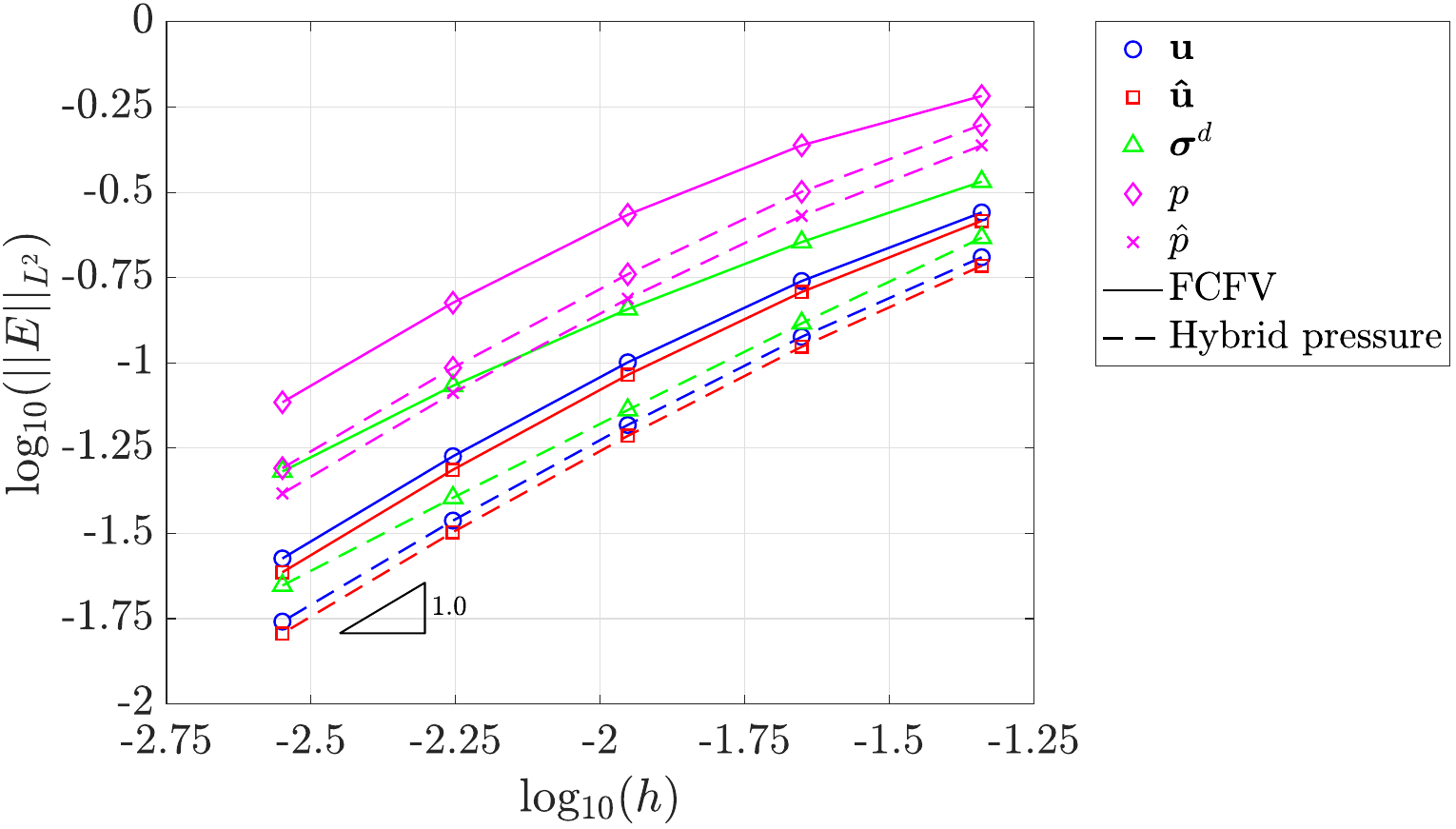}}
	
	\caption{Cavity flow at $Re=1,000$ - Mesh convergence of the error of velocity, hybrid velocity, pressure, hybrid pressure, and deviatoric stress tensor, measured in the $\eltwo$ norm as a function of the cell size $h$ using (a) Lax-Friedrichs and (b) Harten-Lax-van Leer Riemann solvers. }
	\label{fig:convCavity1000}
\end{figure}

It is worth noticing that the quantities of interest reported in Figures~\ref{fig:cavity1000Vel} and~\ref{fig:cavity1000Pres} can be significantly influenced by local modifications of the flow field. 
To evaluate the global accuracy of the method,  the errors of the hybrid pressure and FCFV formulations with respect to the Taylor-Hood solution are reported in Figure~\ref{fig:convCavity1000} as a function of the mesh size $h$.
The results confirm that whilst the errors achieved by the two methods are comparable using the LF stabilisation,  the hybrid pressure clearly outperforms the FCFV scheme in the approximation of all variables when the HLL Riemann solver is employed.
\begin{remark}
The reference Taylor-Hood solution used to compute the errors of the FCFV and hybrid pressure formulation is obtained by simulating the so-called \emph{leaky} cavity flow,  given the strong imposition of the incompatible Dirichlet boundary conditions on the top left and top right corners.
On the contrary, the FCFV and hybrid pressure formulations impose Dirichlet conditions at the barycentres of mesh faces, thus practically circumventing the incompatibility of these boundary conditions.
It is worth noticing that the region affected by the incompatible Dirichlet boundary conditions in the Taylor-Hood solution features multiple cells and is reduced upon mesh refinement. This also affects the size of the region where pressure is non-zero, which shrinks tending to a pointwise singularity when $h {\rightarrow} 0$~\cite{Donea-Huerta-2003}. In order to maintain a comparable reference Taylor-Hood solution for all meshes under analysis and to avoid introducing mesh-dependent artifacts in the error computation,  the regions $[0,0.05]\times[0.95,1]$ and $[0.95,1]\times[0.95,1]$ are excluded from the domain $\Omega$ when computing the $\eltwo(\Omega)$ norm of the error~\cite{Vieira-VGSH-24}, as reported in Figure~\ref{fig:convCavity1000}.
\end{remark}

\begin{figure}[!htb]
	\centering
	\subfigure[]{\includegraphics[width=0.4\textwidth]{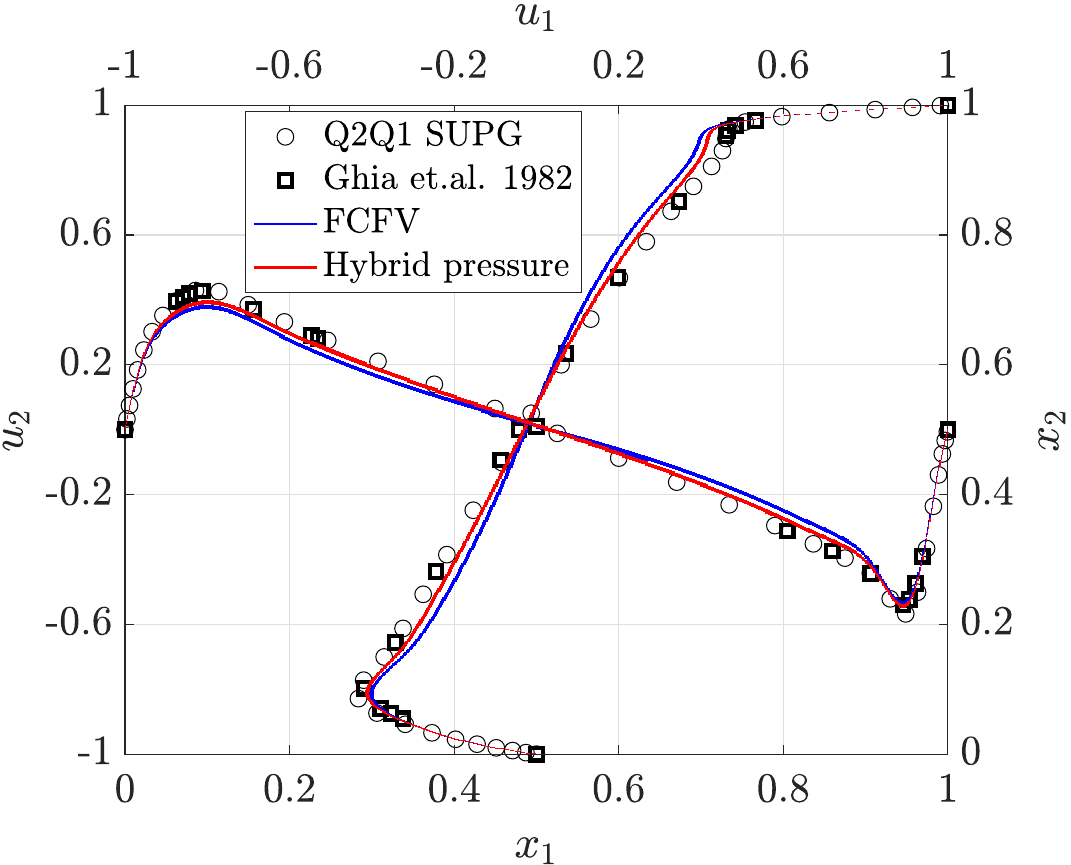}\label{fig:cavity3200Vel}}
	\hspace{5pt}
	\subfigure[]{\includegraphics[width=0.4\textwidth]{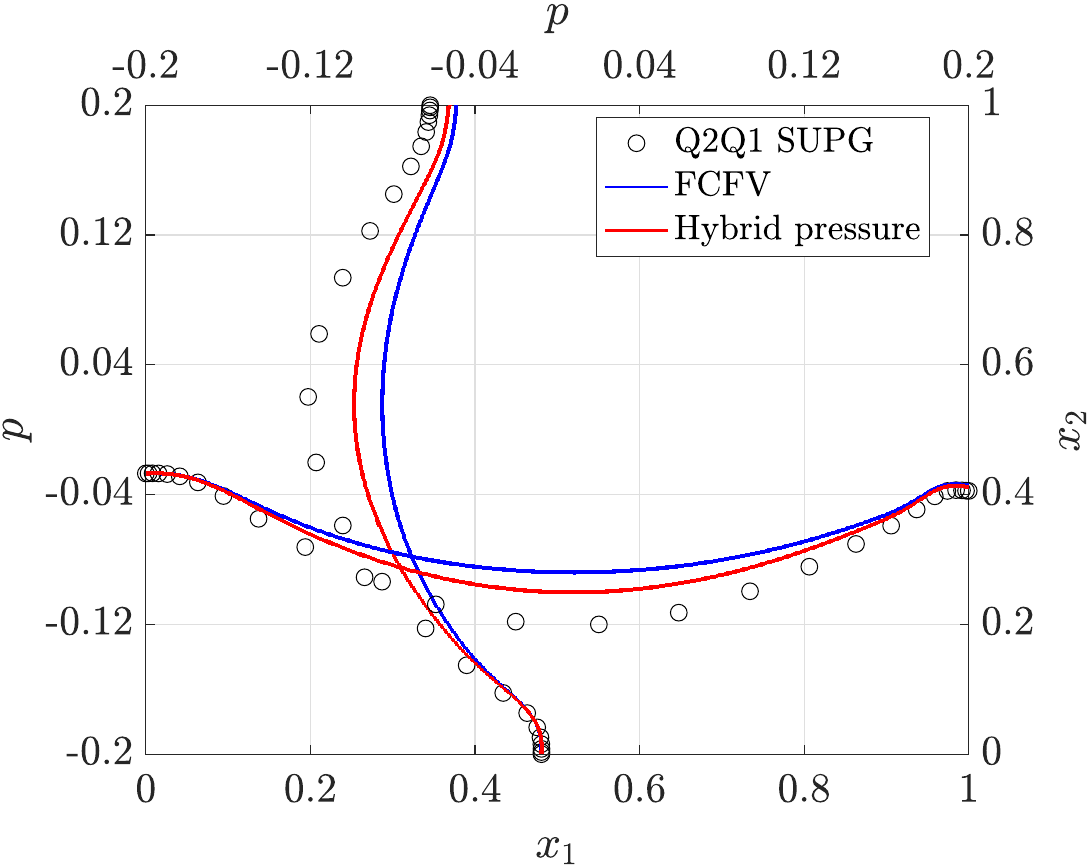}\label{fig:cavity3200Pres}}
	
	\subfigure[]{\includegraphics[width=0.48\textwidth]{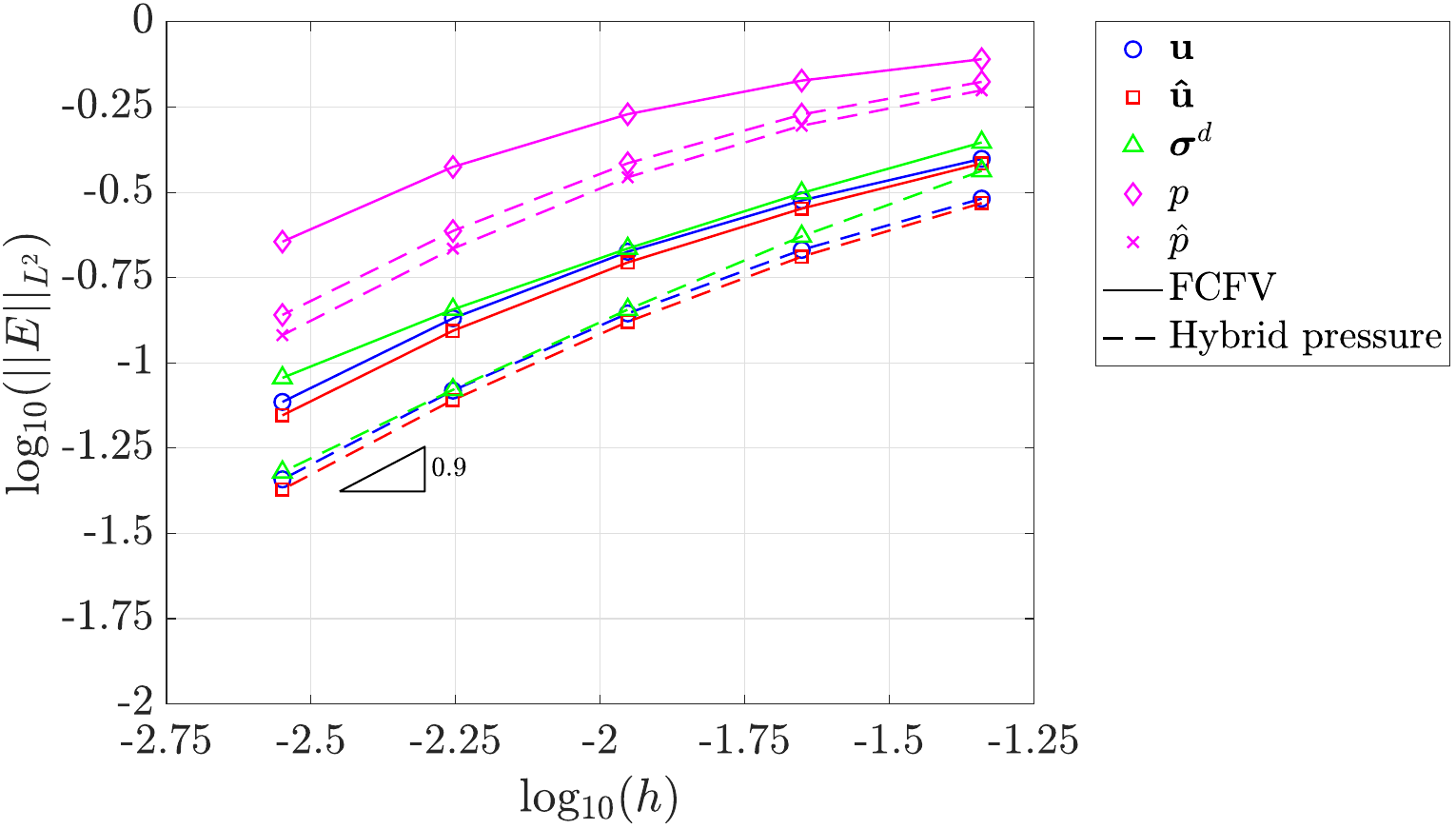}\label{fig:cavity3200Conv}}

	\caption{Cavity flow at $Re=3,200$ - Approximation computed using Harten-Lax-van Leer Riemann solver. Profiles of (a) velocity and (b) pressure along the centrelines of the fourth mesh. (c) Mesh convergence of the error of velocity, hybrid velocity, pressure, hybrid pressure, and deviatoric stress tensor, measured in the $\eltwo$ norm as a function of the cell size $h$.}
	\label{fig:cavity3200}
\end{figure}

Finally,  Figure~\ref{fig:cavity3200} presents the approximation of the cavity flow for Reynolds $3,200$. 
Given the previously observed superiority of HLL with respect to LF,  results are reported only for the HLL Riemann solver using the hybrid pressure formulation and the FCFV method.  
The profiles of velocity and pressure computed on the fourth mesh are displayed along the centrelines. in Figure~\ref{fig:cavity3200Vel} and~\ref{fig:cavity3200Pres}, respectively.
Also in this case, the hybrid pressure formulation outperforms the FCFV approach achieving a slightly more accurate approximation of velocity, while significantly improving the description of pressure. 

The convergence of the relative $\eltwo$ error, reported in Figure~\ref{fig:cavity3200Conv} for all the variables, confirms the superiority of the hybrid pressure formulation to construct an accurate face-based FV paradigm to simulate viscous laminar incompressible flows. 
Moreover, the accuracy appears to slightly increase at higher values of the Reynolds number, when convective effects are more relevant. Further studies, beyond the scope of the present work, need to be performed to evaluate the behaviour of the proposed formulation in the turbulent regime. 
The closely related FCFV method was recently shown to be suitable for simulating turbulent incompressible flows using Reynolds-averaged Navier-Stokes equations with Spalart-Allmaras turbulence model~\cite{Vieira-VGSH-24}.

\subsection{Three-dimensional Navier-Stokes flow past a sphere}
\label{sc:SphereNS}

The last example considers the steady-state incompressible Navier-Stokes flow past a three-dimensional sphere of diameter $D=1$ at $Re=20$. 
The computational domain is defined as $\Omega = \left( [-10D,20D] \times [-10D,10D] \times [-10D,10D] \right)\setminus \mathcal{B}_{D/2,\mathbf{0}}$. 

\begin{figure}[!htb]
	\centering
	\subfigure[]{\label{fig:meshNS}\includegraphics[width=0.45\textwidth]{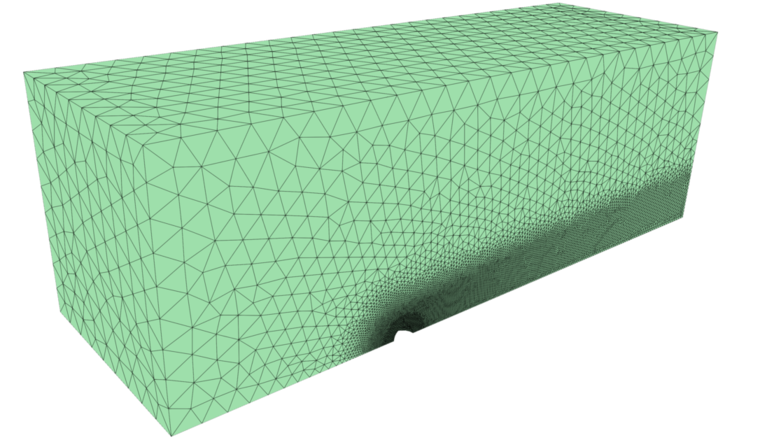}}
	\hspace{5pt}
	\subfigure[Zoom of~\ref{fig:meshNS}]{\includegraphics[width=0.4\textwidth]{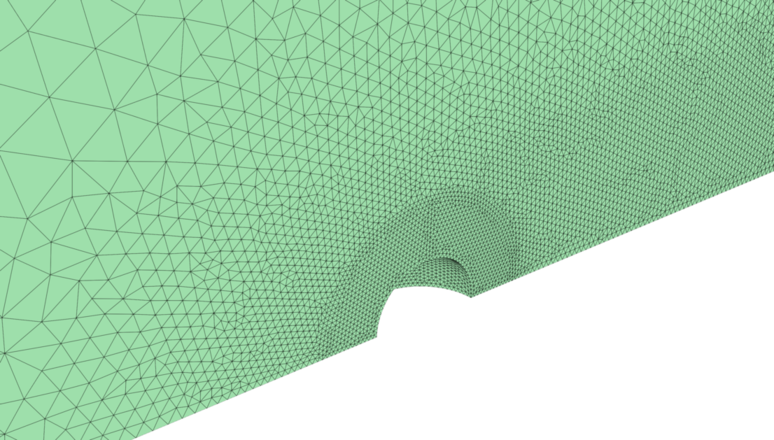}}
		
	\subfigure[]{\label{fig:uNS}\includegraphics[width=0.45\textwidth]{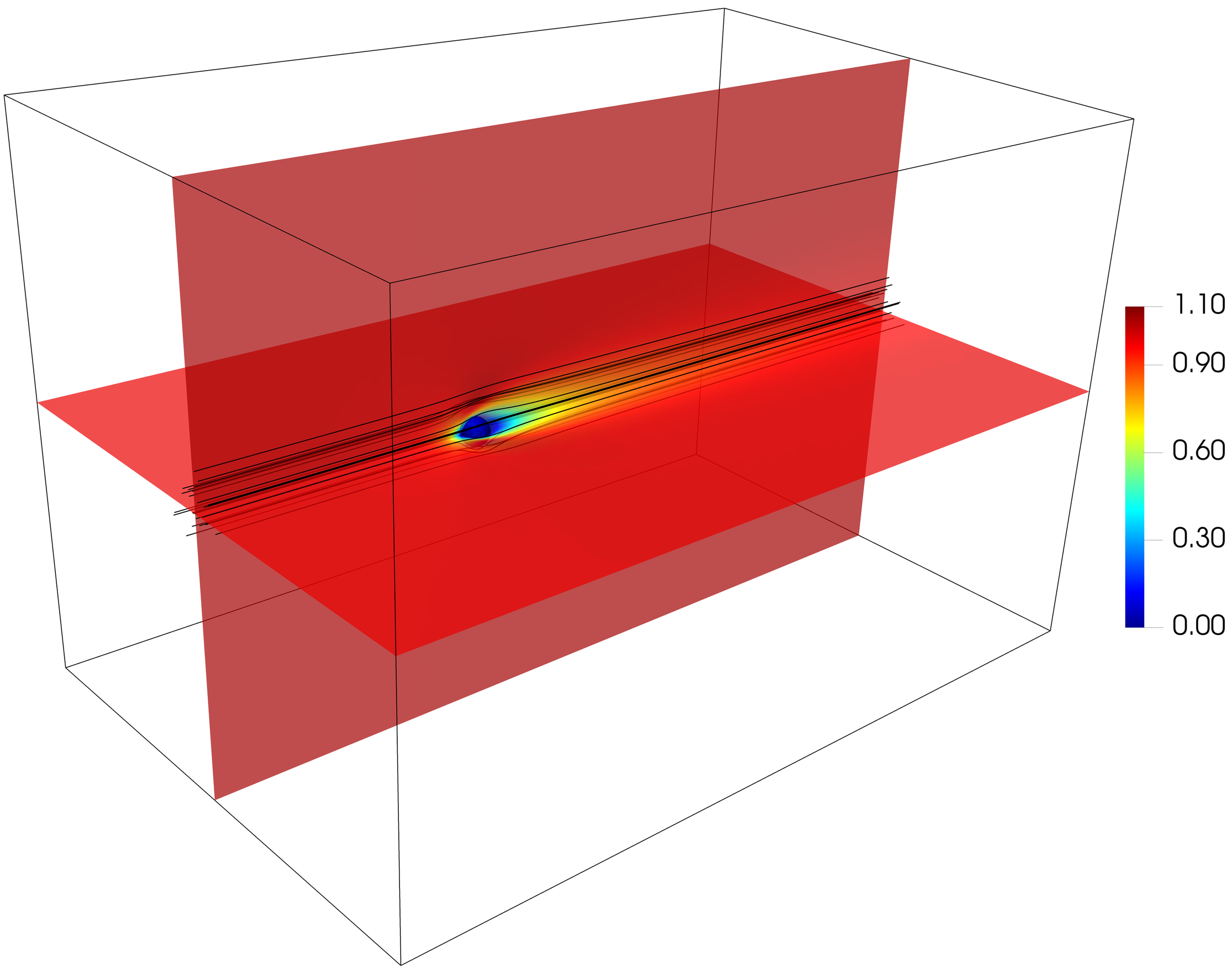}}
	\hspace{5pt}
	\subfigure[Zoom of~\ref{fig:uNS}]{\includegraphics[width=0.4\textwidth]{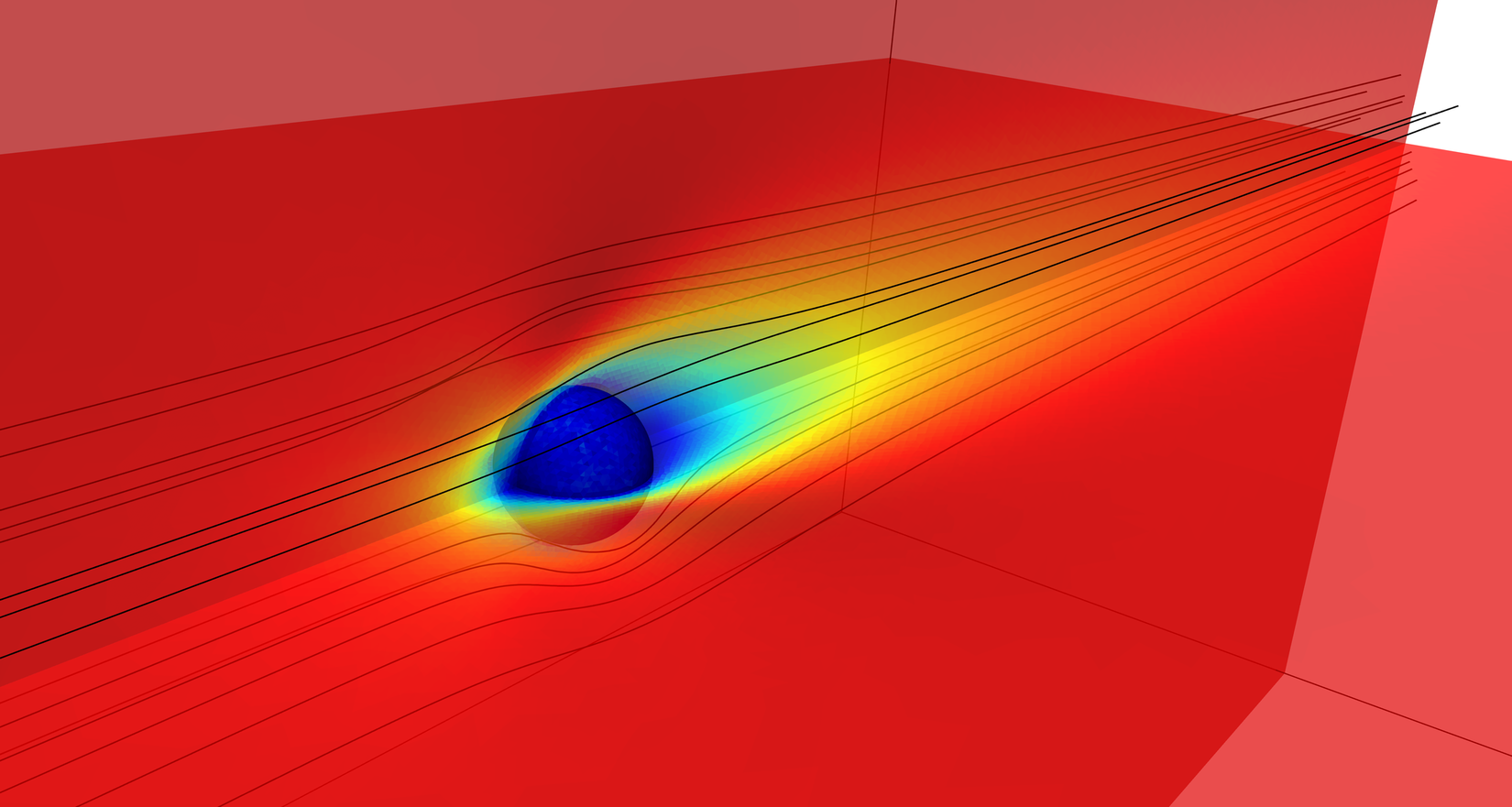}}
		
	\subfigure[]{\label{fig:pNS}\includegraphics[width=0.45\textwidth]{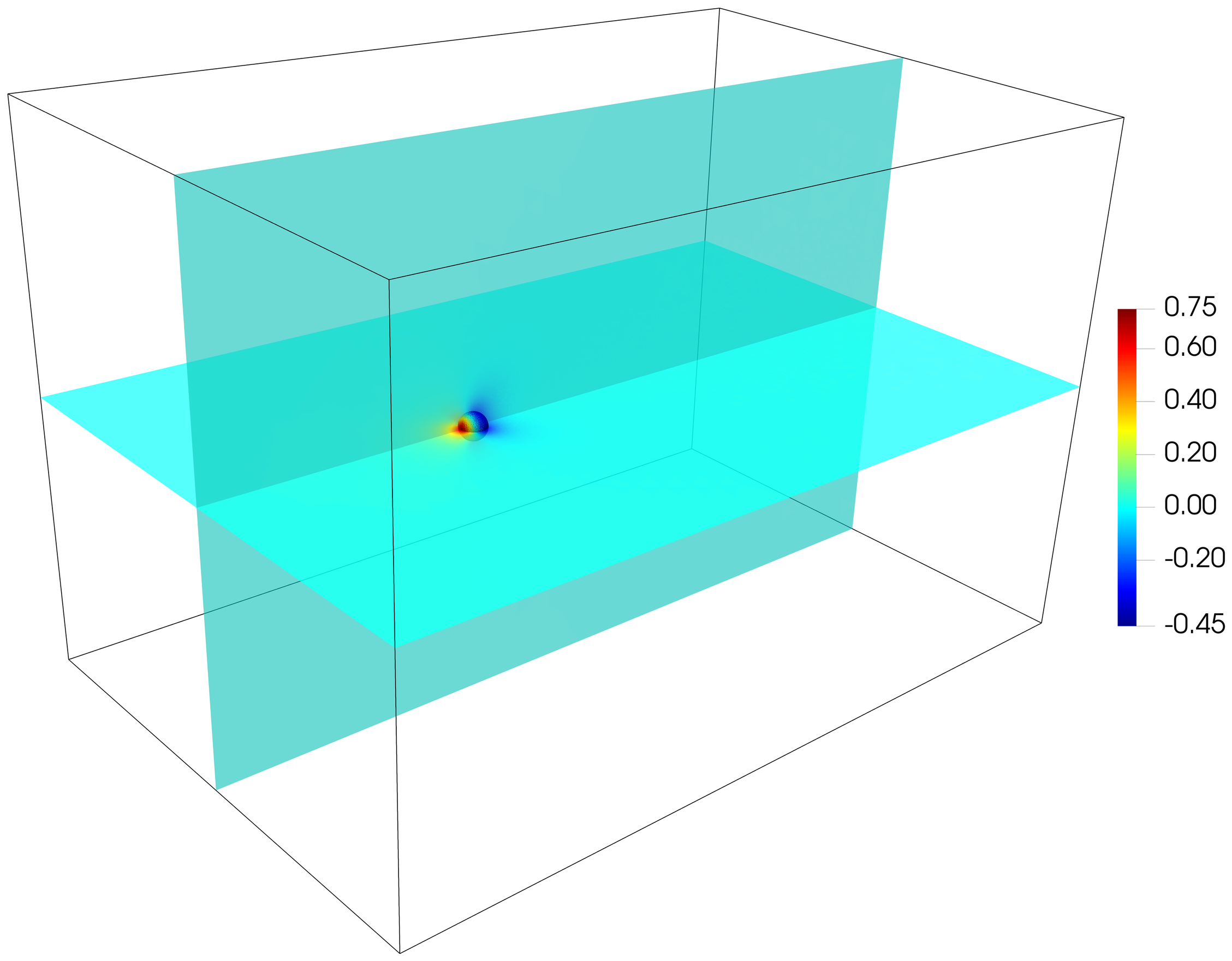}}
	\hspace{5pt}
	\subfigure[Zoom of~\ref{fig:pNS}]{\includegraphics[width=0.4\textwidth]{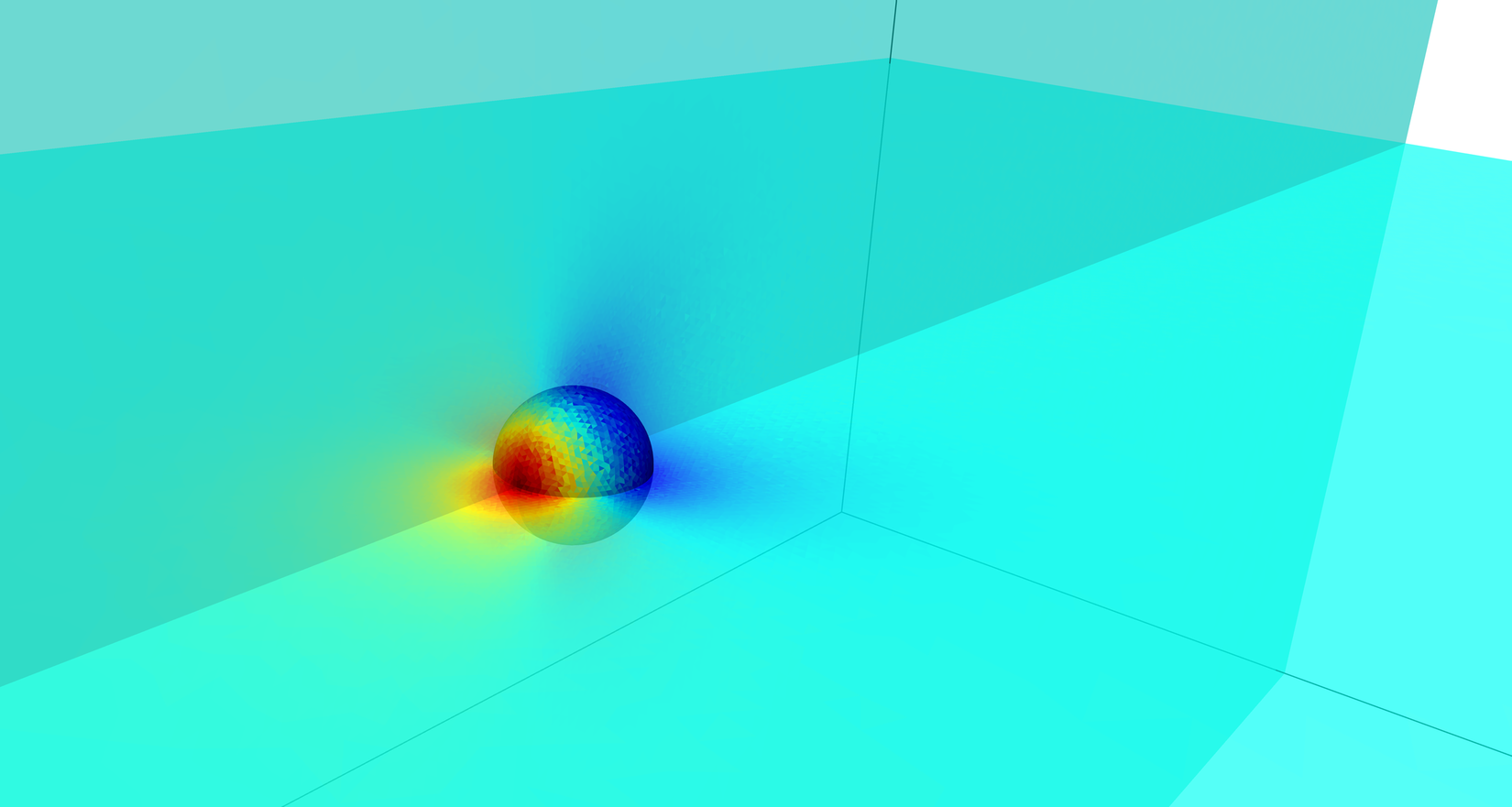}}
	\caption{Flow past a sphere at $Re=20$ - (a,b) Mesh. (c,d) Velocity magnitude and streamlines. (e,f) Pressure. }
	\label{fig:sphereNS}
\end{figure}

The case of $Re=20$ is considered since this benchmark is well documented in the literature~\cite{Schlichting-Gersten-16,Itakura-TB-98,Crivellini-CDB-13} and the symmetry of the flow with respect to the $x_1x_2$ and $x_1x_3$ planes is preserved for low values of the Reynolds number.
Exploiting this symmetry, only one fourth of the domain is simulated. A uniform horizontal velocity $\bu_D = [1,0,0]^T$ is set on the inlet plane $\Ga[D]$ at $x_1=-10D$, a homogeneous Neumann condition is imposed on the outlet plane $\Ga[N]$ at $x_1=20D$, whereas symmetry conditions are enforced on the remaining lateral surfaces $\Ga[S]$. 
The viscosity is selected as $\nu=1/Re$, using $D$ as characteristic length and the magnitude of $\bu_D$ as characteristic velocity of the problem.
It is worth recalling that the convective effects in the Navier-Stokes flow are responsible for breaking the symmetry of the velocity and pressure fields with respect to the $x_2x_3$ plane. In particular, to accurately capture the flow features, an appropriate local refinement of the mesh is required in the wake of the cylinder, as displayed in figure~\ref{fig:meshNS} for the fourth mesh.
In addition,  Figure~\ref{fig:sphereNS} reports the magnitude of the velocity, the streamlines and the pressure field computed using the hybrid pressure formulation with HLL Riemann solver on the fourth level of mesh refinement described below.

To perform a mesh convergence study of the drag coefficient,  five unstructured meshes of tetrahedral cells are constructed using \texttt{Gmsh}~\cite{Gmsh}. 
Table~\ref{tab:meshSphereNS} reports the specifics of the mesh construction, including number of cells, $\numel$, number of faces, $\numfa=\numfa^{\mathcal{I}}+\numfa^{\mathcal{D}}+\numfa^{\mathcal{C}}$ (subdivided in the subset $\mathcal{I}$, $\mathcal{D}$, and $\mathcal{C}$), as well as the local cell size in the first layer near the wall, in the wake of the cylinder, and in the farfield. 
\begin{table}[!htb]
	\centering
	\begin{tabular}{|c|c|c|c|c|c|c|c|c|}
		\hline
		\multirow{2}{*}{Mesh} & $\numel$ & \multicolumn{4}{c|}{Faces} & First wall & Wake & Farfield \\
		\cline{3-6}
		& & $\numfa^{\mathcal{I}}$ & $\numfa^{\mathcal{D}}$ & $\numfa^{\mathcal{C}}$ & $\numfa$ & cell height & cell size & cell size \\
		\hline
		1 & 16,108       & 30,305         & 230      & 3,592    &  34,127       & 0.32 & 1.2 & 1.2 \\
		2 & 32,550      & 62,136          & 280     & 5,648     & 68,064       &0.16 & 0.8 & 1.2 \\
		3 & 110,567     & 214,870       & 516      & 12,012    & 227,398     &0.08 & 0.4 & 1.2 \\
		4 & 473,710    & 931,437       & 1,446  & 30,520  & 963,403     &0.04 & 0.2 & 1.2 \\
		5 & 2,204,960 & 4,364,550  & 5,100  & 85,636  & 4,455,286   &0.02 & 0.1 & 1.2 \\
		\hline		
	\end{tabular}
	\caption{Flow past a sphere at $Re=20$ - Specifics of the computational meshes. }
	\label{tab:meshSphereNS}
\end{table}

The dimensions of the resulting systems are detailed in Table~\ref{tab:dofSphereNS}. For the problem under analysis the number of unknowns for the hybrid velocity are given by $\numdof{\hu} = 3(\numfa^{\mathcal{I}}+\numfa^{\mathcal{C}})$, whereas the unknowns for hybrid pressure and pressure are $\numdof{\hp}=\numfa$, and $\numdof{p}=\numel$, respectively.
Hence, the number of degrees of freedom is $\numdof{\hu} + \numdof{\hp}$ for the hybrid pressure formulation and $\numdof{\hu} + \numdof{p}$ for the FCFV method.
\begin{table}[!htb]
	\centering
	\begin{tabular}{|c|c||c|c||c|c|}
		\hline
		\multirow{2}{*}{Mesh} & \multirow{2}{*}{$\numdof{\hu}$} & \multicolumn{2}{c||}{Hybrid pressure} & \multicolumn{2}{c|}{FCFV} \\
		\cline{3-6}
		 &  & $\numdof{\hp}$ & $\numdof{}$ & $\numdof{p}$ & $\numdof{}$ \\
		\hline
		1 & 101,691          & 34,127       & 135,818       & 16,108         & 117,799 \\
		2 & 203,352       & 68,064     & 271,416        & 32,550        & 235,902 \\
		3 & 680,646       & 227,398    & 908,044      & 110,567      & 791,213 \\
		4 & 2,885,871    & 963,403   & 3,849,274   & 473,710     & 3,359,581  \\
		5 & 13,350,558 & 4,455,286 & 17,805,844 & 2,204,959 & 15,555,517 \\
		\hline		
	\end{tabular}
	\caption{Flow past a sphere at $Re=20$ - Comparison of the number of degrees of freedom of the hybrid pressure and FCFV formulations.}
	\label{tab:dofSphereNS}
\end{table}

\begin{figure}[!htb]
	\centering
	\subfigure[]{\includegraphics[width=0.32\textwidth]{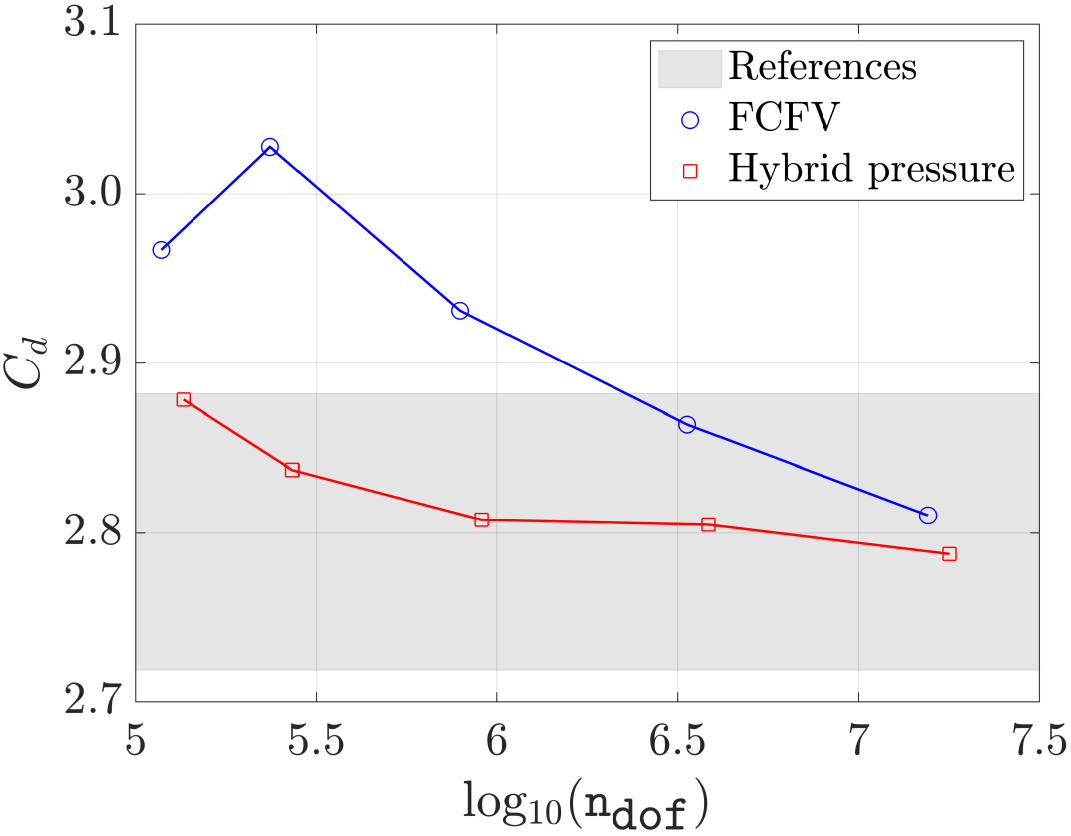}}
	\hspace{2pt}	
	\subfigure[]{\includegraphics[width=0.32\textwidth]{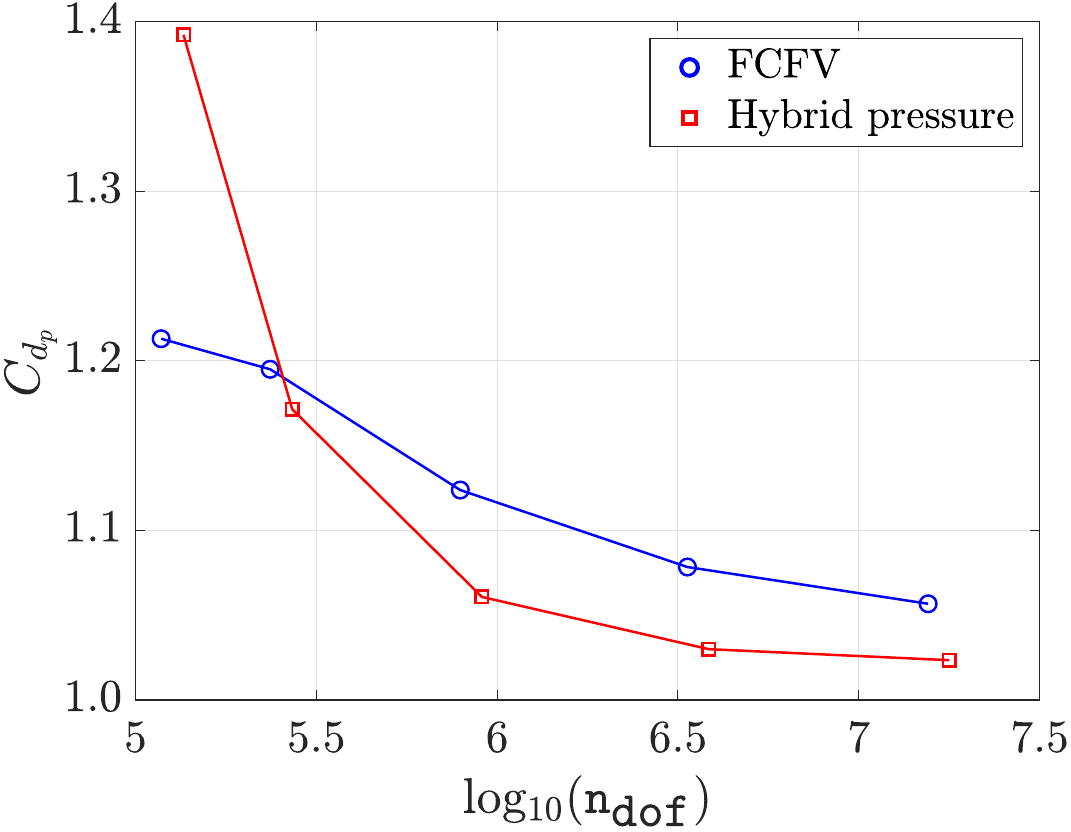}}	
	\hspace{2pt}	
	\subfigure[]{\includegraphics[width=0.32\textwidth]{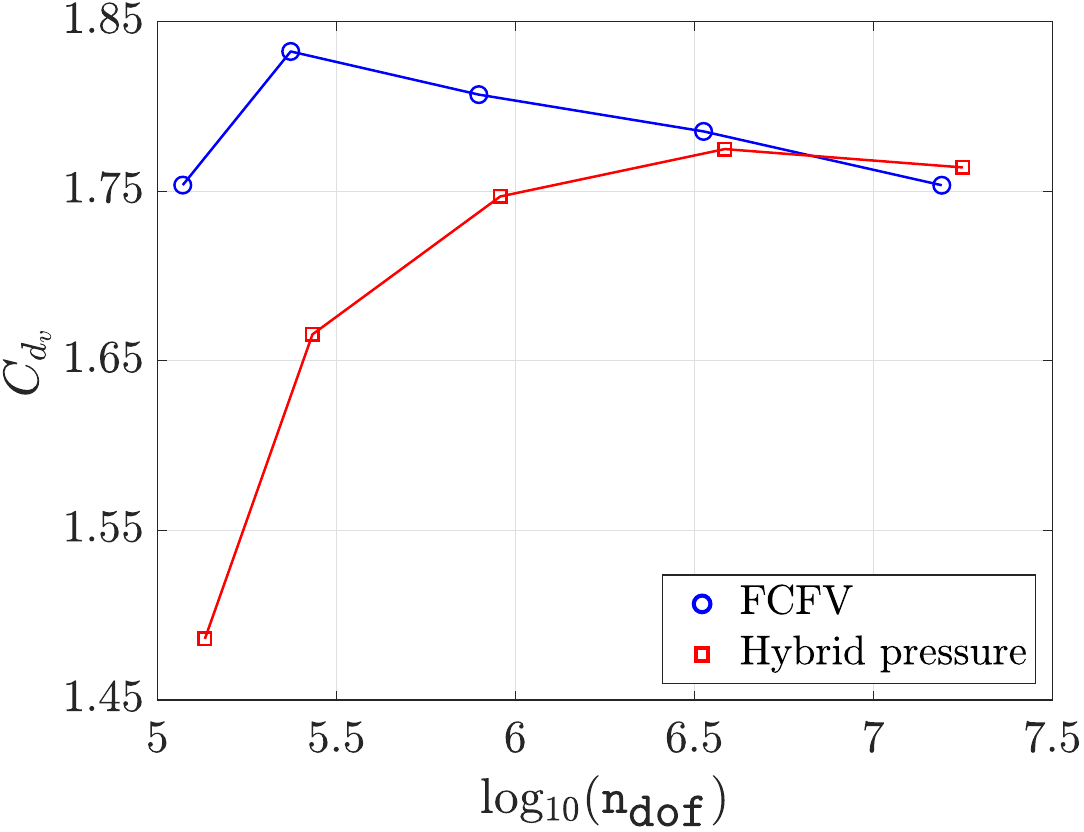}}
	
	\caption{Flow past a sphere at $Re=20$ - Convergence of (a) drag coefficient, (b) pressure part, and (c) viscous part of the drag coefficient as a function of the number of degrees of freedom in the simulation. The shaded area identifies the range of reference values published in the literature for this benchmark.}
	\label{fig:dragSphereNS}
\end{figure}

The convergence of the drag coefficient $C_d$ is displayed in Figure~\ref{fig:dragSphereNS} using the hybrid pressure formulation with HLL stabilisation on the five meshes described above.  The Newton-Raphson algorithm is employed to solve the nonlinear problem, with a tolerance of $10^{-10}$. 
The results show that, even on the coarsest mesh,  the prediction provided by the hybrid pressure formulation lies within the range of values reported in the literature~\cite{Schlichting-Gersten-16,Itakura-TB-98,Crivellini-CDB-13}, whereas the FCFV method requires the fourth level of mesh refinement to achieve a precision in such an interval.  
In addition, although the values of the pressure part, $C_{d_p}$, and the viscous part, $C_{d_v}$,  of the drag coefficient computed using the FCFV method on the first mesh are closer to the corresponding expected values, the hybrid pressure formulation outperforms the FCFV method by providing a faster and almost monotonic convergence of both components.
Table~\ref{tab:dragConvSphereNS} reports the detailed evolution of the drag coefficients and the reference values from the literature, showing excellent agreement of the hybrid pressure formulation.

\begin{table}[!htb]
	\centering
	\begin{tabular}{|c||c|c||c|c||c|c|}
		\hline
		\multirow{2}{*}{Mesh} & \multicolumn{2}{c||}{$C_d$} & \multicolumn{2}{c||}{$C_{d_p}$} & \multicolumn{2}{c|}{$C_{d_v}$} \\
		\cline{2-7}
		 &  Hybrid pressure & FCFV &  Hybrid pressure & FCFV &  Hybrid pressure & FCFV \\
		\hline
		1 & 2.878  & 2.967  & 1.392 & 1.213 & 1.486 & 1.754 \\
		2 & 2.837 & 3.028 & 1.171   & 1.195 & 1.666 & 1.833 \\
		3 & 2.808 & 2.931  & 1.061  & 1.124 & 1.747  & 1.807 \\
		4 & 2.805 & 2.864 & 1.030 & 1.078 & 1.775  & 1.785 \\
		5 & 2.788  & 2.810  & 1.023 & 1.057 & 1.764 & 1.754 \\
		\hline
		\multicolumn{3}{|c|}{Literature} & \multicolumn{4}{c}{} \\
		\cline{1-3}
		\cite{Schlichting-Gersten-16} & \multicolumn{2}{c||}{$2.790$} & \multicolumn{4}{c}{} \\
		\cite{Itakura-TB-98} & \multicolumn{2}{c||}{$2.724$} & \multicolumn{4}{c}{} \\
		\cite{Crivellini-CDB-13} & \multicolumn{2}{c||}{$[2.719,2.882]$} & \multicolumn{4}{c}{} \\
		\cline{1-3}
	\end{tabular}
	\caption{Flow past a sphere at $Re=20$ - Mesh convergence of the drag coefficient $C_d$, its pressure part $C_{d_p}$ and its viscous part $C_{d_v}$.}
	\label{tab:dragConvSphereNS}
\end{table}

To verify the mesh-convergence of the computations, the error on the total drag, as well as its pressure and viscous contribution, is estimated using the last mesh as the reference result. More precisely the estimated error on the $i$-th mesh is measured as
\begin{equation}\label{eq:errQoI}
	\varepsilon_{\!\diamond}^{i} := \frac{|C_{\!\diamond}^{i} - C_{\!\diamond}^{\texttt{5}}|}{|C_{\!\diamond}^{\texttt{5}}|} .
\end{equation}
Table~\ref{tab:dragConvSphereNSerrLast} displays the estimated error for the first four meshes.
\begin{table}[!htb]
	\centering
	\begin{tabular}{|c|c|c||c|c||c|c|}
		\hline
		\multirow{2}{*}{Mesh} & \multicolumn{6}{c|}{Hybrid pressure}	\\
		\cline{2-7}
		& $C_d$ &  $\varepsilon_d$ & $C_{d_p}$ & $\varepsilon_{d_p}$ & $C_{d_v}$ & $\varepsilon_{d_v}$ \\
		\hline
		1 & 2.878 & 0.032 & 1.392 & 0.361 & 1.486 & 0.158 \\
		2 & 2.837 & 0.018 & 1.171 & 0.145 & 1.666 & 0.056 \\
		3 & 2.808 & 0.007 & 1.061 & 0.037 & 1.747 & 0.010 \\
		4 & 2.805 & 0.006 & 1.030 & 0.007 & 1.775 & 0.006 \\	
		5 & 2.788 & -     & 1.023 & -     & 1.764 & -     \\
		\hline		
		\multirow{2}{*}{Mesh} & \multicolumn{6}{c|}{FCFV}	\\
		\cline{2-7}
		& $C_d$ &  $\varepsilon_d$ & $C_{d_p}$ & $\varepsilon_{d_p}$ & $C_{d_v}$ & $\varepsilon_{d_v}$ \\
		\hline
		1 & 2.967 & 0.056 & 1.213 & 0.148 & 1.754 & 0.000 \\
		2 & 3.028 & 0.078 & 1.195 & 0.131 & 1.833 & 0.045 \\
		3 & 2.931 & 0.043 & 1.124 & 0.063 & 1.807 & 0.030 \\
		4 & 2.864 & 0.019 & 1.078 & 0.020 & 1.785 & 0.018 \\
		5 & 2.810 & -     & 1.057 & -     & 1.754 & -     \\		
		\hline
	\end{tabular}
	\caption{Flow past a sphere at $Re=20$ - Drag coefficient and estimated error for the total drag $C_d$, its pressure part $C_{d_p}$ and its viscous part $C_{d_v}$.}
	\label{tab:dragConvSphereNSerrLast}
\end{table}
As qualitatively observed in Figure~\ref{fig:dragSphereNS}, the hybrid pressure formulation achieves convergence on the fourth mesh, assuming a desired error of 1\% for both the pressure and viscous components of the drag.
On the contrary, for the FCFV, the estimated error on the fourth mesh is still above the desired error, showing that further mesh refinement is required and confirming the superiority of the hybrid pressure formulation in approximating pressure in incompressible flows.
It is worth noting that the same conclusions are obtained if the error is measured on the total drag. However, it is preferred to use a more restrictive criteria based on pressure and viscous components, because adding the pressure and viscous component of the drag can lead to error cancellation in the total drag.

Finally, a comparison of the computational cost of the hybrid pressure formulation relative to the FCFV method is provided. From Table\ \ref{tab:dofSphereNS}, it is evident that, for a given mesh, the hybrid pressure formulation has more degrees of freedom than the FCFV approach. 

For instance, consider the fourth level of mesh refinement---the coarsest mesh for which the FCFV method provides a $C_d$ prediction within the range of values reported in the literature. For this mesh (Mesh 4), the global problem in the hybrid pressure formulation involves $17,805,844$ unknowns, compared to $15,555,517$ in the FCFV method. Consequently, for the same mesh, the hybrid pressure formulation is computationally more expensive than the FCFV approach.

Nonetheless, it is more reasonable to compare the methods based on a fixed level of accuracy. For a similar accuracy, the hybrid pressure formulation clearly outperforms the
FCFV approach\ \cite{Vieira-VGSH-24} in terms of computing cost. Specifically, Figure\ \ref{fig:dragSphereNStot} shows that the hybrid pressure formulation can use Mesh 2 to achieve a $C_d$ prediction comparable—if not slightly better—than that provided by the FCFV method on the fourth mesh. Using one computing node with $16$ processors$^\S$, building and assembling the global problem takes $\SI{2.1}{s}$ for the hybrid pressure formulation and $\SI{19.6}{s}$ for the FCFV method. Moreover, as highlighted in Section\ \ref{sc:Computational}, the most computationally expensive task is solving the nonlinear global problem. The hybrid pressure and FCFV methods respectively require $\SI{3.1}{s}$ and $\SI{128.4}{s}$ to solve the global system in each Newton-Raphson iteration.

\let\thefootnote\relax\footnotetext{$^\S$ 2x8-Core Xeon Gold 6134 (3.20 GHz/25MB cache, 2666MHz FSB) with 192 GB RAM}

\section{Concluding remarks}
\label{sc:Conclusion}

This work presented a new hybrid pressure face-centred finite volume solver for viscous laminar incompressible flows. 
The method relies on a mixed hybrid FV formulation in which the cell variables (velocity, pressure and deviatoric strain rate tensor) are expressed in terms of the face velocity and the face pressure by means of a hybridisation procedure. 
In order to do so, incompressibility is enforced weakly and a new definition of the inter-cell mass flux is proposed.
The resulting method has been tested on a suite of 2D and 3D benchmark problems for the steady-state incompressible Navier-Stokes equations. 

First, a linear problem (i.e., Stokes equations) with synthetic solution was used to verify the first-order convergence of velocity, pressure and deviatoric stress tensor on meshes of quadrilateral and triangular cells.  The effect of cell distortion and the sensitity of the method to the choice of the stabilisation coefficient were also assessed.  

Whilst for Stokes flows the hybrid pressure formulation provides results comparable to the FCFV method, its performance are clearly superior in the presence of convection phenomena.
Upon verifying, by means of a synthetic problem,  that optimal convergence of the discrete solution was also achieved in the Navier-Stokes setting, the cavity flow was studied for Reynolds numbers up to $3,200$. These tests showcased the superiority of the hybrid pressure formulation with respect to traditional FCFV schemes, when convective effects are relevant. The hybrid pressure approach with HLL Riemann solver is indeed capable of providing accurate predictions on meshes coarser than the ones required by the FCFV method.
Finally, the flow past a sphere at $Re=20$ was used to confirm the superior performance of the hybrid pressure approach with respect to the FCFV scheme, also in 3D.  Numerical results showed the suitability of the method to compute quantities of engineering interest (i.e., the drag coefficient), using unstructured tetrahedral meshes and achieving sufficient accuracy, in excellent agreement with the reference values published in the literature.

Some computational aspects of the method have also been discussed, numerically studying the properties of the spectrum of the matrix for different cell types,  for uniform and distorted meshes. The experimental observations yielded the conjecture that the matrix obtained from the global problem of the hybrid pressure formulation is negative definite. 
Future work will focus on a more thorough analysis of these properties,  with the goal of developing efficient iterative schemes and preconditioners tailored to the presented hybrid pressure FCFV formulation.

\section*{Acknowledgements}
The authors acknowledge the support of: Generalitat de Catalunya that partially funded the PhD scholarship of DC; H2020 MSCA ITN-EJD ProTechTion (Grant No. 764636) that partially funded the PhD scholarship of LMV; Spanish Ministry of Science, Innovation and Universities and Spanish State Research Agency MICIU/AEI/10.13039/501100011033 (Grants No. PID2020-113463RB-C33 to MG, PID2020-113463RB-C32 to AH and CEX2018-000797-S to MG and AH); Generalitat de Catalunya (Grant No. 2021-SGR-01049 to MG and AH); MG is Fellow of the Serra H\'unter Programme of the Generalitat de Catalunya.

\bibliographystyle{elsarticle-num}
\bibliography{FCFV} 

\appendix

\section{Implementation of the pointwise symmetry of the mixed variable}
\label{app:Voigt}

The symmetry of the mixed variable and the deviatoric stress tensor is imposed pointwise using Voigt notation~\cite{Fish-Belytschko:2007,Tutorial-GSH:2020}. 
More precisely, to treat the $\nsd \times \nsd$ matrix representing $\bLV$,  only $\msd:=\nsd(\nsd+1)/2$ non-redundant components are stored in the vector
\begin{equation} \label{eq:strainVoigt}
\bLVv :=\begin{cases}
\bigl[\Lcomp{11} ,\; \Lcomp{22} ,\; \Lcomp{12} \bigr]^T
&\text{in 2D,} \\
\bigl[\Lcomp{11} ,\; \Lcomp{22} ,\; \Lcomp{33} ,\; \Lcomp{12} ,\; \Lcomp{13} ,\; \Lcomp{23} \bigr]^T
&\text{in 3D.} 
\end{cases}
\end{equation}

Similarly, to represent the normal direction to a surface, the matrix $\bNv$ is defined as
\begin{equation} \label{eq:normalVoigt}
\bNv :=\begin{cases}
\begin{bmatrix}
\ncomp{1} & 0 & \ncomp{2} \\
0 & \ncomp{2} & \ncomp{1}
\end{bmatrix}^T
&\text{in 2D,} \\
\begin{bmatrix}
\ncomp{1} & 0 & 0 & \ncomp{2} & \ncomp{3} & 0\\
0 & \ncomp{2} & 0 & \ncomp{1} & 0 & \ncomp{3} \\
0 & 0 & \ncomp{3} & 0 & \ncomp{1} & \ncomp{2}
\end{bmatrix}^T
&\text{in 3D,} 
\end{cases}
\end{equation}
with $\ncomp{k}$ being the $k$-th component of the discrete outward unit normal vector $\bnV$.

It follows that the second-order tensor $\bhuV \bnV^T  {+}  \bnV  \bhuV^T {-} \frac{2}{3} (\bnV^T\bhuV)\Insd$ associated with the computation of the mixed variable can be rewritten at the discrete level using Voigt notation as $\bDv \bNv \vect{\hat{u}}$,  where the matrix $\bDv$ is given by
\begin{equation}\label{eq:DVoigt}
	\bDv := \begin{bmatrix} \displaystyle	 2\Id{\nsd}  - \frac{2}{3} \Jd{\nsd} & \mat{0}_{\nsd \times \nrr}  \\ \mat{0}_{\nrr \times \nsd} & \Id{\nrr} \end{bmatrix} ,
\end{equation}
$\Id{m}$ is the $m \times m$ identity matrix, $\Jd{\ell}$ the $\ell \times \ell$ matrix with all components equal to 1, and $\nrr := \nsd(\nsd - 1)/2$ denotes the number of rigid body rotations as a function of the number of spatial dimensions $\nsd$.

Hence, the vectorised form of the first equation in~\eqref{eq:localRes} yields the residual
\begin{subequations}\label{eq:resVoigt}
\begin{equation} \label{eq:localResEpsV}
		\bR_{e,L} =   
		\volE \bLVve 
		+ \sum_{j \in \Dset}\! \areaFj \bDv \bNvj \buDj
		+ \sum_{j \in \Bset}\!  \areaFj \bDv \bNvj \bhuj
		,
\end{equation}
for the computation of the mixed variable,  whereas the residual of the global problem associated with the inter-cell continuity of the momentum flux (see first equation in~\eqref{eq:globalRes}) is given by
\begin{equation} \label{eq:globalResUhatV}
		\bR_{e,i,\hu} =  
		\areaFi \Bigl\{ 
		\Bigl( \nu \bNvi^T\bLVve + \btau_i \bue - \btau_i \bhui \Bigr)\!  \chi_{\Iset}(i) 
		+ \bhBi\chi_{\Cset}(i) \Bigr\} ,
		\qquad \text{for all $i \in \Bset$} .
\end{equation}
\end{subequations}
The remaining residuals $\bR_{e,u}$, $\bR_{e,p}$, and $\bR_{e,i,\hp}$ maintain the definitions in equation~\eqref{eq:discreteRes}. Finally,  the discrete boundary operator~\eqref{eq:FCFV_BC} is rewritten using Voigt notation as
\begin{align} \label{eq:FCFV_BCv}
\bhBi =
\begin{cases}
	\hpi\bni + \nu \bNvi^T\bLVve + \btau_i^d(\bue-\bhui) + \bgi
	& \text{ for  $i \in \Nset$,} \\
	\left\{\begin{aligned} & \btki^T \bigl[ \nu\bNvi^T\bLVve + \btau_i^d(\bue-\bhui) \bigr] \\ & \bni^T\bhui \end{aligned} \right\}
	& \text{ for  $i \in \Sset$, and for $k=1,\ldots,\nsd-1$.}
\end{cases}
\end{align}

\section{Hybrid pressure FCFV method for the Stokes equations}
\label{app:Stokes}

For highly viscous flows, the convection term in equation~\eqref{eq:NS} is negligible and the resulting Stokes equations can be written as
\begin{equation} \label{eq:Stokes}
 \left\{\begin{aligned}
 -\Div\bsigma &= \bs       &&\text{in $\Omega$,}\\
   \Div\bu &= 0  &&\text{in $\Omega$,} \\
   \bB(\bepsD,\bu,p) &= \bm{0} &&\text{on $\partial\Omega$,}
 \end{aligned}\right.
\end{equation}
with $\nu=1$, while maintainig the definitions of the Cauchy stress tensor and the boundary conditions introduced in Section~\ref{sc:NS}. 

In this case, the local and global problems~\eqref{eq:FCFV_Local_Solve} and~\eqref{eq:FCFV_Global_Solve} respectively reduce to
\begin{subequations}\label{eq:FCFV_Stokes}
\begin{equation}\label{eq:FCFV_Local_Stokes}
\left\{\begin{aligned} 
 -\int_{\Omega_e} \bL \, d\Omega  
 =& 
 \int_{\partial\Omega_e\cap\Ga[D]} \left(\bu_D{\otimes}\bn{+}\bn{\otimes}\bu_D {-} \frac{2}{3} (\bu_D{\cdot}\bn)\Insd \right) d\Gamma \\
 &+ \int_{\partial \Omega_e \setminus \Ga[D] } \left(\bhu   {\otimes}\bn{+}\bn{\otimes}\bhu {-} \frac{2}{3}(\bhu   {\cdot}\bn)\Insd \right) d\Gamma ,
\\[1ex]
\int_{\partial \Omega_e} \btauD \bu \, d\Gamma 
 =& 
 \int_{\Omega_e}  \bs \, d\Omega
 +\int_{\partial \Omega_e\cap\Ga[D]} \btauD \bu_D \, d\Gamma
 +\int_{\partial \Omega_e\setminus\Ga[D]} \btauD \, \bhu \, d\Gamma
 -\int_{\partial \Omega_e} \hp\bn \, d\Gamma ,
\\[1ex]
-\int_{\partial \Omega_e} \tauP p \, d\Gamma
=&
 \int_{\partial \Omega_e\cap\Ga[D]} \bu_D \cdot \bn \, d\Gamma
 +\int_{\partial \Omega_e\setminus\Ga[D]} \bhu \cdot \bn \, d\Gamma
-\int_{\partial \Omega_e} \tauP \hp \, d\Gamma ,
\end{aligned}\right.
\end{equation}
and 
\begin{equation}\label{eq:FCFV_Global_Stokes}
\left\{\begin{aligned} 
  \sum_{e=1}^{\numel} \Bigg\{ 
	\int_{\partial\Omega_e \setminus \partial\Omega} \nu\bL\bn \, d\Gamma
	+\int_{\partial\Omega_e \setminus \partial\Omega} \btauD \bu \, d\Gamma
	- \int_{\partial\Omega_e \setminus \partial\Omega} \btauD \,\bhu \, d\Gamma \Bigg.  \hspace{70pt}& \\
	+ \Bigg. \int_{\partial\Omega_e \cap \partial\Omega\setminus\Ga[D]} \bhB(\bL,\bu,\bhu,\hp) \, d\Gamma   
\Bigg\} &= \bm{0} ,
\\
  \sum_{e=1}^{\numel} \Bigg\{ 
		\int_{\partial\Omega_e} \tauP p \, d\Gamma
		-\int_{\partial\Omega_e} \tauP \hp \, d\Gamma
\Bigg\} &= 0 ,
\end{aligned}\right.
\end{equation}
\end{subequations}
where the trace of the boundary operator is defined in equation~\eqref{eq:hybridBC}.

The discretisation of equation~\eqref{eq:FCFV_Stokes} following the strategy in Section~\ref{sc:Discrete} and exploiting Voigt notation for the mixed variable as described in Appendix~\ref{app:Voigt} yields the linear local problem
\begin{subequations}\label{eq:linearStokes}
\begin{equation} \label{eq:localStokes}
	\left\{
	\begin{aligned}
		-\volE \bLVve &=   
		\sum_{j \in \Dset}\! \areaFj \bDv \bNvj \buDj
		+ \sum_{j \in \Bset}\!  \areaFj \bDv \bNvj \bhuj		
		,	\\
		\sum_{j \in \Aset}\! \areaFj \btauD_j \bue &=  
		 \volE \bse
		+ \sum_{j \in \Dset}\!  \areaFj  \btauD_j  \buDj
		+ \sum_{j \in \Bset}\!  \areaFj  \btauD_j \, \bhuj
		- \sum_{j \in \Aset}\! \areaFj \hpj \bnj
		,	\\
		-\sum_{j \in \Aset}\! \areaFj \tauP_j \pe &=  
		\sum_{j \in \Dset}\! \areaFj \bnj^T\buDj
		+ \sum_{j \in \Bset}\!  \areaFj \bnj^T\bhuj
		- \sum_{j \in \Aset}\! \areaFj \tauP_j \hpj
		,
	\end{aligned}
	\right.
\end{equation}
for $e=1,\ldots,\numel$ and the linear global problem
\begin{equation} \label{eq:globalStokes}
	\left\{
	\begin{aligned}
	\sum_{e=1}^{\numel}   
		\areaFi \Bigl\{ 
		\Bigl( \nu \bNvi^T\bLVve + \btauD_i \bue - \btauD_i \, \bhui \Bigr)\!  \chi_{\Iset}(i) 
		+ \bhBi\chi_{\Cset}(i) \Bigr\} 
		&= \bm{0} ,
		&&\text{for all $i \in \Bset$}
		,	\\
		\sum_{e=1}^{\numel} 
		\areaFi \tauP_i (\pe - \hpi) 
		&= 0 ,
		&&\text{for all $i \in \Aset$}
		,
	\end{aligned}
	\right.
\end{equation}
\end{subequations}
with the discrete boundary operator $\bhBi$ introduced in~\eqref{eq:FCFV_BCv}. 

For the sake of readability, let us introduce the following quantities, exclusively depending on geometric information and user-defined data:
\begin{subequations}\label{eq:dataStokes}
\begin{equation}\label{eq:coeffStokes}
\mat{A}_u^e :=  \sum_{j \in \Aset}\! \areaFj \btauD_j
, \qquad
\mathrm{a}_p^e := \sum_{j \in \Aset}\! \areaFj \tauP_j ,
\end{equation}
\begin{equation}\label{eq:rhsStokes}
\begin{gathered}
\bfV{L}^e :=  \sum_{j \in \Dset}\! \areaFj \bDv \bNvj \buDj
, \qquad
\bfV{u}^e := \volE \bse + \sum_{j \in \Dset}\!  \areaFj  \btauD_j  \buDj
, \\
\fV{p}^e := \sum_{j \in \Dset}\! \areaFj \bnj^T \buDj .
\end{gathered}
\end{equation}
\end{subequations}
Replacing the expressions of the cell unknowns $\bLVv$, $\buV$, and $\bpV$ from equation~\eqref{eq:localStokes} into~\eqref{eq:globalStokes},  it follows the linear system
\begin{equation} \label{eq:linearSysStokes}
\begin{bmatrix}
	\Khuhu & \Khuhp \\
	\Khphu & \Khphp
\end{bmatrix}
\begin{Bmatrix}
	\bhuV \\
	\bhpV 
\end{Bmatrix}
=
\begin{Bmatrix}
	\bfV{\hu} \\
	\bfV{\hp} 
\end{Bmatrix} .
\end{equation}
In order to compare the structure of system~\eqref{eq:linearSysStokes} to equation~38 and~39 of the FCFV formulation~\cite{RS-SGH:2018_FCFV1},  assume $\Ga[S] = \emptyset$. Hence,  the blocks composing the matrix and the right-hand side are obtained by assembling the cell contributions 
\begin{subequations}\label{linearSysBlocks}
\begin{align}
	[\Khuhu]^e_{i,j} :=&  	
	\areaFi \areaFj 
	\left( 
	- \nu \volE^{-1} \bNvi^T  \bDv \bNvj 
	+  \btauD_i [\mat{A}_u^e]^{-1} \,  \btauD_j
	\right)
	- \areaFi \btauD_i \delta_{ij} \\
	&\text{for all $i,j \in \Bset$} ,  \notag\\[1ex]
	[\Khuhp]^e_{i,j} :=& 
	-\areaFi \areaFj 
	\btauD_i [\mat{A}_u^e]^{-1} \bnj
	+ \areaFi \bni \delta_{ij} \chi_{\Nset}(i) \\
	&\text{for all $i \in \Bset$ and $j \in \Aset$} , \notag\\[1ex]
	[\Khphu]^e_{i,j} :=&
	- \areaFi \areaFj \tauP_i \, [\mathrm{a}_p^e]^{-1} \bnj^T \\
	&\text{for all $i \in \Aset$ and $j \in \Bset$} ,  \notag\\[1ex]
	[\Khphp]^e_{i,j} :=&
	\areaFi \areaFj \tauP_i \, [\mathrm{a}_p^e]^{-1} \tauP_j
	- \areaFi \tauP_i \delta_{ij} \\
	&\text{for all $i \in \Aset$ and $j \in \Aset$} ,  \notag\\[1ex]
	[\bfV{\hu}]^e_i :=&  
	\areaFi \left( \nu \volE^{-1} \bNvi^T \bfV{L}^e  - \btauD_i [\mat{A}_u^e]^{-1}  \bfV{u}^e \right) - \areaFi \bgi \chi_{\Nset}(i)  \\
	&\text{for all $i \in \Bset$} ,  \notag\\[1ex]
	[\bfV{\hp}]^e_i :=& 
	\areaFi \tauP_i \, [\mathrm{a}_p^e]^{-1} \fV{p}^e , \\
	&\text{for all $i \in \Aset$} , \notag
\end{align}
\end{subequations}
for $e=1,\ldots,\numel$,  with $\delta_{ij}$ being a Kronecher delta. 
Note that while $\Khuhu$ and $\Khphp$ are symmetric,  the global matrix in equation~\eqref{eq:linearSysStokes} is not, contrary to the FCFV formulation~\cite{RS-SGH:2018_FCFV1}. Moreover, following the conjecture introduced in section~\ref{sc:Computational},  the matrix in~\eqref{eq:linearSysStokes} is expected to be negative definite, whereas the matrix of the global problem in the FCFV formulation~\cite{RS-SGH:2018_FCFV1} features a saddle point structure.

Finally, if a pure Dirichlet boundary value problem is considered, constraint~\eqref{eq:constraintDirichlet} needs to be included in the above system.  At the discrete level, this yields the additional equation
\begin{equation}\label{eq:constraintDiscrete}
\sum_{e=1}^{\numel} \volE \pe = 0 ,
\end{equation}
which is accounted for by introducing a Lagrange multiplier $\lambda$. Hence, the resulting system is
\begin{equation} \label{eq:linearSysStokesLagMult}
\begin{bmatrix}
	\Khuhu & \Khuhp & \Klhu^T \\
	\Khphu & \Khphp & \Klhp^T \\
	\Klhu     & \Klhp     & 0
\end{bmatrix}
\begin{Bmatrix}
	\bhuV \\
	\bhpV \\
	\lambda
\end{Bmatrix}
=
\begin{Bmatrix}
	\bfV{\hu} \\
	\bfV{\hp} \\
	\fV{\lambda}
\end{Bmatrix} ,
\end{equation}
with 
\begin{subequations}\label{linearSysBlocksLagMult}
\begin{align}
	[\Klhu]^e_{j} :=&  	
	- \areaFj \volE [\mathrm{a}_p^e]^{-1} \bnj^T
	&&\text{for all $j \in \Bset$} , \\
	[\Klhp]^e_{j} :=& 
	\areaFj \volE [\mathrm{a}_p^e]^{-1} \tauP_j
	&&\text{for all $j \in \Aset$} , \\
	[\fV{\lambda}]^e :=&
	\volE	 [\mathrm{a}_p^e]^{-1} \fV{p}^e .
\end{align}
\end{subequations}

\end{document}